\pdfoutput=1

\RequirePackage{fix-cm}

\documentclass[leqno, fleqn, centertags, 12pt]{article}

\usepackage{latexsym}
\usepackage{amsmath}
\usepackage{amssymb}
\usepackage{verbatim}

\usepackage[T1]{fontenc}

\usepackage[upright, widespace]{fourier}

\usepackage{bm}

\usepackage{fullpage}

\usepackage{tikz-cd}

\tikzcdset{row sep/special/.initial=12ex}
\tikzcdset{row sep/spe/.initial=10.5ex}
\tikzcdset{row sep/speh/.initial=13ex}

\tikzcdset{column sep/spec/.initial=5em}
\tikzcdset{column sep/spech/.initial=6em}
\tikzcdset{column sep/speci/.initial=6em}
\tikzcdset{column sep/specc/.initial=8em}
\tikzcdset{column sep/speccc/.initial=10em}

\tikzcdset{row sep/red/.initial=3.79em}

\tikzcdset{row sep/tri/.initial=5em}
\tikzcdset{column sep/tri/.initial=2.4em}

\tikzcdset{row sep/boom/.initial=5em}
\tikzcdset{column sep/boom/.initial=5.1em}

\tikzcdset{row sep/booom/.initial=5.5em}
\tikzcdset{column sep/booom/.initial=5.712em}

\tikzcdset{row sep/boomm/.initial=5.6em}
\tikzcdset{column sep/boomm/.initial=4.25em}

\tikzcdset{row sep/bboom/.initial=8em}
\tikzcdset{column sep/bboom/.initial=8em}

\tikzcdset{row sep/booms/.initial=4.25em}
\tikzcdset{column sep/booms/.initial=4.335em}

\tikzcdset{row sep/boomsss/.initial=4.45em}
\tikzcdset{column sep/boomsss/.initial=4.54em}

\tikzcdset{row sep/oom/.initial=4.5em}
\tikzcdset{column sep/oom/.initial=4.59em}

\tikzcdset{row sep/boomss/.initial=4.75em}
\tikzcdset{column sep/boomss/.initial=4.845em}

\tikzcdset{row sep/free/.initial=12ex}
\tikzcdset{column sep/freee/.initial=4em}

\newskip\nineskipamount \nineskipamount=9pt plus 0pt minus 0pt
\newskip\zeroskipamount \zeroskipamount=0pt plus 0pt minus 0pt
\usepackage[noorphans,vskip=\nineskipamount]{quoting}

\makeatletter
\renewcommand{\@makefntext}[1]{\vspace*{0.5ex}\parindent=0em
\hspace*{-0.4em}
\hbox to 0.4em{\hss\@makefnmark}\hspace*{0.4em}{#1}
}
\makeatother

\newcounter{mysectionnumber}
\newcommand{\mysection}[2]{\setcounter{footnote}{0}
\setcounter{equation}{0}
\setcounter{myparnum}{0}
\refstepcounter{mysectionnumber}
\vspace{27pt}{\Large {\themysectionnumber.} {#1}}\label{#2}\vspace*{15pt}}

\numberwithin{equation}{section}

\newcommand{\myuppar}[1]{\vspace{\medskipamount}\textbf{#1}\hspace*{0.5em}}

\newcommand{\myit}[1]{\textbf{\textit{#1}}\hspace{0.0em}}

\newcounter{myparnum}[mysectionnumber]
\renewcommand{\themyparnum}{\arabic{mysectionnumber}.\arabic{myparnum}}
\newcommand{\mypar}[2]{\refstepcounter{myparnum}{\vspace{\medskipamount}\textbf{{\themyparnum. #1}\label{#2}}\hspace{0.5em}}}

\newcounter{mylemmanum}[myparnum]
\newcounter{myappendnumber}
\newcounter{myaparnum}[myappendnumber]
\newcommand{\myappend}[2]{\setcounter{footnote}{0}
\setcounter{myaparnum}{0}
\setcounter{myparnum}{0}
\refstepcounter{myappendnumber}
\vspace{27pt}{\Large A\dff.{\themyappendnumber.}\oss {#1}}\label{#2}\vspace*{15pt}}

\newcommand{\myapar}[2]{\refstepcounter{myaparnum}{\vspace{\medskipamount}\textbf{{\themyaparnum. #1}\label{#2}}\hspace{0.5em}}}
\renewcommand{\themyaparnum}{A\halfff\fff.\fff\themyappendnumber.\arabic{myaparnum}}

\newcounter{myapparnum}[mysectionnumber]
\newcommand{\proof}{\vspace{\medskipamount}{\textbf{{\emph{Proof}.}}\hspace*{1em}}}

\newcommand{\eproof}{ $\blacksquare$}

\newcommand{\dis}{\displaystyle}

\def\sss{\hspace{0.05em}\ }
\def\dss{\hspace{0.1em}\ }
\def\trs{\hspace{0.15em}\ }
\def\qss{\hspace{0.2em}\ }
\def\pss{\hspace{0.3em}\ }
\def\oss{\hspace{0.4em}\ }

\def\halfff{\hspace*{0.025em}}
\def\fff{\hspace*{0.05em}}
\def\dff{\hspace*{0.1em}}
\def\trf{\hspace*{0.15em}}
\def\qff{\hspace*{0.2em}}
\def\pff{\hspace*{0.3em}}
\def\off{\hspace*{0.4em}}

\def\ttff{{\hspace*{-0.05em}--\hspace*{0.15em}}}

\newcommand{\hnsp}{\hspace*{-0.05em}}
\newcommand{\nsp}{\hspace*{-0.1em}}
\newcommand{\nnsp}{\hspace*{-0.15em}}
\newcommand{\snsp}{\hspace*{-0.175em}}
\newcommand{\dnsp}{\hspace*{-0.2em}}

\renewcommand{\leq}{\leqslant}
\renewcommand{\geq}{\geqslant}

\newcommand{\rrr}{\mathbf{R}}
\newcommand{\nnn}{\mathbf{N}}

\newcommand{\pr}{\operatorname{p{\fff}r}}

\newcommand{\core}{\operatorname{core}\trf}
\newcommand{\sk}{\operatorname{sk}}
\newcommand{\id}{\operatorname{id}}

\newcommand{\num}[1]{|\qff #1 \qff|}

\newcommand{\norm}[1]{\|\qff #1 \qff\|}
\newcommand{\sco}[1]{\langle\dff #1 \dff\rangle}

\newcommand{\ttoo}{\hspace*{0.2em}\longrightarrow\hspace*{0.2em}}

\begin{document}

\setlength{\baselineskip}{12pt plus 0pt minus 0pt}
\setlength{\parskip}{12pt plus 0pt minus 0pt}
\setlength{\abovedisplayskip}{12pt plus 0pt minus 0pt}
\setlength{\belowdisplayskip}{12pt plus 0pt minus 0pt}

\newskip\smallskipamount \smallskipamount=3pt plus 0pt minus 0pt
\newskip\medskipamount   \medskipamount  =6pt plus 0pt minus 0pt
\newskip\bigskipamount   \bigskipamount =12pt plus 0pt minus 0pt

\author{Nikolai\qss V.\qss Ivanov}
\title{Simplicial\qss sets,\oss Postnikov\qss systems,\\ and\qss bounded\qss cohomology}
\date{}

\footnotetext{\hspace*{-0.65em}\copyright\oss 
Nikolai\qss V.\qss Ivanov,\oss 2020.\trs 
Neither\sss the work reported\sss in\dss the present paper\halfff,\qss
nor\dss its preparation were supported\dss by\dss any corporate entity.}

\maketitle

\renewcommand{\baselinestretch}{1}
\selectfont

\vspace*{12ex}

\myit{\hspace*{0em}\large Contents}\vspace*{1ex} \vspace*{\bigskipamount}\\ 
\hbox to 0.8\textwidth{\myit{\phantom{A.}1.}\hspace*{0.5em} Introduction\hfil  2}\hspace*{0.5em} \vspace*{0.25ex}\\
\hbox to 0.8\textwidth{\myit{\phantom{A.}2.}\hspace*{0.5em} Simplicial\sss sets and $\Delta$\dnsp-sets\hfil 5}\hspace*{0.5em} \vspace*{0.25ex}\\
\hbox to 0.8\textwidth{\myit{\phantom{A.}3.}\hspace*{0.5em} Postnikov\dss systems\sss and\sss minimality\hfil 9}\hspace*{0.5em} \vspace*{0.25ex}\\
\hbox to 0.8\textwidth{\myit{\phantom{A.}4.}\hspace*{0.5em} Classifying\sss spaces of\dss categories and\dss groups\hfil 11}\hspace*{0.5em} \vspace*{0.25ex}\\
\hbox to 0.8\textwidth{\myit{\phantom{A.}5.}\hspace*{0.5em} Bundles\sss with\dss Eilenberg--MacLane\dss fibers\hfil 14}\hspace*{0.5em} \vspace*{0.25ex}\\
\hbox to 0.8\textwidth{\myit{\phantom{A.}6.}\hspace*{0.5em} Unraveling\sss simplicial\sss sets\hfil 22}\hspace*{0.5em} \vspace*{0.25ex}\\
\hbox to 0.8\textwidth{\myit{\phantom{A.}7.}\hspace*{0.5em} Isometric\dss isomorphisms\sss in\dss bounded cohomology\hfil 28}\hspace*{0.5em} \vspace*{1ex}\\
\myit{Appendices}\hspace*{0.5em}  \hspace*{0.5em} \vspace*{1ex}\\
\hbox to 0.8\textwidth{\myit{\phantom{A.}A.1.}\hspace*{0.5em} The constructions of\pss Milnor\dss and\qss Segal\hfil 33}\hspace*{0.5em}\vspace*{0.25ex}\\
\hbox to 0.8\textwidth{\myit{\phantom{A.}A.2.}\hspace*{0.5em} Few\dss technical\dss lemmas\hfil 35}\hspace*{0.5em} \vspace*{1ex}\\
\hbox to 0.8\textwidth{\myit{References}\hspace*{0.5em}\hfil 39}\hspace*{0.5em}  \vspace*{0.25ex}

\renewcommand{\baselinestretch}{1}
\selectfont

\newpage
\mysection{Introduction}{introduction}\vspace{-1.35pt}

\myuppar{Bounded cohomology\sss of\trs topological\sss spaces.}
The bounded cohomology\dss groups\sss
$\widehat{H}^{\fff *}(\trf X \trf)$
of\dss a\sss topological\sss space\sss $X$\sss
were\sss introduced\sss by\qss Gromov\qss \cite{gro}.\oss
The definition of\dss  
$\widehat{H}^{\fff *}(\trf X \trf)$
is\dss almost\dss the same as\sss the definition of\trs the singular
cohomology\sss $H^{\fff *}(\trf X\fff,\qff \rrr \trf)$
of\dss $X$\sss 
with\sss real\sss coefficients.\oss
Namely,\oss in order\sss to define\dss $\widehat{H}^{\fff *}(\trf X \trf)$\sss
one needs only\dss to replace
arbitrary\sss singular $n$\dnsp-cochains by\sss singular $n$\dnsp-cochains
which are\sss bounded as real-valued\dss functions on\sss the set\sss of\dss
singular $n$\dnsp-simplices.\oss
The effect\sss of\trs this change\dss is\dss rather dramatic.\oss
It\dss turns out\dss that\sss $\widehat{H}^{\fff *}(\trf X \trf)$
depends only\sss on\sss the fundamental\sss group of\dss $X$\nnsp,\oss
or\halfff,\oss what\dss is\dss the same,\pss
$\widehat{H}^{\fff *}(\trf X \trf)$
does not\sss depend on\sss the higher\sss homotopy\sss groups of\dss $X$\nnsp.\oss
At\dss the same\sss time $\widehat{H}^{\fff *}(\trf X \trf)$
carries an additional\sss structure,\pss a canonical\sss semi-norm,\oss
and\dss this semi-norm\dss is\dss the\qss \emph{raison\dss d'\^{e}tre}\pss of\dss
the\sss theory.\oss\vspace{-1.35pt}

Gromov's\dss exposition\qss \cite{gro}\qss of\trs the bounded cohomology\dss theory\dss
is\dss rather cryptic.\oss
For more\sss than\dss three decades\sss the only\sss available detailed\sss proofs
of\trs the main\sss results of\trs the bounded cohomology\dss theory\sss
were author's\dss proofs\qss \cite{i1},\oss
under a\sss technical\sss assumption\dss removed\sss in\qss \cite{i2}.\oss
Only\dss recently\qss R.\dss Frigerio\dss and\dss M.\dss Moraschini\qss \cite{fm}\qss
reconstructed\trs Gromov's\dss proofs and\dss provided a detailed\sss exposition
of\trs the\sss theory\dss following\dss Gromov's\dss outline.\oss\vspace{-1.35pt}

\myuppar{Simplicial\sss sets and\dss Postnikov\dss systems.}
The modern singular\sss homology and cohomology\dss theory\sss was
created\dss by\qss S.\dss Eilenberg\sss in\sss his paper\qss \cite{e}.\oss
Later on\trs Eilenberg,\oss in collaboration\dss with\trs S.\dss MacLane\dss
and\trs J.\dss Zilber,\oss undertook a detailed analysis of\trs this\sss theory.\oss 
In\sss particular\halfff,\oss Eilenberg\dss and\dss MacLane\pss 
\cite{em1},\qss \cite{em2}\qss studied\dss the influence of\trs the homotopy\dss
groups of\dss spaces on\dss their\sss homology\sss and cohomology\sss groups.\oss
The second\sss paper\qss \cite{em2}\qss in\sss this series relied\sss on\sss
the notion of\qss \emph{complete semi-simplicial\sss complexes}\pss just\dss introduced\dss
by\dss Eilenberg\dss and\dss Zilber\qss \cite{ez}.\oss
Nowadays\sss complete semi-simplicial\sss complexes are known as\qss
\emph{simplicial\sss sets}.\oss\vspace{-1.35pt}

The problem of\dss influence of\dss homotopy\sss groups on\sss
homology\sss and cohomology\sss groups was addressed\dss by\qss 
Eilenberg\dss and\dss MacLane\dss only\dss in\sss fairly\sss special\sss cases.\oss
Their\dss results were subsumed\sss by\sss 
the\sss theory of\qss \emph{natural\sss systems}\pss 
of\pss M.M.\dss Postnikov\qss
\cite{po1},\qss \cite{po2}.\oss
Natural\sss systems were quickly\dss renamed\sss into\qss
\emph{Postnikov\dss systems}.\oss
The\sss theory\sss of\qss Postnikov\dss systems allows,\oss at\dss least\sss
in\sss principle,\oss to determine\sss the homology\sss and cohomology\dss groups
of\dss a space starting\sss with\sss its\sss homotopy\sss groups
and some additional\sss invariants,\oss known as\qss
\emph{Postnikov\dss invariants}.\oss
Simplicial\sss sets served as\sss the natural\dss framework\dss for\trs
Postnikov's\trs theory\qss \cite{po2}.\oss\vspace{-1.35pt}

\myuppar{Simplicial\sss sets and\dss bounded cohomology\dss theory.}
The bounded cohomology\dss groups\sss
$\widehat{H}^{\fff *}(\trf K \trf)$ of\dss a simplicial\sss set\sss $K$\sss
are defined\dss in an obvious manner\halfff.\oss
Namely,\oss the simplicial $n$\dnsp-cochains of\sss $K$ with\sss real\sss coefficients
are\sss the real-valued\sss functions on\sss the set $K_{\dff n}$
of\sss $n$\dnsp-simplices of\dss $K$\nnsp.\oss
A simplicial $n$\dnsp-cochain\dss is\qss \emph{bounded}\pss
if\trs it\dss is\dss bounded as\sss a real-valued\sss function on $K_{\dff n}$\nsp,\oss
and\sss the coboundary of\dss a bounded cochain\dss is\dss obviously\dss bounded.\oss
The bounded cohomology\sss 
$\widehat{H}^{\fff *}(\trf K \trf)$\sss are\dss defined as\sss the
cohomology\sss of\trs the cochain complex of\trs bounded cochains.\oss
Homotopies of\dss simplicial\dss maps\sss lead\dss to cochain\sss homotopies
of\dss complexes of\dss bounded cochains and,\oss
similarly\dss to\sss the case of\dss spaces,\pss
$\widehat{H}^{\fff *}(\trf K \trf)$ 
depends only\sss on\sss the homotopy\dss type of\dss $K$\nnsp.\oss

It\dss is\dss only\dss natural\dss to expect\dss that\dss
the\sss theories of\dss simplicial\sss sets and\qss Postnikov\dss systems
can\sss be adapted\dss to\sss the bounded version of\trs 
the cohomology\dss theory\sss and,\oss
in\dss particular\halfff,\oss 
used\dss to prove\sss that\dss the bounded cohomology\sss groups\sss
$\widehat{H}^{\fff *}(\trf K \trf)$\sss
do not\sss depend on\sss the higher\sss homotopy\sss groups of\dss $K$\nnsp.\oss
The goal\sss of\trs the present\dss paper\dss is\dss to show\sss 
that\dss this\dss is\dss indeed\dss the case and,\oss
moreover\halfff,\oss that\sss the\sss tools provided\dss by\dss the\sss theories\sss of\dss
simplicial\sss sets and\qss Postnikov\sss systems
are nearly\dss ready\dss for using\dss in\dss 
the bounded cohomology\dss theory.\oss\vspace{-0.25pt}

\myuppar{Kan\sss extension\sss property.}
It\dss is\dss well\dss known\sss that\sss in\sss the\sss the\sss theory
of\dss simplicial\sss sets\sss the internal\dss notions of\dss
homotopy,\pss homotopy\dss type,\pss and\dss homotopy\dss groups
are reasonable only\dss for simplicial\sss sets satisfying\sss an additional\sss
condition\sss known as\sss the\dss Kan\sss extension\sss property.\oss
Such simplicial\sss sets are called\qss
\emph{Kan\dss simplicial\sss sets},\pss
or\halfff,\pss more recently,\dff\oss
\emph{fibrant\sss simplicial\sss sets}.\oss
The need\dss to restrict\dss the class of\dss considered simplicial\sss sets\dss
is\dss clear\sss in\sss the bounded cohomology\dss theory.\oss
Indeed,\oss the bounded cohomology\dss groups of\dss simplicial\sss sets arising\sss
from\sss finite simplicial\sss complexes are\sss the same as\sss the usual\sss
real\sss cohomology\dss groups.\oss
It\dss turns out\dss the\dss Kan\dss extensions condition\dss is\dss
exactly\dss what\dss is\dss needed\sss for\sss 
the bounded cohomology\dss theory.\oss\vspace{-0.25pt}

\myuppar{The main\dss theorems.}
The first\dss main\dss theorem\sss of\trs this paper\dss is\dss concerned\dss
with\dss locally\dss trivial\dss bundles of\dss simplicial\sss sets.\oss
The above notwithstanding,\oss the base\dss is\dss allowed\dss to be
an arbitrary\sss simplicial\sss set.\oss
In\dss this\sss theorem $K\dff(\dff \pi\fff,\qff n\trf)$ denotes\sss 
the\trs Eilenberg--MacLane\trs simplicial\sss set,\oss
as\sss in\dss the classical\dss papers of\qss Eilenberg\dss and\dss MacLane.\oss
See\dss Section\qss \ref{bundles}\qss for\dss the details.\oss\vspace{-0.25pt}

\myuppar{Theorem\qss A.}
\emph{Let\qss $E\fff,\pff B$\sss be simplicial\sss sets,\oss
and\dss let\dss 
$p\dff \colon\dff E\qff \ttoo\qff B$\dss
be a\sss locally\dss trivial\dss bundle with\dss 
the\dss Eilenberg--MacLane\trs simplicial\sss set\qss 
$K\dff(\dff \pi\fff,\qff n\trf)$\sss
as\sss the fiber\halfff.\oss
If\qss $n\qff >\qff 1$\nnsp,\oss
then\dss the map induced\dss by $p$\sss in\sss bounded cohomology\dss 
is\dss an\dss isometric\sss isomorphism.\oss}\vspace{6pt}\vspace{-0.25pt}

See\trs Theorem\qss \ref{bundle-isomorphism}.\oss
This\sss result\sss should\dss be considered as a\sss relative property\sss
of\dss $E$\sss with\sss respect\dss to\sss $B$\nnsp.\oss
Since $K\dff(\dff \pi\fff,\qff n\trf)$
has\sss the\dss Kan\dss extension property,\oss 
the\sss latter\dss is\dss present,\oss
albeit\sss implicitly.\vspace{-0.25pt}

The second\dss main\dss theorem\dss is\dss concerned\sss with\dss the\qss
``classifying\sss spaces''\dss
$\mathit{B}\trf \pi$\sss of\dss discrete groups $\pi$\nnsp,\oss
which are actually\dss not\sss spaces,\oss but\sss another\sss incarnation of\dss 
simplicial\sss sets\sss
$K\dff(\dff \pi\fff,\qff 1\trf)$\nnsp.\oss
See\dss Section\qss \ref{classifying-spaces}\qss for\dss a discussion of\dss
classifying\sss spaces.\oss\vspace{-0.25pt}

\myuppar{Theorem\qss B.}
\emph{Let $\pi$ be a discrete group
and\dss $\kappa\qff \subset\qff \pi$ be a normal\sss
amenable subgroup of\dss $\pi$\nnsp.\oss
Let\dss 
$p\dff \colon\dff
\pi\qff \ttoo\qff \pi/\kappa$\dss
be\sss the quotient\dss homomorphism.\oss
Then\dss
$\mathit{B}\fff p\dff \colon\dff
\mathit{B}\trf \pi
\qff \ttoo\qff
\mathit{B}\trf (\dff \pi/\kappa\dff)$\dss
induces\sss isometric\sss isomorphism\sss 
in\sss bounded cohomology.\oss}\vspace{6pt}\vspace{-0.25pt}\vspace{-0.125pt}

See\trs Theorem\qss \ref{quotient-by-amenable}.\oss
Note\sss that\dss the classifying spaces\sss $\mathit{B}\trf \pi$\sss are\dss
Kan\dss simplicial\sss sets.\oss\vspace{-0.25pt}

In addition\sss to\sss the classical\dss theory of\dss simplicial\sss sets,\oss
the proofs of\trs these\sss theorems\sss involve\sss the notion of\qss
\emph{$\Delta$\dnsp-sets},\oss introduced\dss under\dss the name\qss
\emph{semi-simplicial\dss complexes}\pss
by\dss Eilenberg\dss and\dss Zilber\dss in\sss the same paper\qss \cite{ez}\qss
in\sss which\sss they\sss introduced simplicial\sss sets.\oss
The $\Delta$\dnsp-sets differ\sss from simplicial\sss sets by\dss the\sss lack\sss
of\dss degeneracy\sss operators.\oss
See\dss Section\qss \ref{simplicial-sets}\qss for\sss the details.\oss
The modern\dss term\qss ``\dnsp$\Delta$\dnsp-set\fff''\qss goes back\dss to\dss
C.\dss Rourke\dss and\dss B.\dss Sanderson\qss \cite{rs}.\oss
Our\sss main example of\dss a $\Delta$\dnsp-set\dss is\dss the
infinitely\sss dimensional\sss simplex $\Delta\fff [\dff \infty\dff]$\nnsp,\oss
the union over $n\qff \in\pff \nnn$\sss 
of\trs the standard $n$\dnsp-dimensional\sss simplices
$\Delta\fff [\halfff n\dff]$\sss
considered as $\Delta$\dnsp-sets.\oss
Again,\oss see\dss Section\qss \ref{simplicial-sets}\qss for\sss the details.\oss
The $\Delta$\dnsp-sets are used\dss to construct\qss \emph{unravelings}\pss of\dss
simplicial\sss sets.\oss
The idea of\trs the unraveling\sss goes back\dss to\trs G.\dss Segal\pss \cite{s}\qss 
and\qss J.\dss Milnor's\dss construction of\dss classifying\sss spaces\qss \cite{mi}.\oss
The\qss \emph{unraveling}\pss of\dss a simplicial\sss set\sss $K$\sss
is\dss simply\dss the $\Delta$\dnsp-set\sss 
$K\dff \times\dff \Delta\fff [\dff \infty\dff]$\nnsp.\oss
In\dss this product\dss the degeneracy operators\sss of\dss $K$
are ignored and\dss the product\dss is\dss taken\qss 
\emph{dimension-wise},\oss exactly\sss as\sss
the products of\dss simplicial\sss sets.\oss
The main\sss result\sss about\sss unravelings\dss 
is\dss the following.\oss

\myuppar{Theorem\qss C.}
\emph{Let\qss $K$\sss be a simplicial\sss set.\oss
The projection\dss
$p\dff \colon\dff
K\dff \times\dff \Delta\fff [\dff \infty\dff]
\qff \ttoo\qff
K$\dss
induces\sss isometric\dss isomorphisms\dss in\dss the  
bound\-ed cohomology\dss groups.\oss}\vspace{6pt}

See\trs Theorem\qss \ref{blow-up-isomorphism}.\oss
The proof\dss is\dss based on a modification of\trs 
the method of\dss acyclic models.\oss

\myuppar{Applications.}
When combined\dss with\dss the basic facts of\trs the\sss theory\sss of\qss
Postnikov\dss systems,\oss
Theorem\qss A\qss easily\dss implies\sss the following\dss theorem.\oss

\myuppar{Theorem\qss D.}
\emph{Let\qss $K$\sss be\dss a connected\trs Kan\dss
simplicial\sss set\sss
and\dss
$f\dff \colon\dff
K\qff \ttoo\qff
\mathit{B}\qff \pi_{\fff 1}\dff(\trf K\fff,\qff v\trf)$\nnsp,\oss
where $v$\sss is\dss a vertex of\pss $K$\nnsp,\oss
be a simplicial\dss map\sss inducing isomorphism of\trs fundamental\dss groups.\oss
Then\sss $f$\sss induces 
an\sss isomorphism\dss in\dss bounded cohomology.\oss}\vspace{6pt}

See\trs Theorem\qss \ref{bounded-and-pione}.\oss
This\dss theorem\dss together\sss with\trs
Theorem\qss B\qss easily\dss implies\sss the following.\oss

\myuppar{Theorem\qss E.}
\emph{Let\qss $K\fff,\pff L$\dss be connected\dss Kan\qss simplicial\sss sets
and\dss let\sss $v$ be a vertex of\qss $K$\nnsp.\oss
Let\qss $f\dff \colon\dff K\qff \ttoo\qff L$\qss be a simplicial\dss map.\oss
If\qss 
$f_{\dff *}\dff \colon\dff
\pi_{\fff 1}\dff(\trf K\fff,\qff v\trf)
\qff \ttoo\qff
\pi_{\fff 1}\dff(\trf L\fff,\qff f\dff(\dff v\trf)\trf)$\dss
is\dss surjective and\dss has amenable kernel,\oss
then\dss $f$\dss induces an\sss isometric\sss isomorphisms\sss 
in\dss bounded cohomology.\oss}\vspace{6pt}

See\trs Theorem\qss \ref{iso-amenable-kernel}.\oss
When applied\dss to\sss the singular\sss simplicial\sss sets of\trs
topological\sss spaces,\oss
Theorems\qss E\qss and\qss D\qss turn\sss into\trs
Gromov's\pss \emph{Mapping\dss theorem}\qss and\dss its\trs
Corollary\qss (A)\qss respectively.\oss
See\qss \cite{gro},\oss p.\qss 40.\oss
In\dss fact,\oss Gromov\trs deduces\sss his\dss Mapping\dss theorem\dss
from\dss his\dss Corollary\qss (A).\oss
Observing\sss some similarity\dss between\dss this deduction\qss
(see\qss \cite{gro},\oss the\sss top of\qss p.\qss 47)\qss
and\trs Segal's\dss unraveling\dss of\dss categories\qss \cite{s}\qss
was\sss the starting\dss point\sss of\trs the present\dss paper\halfff.\oss

\myuppar{The structure of\trs the\sss paper\halfff.}
Sections\qss \ref{simplicial-sets},\qss \ref{postnikov},\qss
and\dss the first\dss half\dss of\pss Section\qss \ref{classifying-spaces}\qss
are devoted\dss to\sss the basic definitions and\sss a\sss review of\trs theories used\sss
in\dss the paper\halfff.\oss
The second\dss half\dss of\pss Section\qss \ref{classifying-spaces}\qss
introduces some ideas behind\dss the proof\dss of\qss
Theorem\qss B\qss and\dss the definition of\trs unravelings.\oss
Appendix\qss \ref{milnor-segal}\qss provides additional\dss motivation,\oss
but\dss is\dss not\dss used\sss in\sss the main\sss part\sss of\trs the\sss paper\halfff.\oss
Section\qss \ref{bundles}\qss is\dss the\sss technical\dss heart\sss of\trs the paper\halfff,\oss
laying the groundwork\sss for\sss the proof\dss of\qss Theorem\qss A.\oss
Theorem\qss C\qss is\dss proved\sss in\trs Section\qss \ref{unraveling},\oss
which does not\sss depends on\sss the rest\sss of\trs the paper.\oss
Section\qss \ref{main-theorems}\qss is\dss devoted\dss to\sss the proofs of\qss
Theorems\qss A\qss and\qss B\qss and deducing\qss
Theorems\qss D\qss and\qss E\qss from\dss them.\oss
Appendix\qss \ref{lemmas}\qss is\dss devoted\dss to\sss the proofs
of\dss several\dss technical\dss lemmas.\oss

\newpage
\mysection{Simplicial\sss sets\qss and\qss $\Delta$\dnsp-sets}{simplicial-sets}

\myuppar{The categories\dss $\bm{\Delta}$\dss and\dss $\Delta$\nnsp.}
We will\dss include\sss $0$\sss in\sss the set\trs $\nnn$\dss of\dss natural\dss numbers.\oss
For every\dss $n\qff \in\pff \nnn$\dss let\dss
$[\halfff n\dff]$\dss be\sss the set\dss
$\{\trf 0\fff,\pff 1\fff,\pff \ldots\fff,\pff n\trf\}$\nnsp.\oss
The category\dss $\bm{\Delta}$\dss has sets\dss $[\halfff n\dff]$\dss
as objects and\dss non-decreasing\dss maps\dss
$[\halfff m\dff]\qff \ttoo\qff [\halfff n\dff]$\dss
as morphisms from\dss $[\halfff m\dff]$\dss to\dss $[\halfff n\dff]$\nnsp,\oss
with\dss the composition\dss being\dss the composition of\dss maps.\oss
The category\dss $\Delta$\dss is\sss the subcategory\dss of\dss $\bm{\Delta}$\dss
having\dss the same objects and\dss strictly\dss increasing\sss maps\dss
$[\halfff m\dff]\qff \ttoo\qff [\halfff n\dff]$\dss
as morphisms from\dss $[\halfff m\dff]$\dss to\dss $[\halfff n\dff]$\nnsp.\oss

\myuppar{Simplicial\sss sets and $\Delta$\dnsp-sets.}
A\qss \emph{simplicial\sss set}\pss is\dss a contravariant\dss functor\dss
from\dss $\bm{\Delta}$\dss to\sss the category\sss of\dss sets.\oss
Similarly\halfff,\pss a\qss \emph{$\Delta$\dnsp-set}\pss is\dss a contravariant\dss functor\dss
from\dss $\Delta$\dss to\sss the category\sss of\dss sets.\oss
So,\oss a simplicial\sss set\sss $K$\sss consists of\dss a\sss set\sss $K_{\dff n}$\sss
for every\dss $n\qff \in\pff \nnn$\dss and a map\dss
$\theta^{\fff *}\dff \colon\dff K_{\dff n}\qff \ttoo\qff K_{\dff m}$\dss 
for every\dss non-decreasing\sss map\dss
$\theta\dff \colon\dff [\halfff m\dff]\qff \ttoo\qff [\halfff n\dff]$\nnsp.\oss
For $\Delta$\dnsp-set\sss $K$\sss the map\dss $\theta^{\fff *}$ is\dss defined
only\dss if\dss $\theta$\sss is\dss strictly\dss increasing.\oss
If\trs $K$\dss is\dss a\sss simplicial\sss set\halfff,\oss 
then\dss the restriction of\trs the functor\dss $K$\dss to\sss the subcategory\sss $\Delta$\sss
of\dss $\bm{\Delta}$\dss
is\dss a $\Delta$\dnsp-set\halfff,\oss
which\sss we will\sss denote by\dss $\Delta\dff K$\sss
or\sss simply\dss by\sss $K$\nnsp.\oss

The elements of\trs $K_{\dff n}$\sss are called\dss \emph{$n$\dnsp-simplices},\oss
or\qss \emph{simplices\sss of\dss dimension\sss $n$}\dss of\dss $K$\nnsp,\oss
and\dss the maps\dss $\theta^{\fff *}$\sss the\qss 
\emph{structure maps}\qss of\trs $K$\nnsp.\oss
The $0$\dnsp-simplices are also called\qss \emph{vertices}\qss
and\dss if\dss $\sigma\qff \in\pff K_{\dff n}$\nsp,\oss 
then\dss the\qss \emph{vertices}\qss of\sss $\sigma$ are\sss
$0$\dnsp-simplices of\trs the form $\theta^{\fff *}\dff(\dff \sigma\trf)$
with\sss $\theta$\sss being a map\dss
$[\dff 0\dff]\qff \ttoo\qff [\halfff n\dff]$\nnsp.

If\trs $K\fff,\pff L$\dss are either 
simplicial\sss  or\sss $\Delta$\dnsp-sets,\oss
then a\qss \emph{simplicial\dss map}\pss $K\qff \ttoo\qff L$\dss 
is\dss a natural\dss transformation of\dss functors,\oss
i.e.\qss as a sequence of\dss maps\dss
$K_{\dff n}\qff \ttoo\qff L_{\dff n}$\dss such\dss that\vspace{0pt}
\[
\quad
\begin{tikzcd}[column sep=boomss, row sep=boomss]
K_{\dff n}
\arrow[r]
\arrow[d, "\dis \theta^{\fff *}"']
&
L_{\dff n}
\arrow[d, "\dis \theta^{\fff *}"']
\\
K_{\dff m} 
\arrow[r]
&
L_{\dff m}
\end{tikzcd}
\]

\vspace{-10pt}
is\dss a commutative diagram\dss for every\sss morphism\dss
$\theta\dff \colon\dff [\halfff m\dff]\qff \ttoo\qff [\halfff n\dff]$\dss
of\dss $\bm{\Delta}$\sss or\sss $\Delta$\sss respectively.\oss

\myuppar{The face and degeneracy\sss operators.}
For every\dss $n\qff \in\pff \nnn$\dss
and\sss $i\qff \in\qff [\halfff n\dff]$\dss
there\dss is\dss a unique surjective non-decreasing\dss map 
$s\dff(\dff i\trf)\dff \colon\dff
[\halfff n\qff +\qff 1\dff]\qff \ttoo\qff [\halfff n\dff]$\dss
taking\dss the value $i$ twice.\oss
If\trs $n\qff >\qff 0$\nnsp,\oss then
then\dss there is\dss a unique strictly\dss increasing\dss map
$d\dff(\dff i\trf)\dff \colon\dff
[\halfff n\qff -\qff 1\dff]\qff \ttoo\qff [\halfff n\dff]$\dss
not\dss taking\dss the value\dss $i$\nnsp.\oss
The structure maps\dss
$s_{\dff i}
\off =\off
s\dff(\dff i\trf)^{\fff *}$\sss
and\dss
$\partial_{\fff i}
\off =\off
d\dff(\dff i\trf)^{\fff *}$\sss
are\dss called\dss the $i${\dnsp}th\qss \emph{degeneracy}\qss and\qss \emph{face\sss operators}\pss
respectively\halfff.\oss
If\sss $\sigma$\sss is\dss a simplex\sss of\dss $K$\nnsp,\oss
then\dss $\partial_{\dff i}\dff \sigma$\dss is\dss called\dss the\sss $i${\dnsp}th\qss
\emph{face}\qss of\sss $\sigma$\nnsp.

Clearly,\oss every\dss non-decreasing\dss map\dss
$[\halfff m\dff]\qff \ttoo\qff [\halfff n\dff]$\dss
admits a unique presentation as a composition\dss
$[\halfff m\dff]\qff \ttoo\qff [\fff k\dff]\qff \ttoo\qff [\halfff n\dff]$\dss
of\dss a surjective non-decreasing map\dss
$[\halfff m\dff]\qff \ttoo\qff [\fff k\dff]$\dss
and\dss a strictly\dss increasing map\dss
$[\fff k\dff]\qff \ttoo\qff [\halfff n\dff]$\nnsp.\oss
On\dss the other\sss hand,\oss
every\dss strictly\dss increasing\dss map\dss
is\dss a composition of\dss several\dss maps of\trs the form\sss $d\dff(\dff i\trf)$\nnsp,\oss
and every\sss surjective non-decreasing map\dss is\dss a composition of\dss
several\dss maps of\trs the form\sss $s\dff(\dff i\trf)$\nnsp.\oss

It\dss follows\dss that\sss every\sss structure map\dss is\dss
a composition of\dss several\dss face\sss and\sss degeneracy\sss operators.\oss
These operators satisfy\sss some simple and\dss well\dss known\dss relations\sss 
implied\dss by\sss relations between\dss maps\dss 
$s\dff(\dff i\trf)\dff,\off d\dff(\trf j\trf)$\nsp,\oss
which\sss we do\dss not\dss reproduce here.\oss
Conversely\halfff,\oss the face and degeneracy operators\dss 
$\partial_{\fff j}\qff,\off s_{\dff i}$\dss
satisfying\dss these relations can\dss be extended\dss to a contravariant\dss
functor\sss from\dss $\bm{\Delta}$\dss to\sss the category\sss of\dss sets,\oss
i.e.\qss to a simplicial\sss set\halfff.\oss
Similarly,\oss face operators\dss $\partial_{\fff j}$\dss satisfying\dss
the relations involving only\dss face operators can\sss be extended\dss
to a functor\sss from\sss $\Delta$\sss to\sss the category\sss of\dss sets,\oss
i.e.\qss to a $\Delta$\dnsp-set\halfff.\oss
If\trs $K$\dss is\dss a\sss simplicial\sss set\halfff,\oss
then\dss the $\Delta$\dnsp-set\sss $\Delta\dff K$\sss is\dss 
the result\sss of\dss ignoring\dss the degeneracy\sss operators 
of\trs $K$\nnsp.\oss\vspace{1.45pt}

\myuppar{Non-degenerate simplices.}
A simplex $\sigma$ of\dss a simplicial\sss set $K$
is\dss said\dss to be\qss \emph{degenerate}\pss if\trs it\dss
belongs\sss to\sss the image of\dss some $s_{\dff i}$\nsp,\oss
and\qss \emph{non-degenerate}\pss otherwise.\oss
An $n$\dnsp-simplex $\sigma$ 
is\dss degenerate\dss if\trs and\dss only\trs if\trs
$\sigma\off =\off \theta^{\fff *}\dff(\dff \tau\trf)$\sss
for an $m$\dnsp-simplex\ $\tau$ with $m\qff <\qff n$
and a surjective non-decreasing\dss
$\theta\dff \colon\dff [\halfff n\dff]\qff \ttoo\qff [\halfff m\dff]$\nnsp.\oss
By\sss a\dss lemma of\trs Eilenberg\dss and\dss Zilber\qss \cite{ez},\oss
if\dss $\tau$\sss is\dss required\dss to be non-degenerate,\oss
then\sss the presentation\sss
$\sigma\off =\off \theta^{\fff *}\dff(\dff \tau\trf)$\sss
is\dss unique.\oss
See\trs  Lemma\qss \ref{eilenberg-zilber}.\oss\vspace{1.45pt}

\myuppar{Simplicial\sss sets\sss from $\Delta$\dnsp-sets.}
Let\dss $D$\dss be a $\Delta$\dnsp-set\halfff.\oss 
It\dss gives rise\sss to simplicial\sss set
$\bm{\Delta}\dff D$\dss defined as\sss follows.\oss
The\dss $n$\dnsp-simplices of\trs $\bm{\Delta}\dff D$\dss are\sss the pairs\dss
$(\dff \sigma\fff,\pff \rho\trf)$\dss such\dss that\dss
$\sigma$\dss is\dss an $l$\dnsp-simplex of\trs $D$\dss 
for some\dss $l\qff \leq\qff n$\dss and\dss
$\rho\dff \colon\dff
[\halfff n\dff]\qff \ttoo\qff [\dff l\qff]$\dss
is\dss a surjective non-decreasing\dss map.\oss
In order\dss to define $\theta^{\dff *}$\dss
for a non-decreasing\dss map\trs
$\theta\dff \colon\dff [\halfff m\dff]\qff \ttoo\qff [\halfff n\dff]$\nnsp,\oss
we represent\sss $\theta$\sss as\sss the composition\dss
$\theta\off =\off \tau\dff \circ\trf \varphi$\nnsp,\oss
where\dss $\tau$\dss is\dss a\sss strictly\dss increasing\dss map and\sss
$\varphi$\dss is\dss a\sss surjective non-decreasing\dss map,\pss
and\dss set\dss 
$\theta^{\fff *}\dff(\dff \sigma\fff,\pff \rho\trf)
\off =\off
(\trf \tau^{\dff *}\dff(\dff \sigma\trf)\fff,\pff \varphi\trf)$\nnsp.\oss
One can easily\sss check\dss that\dss 
$(\dff \theta\dff \circ\trf \eta\trf)^{\dff *}
\off =\off
\theta^{\dff *}\fff \circ\qff\fff \eta^{\dff *}$\dss
and\dss hence\sss $\bm{\Delta}\fff K$\dss is\dss indeed a simplicial\sss set\halfff.\oss
The correspondence\dss
$D\off \longmapsto\off \bm{\Delta}\dff D$\dss
naturally\sss extends\sss to simplicial\sss maps,\oss
i.e.\qss leads\sss to a functor from\dss the category\sss of\dss $\Delta$\dnsp-sets\sss 
to\sss the category\sss of\dss simplicial\sss sets.\oss\vspace{1.45pt}

\myuppar{Simplicial\sss sets and $\Delta$\dnsp-sets from simplicial\sss complexes.}
Recall\dss that\sss a\qss \emph{simplicial\sss complex}\pss $S$\dss
is\dss a set\sss of\qss \emph{vertices}\qss $V\off =\off V_{\dff S}$\dss
together\sss with a collection of\dss finite subsets of\trs $V$\dnsp,\oss
called\qss \emph{simplices}\pss of\dss $S$\nnsp,\oss 
subject\dss to\sss the condition\dss that\sss 
a subset\sss of\dss a simplex\dss is\dss also a simplex.\oss
Elements of\dss a simplex\dss $\sigma\qff \subset\qff V$\dss
are called\dss the\qss \emph{vertices of}\dss $\sigma$\nnsp.\oss
A simplex $\tau$\sss is\dss said\dss to be a\qss \emph{face}\qss
of\dss a simplex $\sigma$\sss if\trs $\tau\qff \subset\qff \sigma$\nnsp.\oss
A\qss \emph{local\sss order}\qss on a simplicial\sss complex\sss $S$\sss
is\dss an assignment\sss of\dss a\sss 
linear order\dss $<_{\dff \sigma}$\trs
on $\sigma$ for each simplex $\sigma$\nnsp.\oss
These orders are required\dss to agree in\dss the sense\sss that\sss
$<_{\dff \tau}$\dss is\dss the restriction of\dss $<_{\dff \sigma}$\sss
if\dss $\tau$\sss is\dss a\sss face of\dss $\sigma$\nnsp.\oss
For example,\oss if\dss $<$\dss is\dss a\sss linear order on\sss $V_{\dff S}$\nsp,\oss
then\dss the restrictions of\dss $<$\dss to simplices form a\sss local\sss order on\sss $S$\nnsp.\oss
The simplest\sss examples of\dss simplicial\sss sets and $\Delta$\dnsp-sets\sss 
are provided\dss by\dss the following construction.\oss\vspace{1.45pt}

A\sss locally\sss ordered simplicial\sss complex\sss $S$\sss gives rise\sss
to a $\Delta$\dnsp-set\sss $\Delta\fff S$\sss and\dss a
simplicial\dss set $\bm{\Delta}\fff S$\nnsp.\oss
The $n$\dnsp-simplices of\sss $\Delta\fff S$\sss
are injective maps\dss 
$\sigma\dff \colon\dff
[\halfff n\dff]\qff \ttoo\qff V$\dss
such\dss that\dss the image\dss $\sigma\dff(\trf [\halfff n\dff]\trf)$\dss
is\dss a\sss simplex\sss of\dss $S$\sss
and $\sigma$ is\dss order-preserving.\oss
Of\dss course,\oss the sets\dss $[\halfff n\dff]$\dss are\sss considered\dss with\dss their\sss
natural\sss order\sss induced\dss from\dss $\nnn$\nnsp.\oss
The $n$\dnsp-simplices of\sss $\bm{\Delta}\fff S$\sss
are maps\dss 
$\sigma\dff \colon\dff
[\halfff n\dff]\qff \ttoo\qff V$\dss
such\dss that\dss $\sigma\dff(\trf [\halfff n\dff]\trf)$\dss
is\dss a\sss simplex\sss of\dss $S$\sss
and $\sigma$\sss in\dss non-decreasing\sss with\sss respect\dss to\sss
the orders on $[\halfff n\dff]$ and\dss this simplex.\oss
In\dss both cases\sss the structure maps $\theta^{\fff *}$\dss 
are defined\dss by\trs
$\theta^{\fff *}\dff(\dff \sigma\trf)
\off =\off 
\sigma\dff \circ\trf \theta$\nnsp.\oss
An\dss easy\sss check\sss shows\sss that\dss
$\bm{\Delta}\fff S\off =\off \bm{\Delta}\fff \Delta\fff S$\nnsp.\oss
The\sss local\sss order\sss involved\sss in\sss this construction almost\dss
never\sss matters and\dss is\dss rarely\dss mentioned.\oss

\myuppar{Basic examples.}
For\dss $n\qff \in\pff \nnn$\dss the set $[\halfff n\dff]$
can\sss be considered as a simplicial\sss complex\sss having $[\halfff n\dff]$
as its set\sss of\dss vertices and all\sss subsets of\sss $[\halfff n\dff]$ as simplices.\oss
The usual\sss order on $\nnn$\sss turns $[\halfff n\dff]$\sss into
a\sss locally\sss ordered simplicial\sss complex.\oss
The $\Delta$\dnsp-set $\Delta\fff [\halfff n\dff]$\sss
has as $k$\dnsp-simplices strictly\dss increasing\sss maps\dss
$[\dff k\trf]\qff \ttoo\qff [\halfff n\dff]$\nnsp,\oss
and\dss the simplicial\sss set $\bm{\Delta}\fff [\halfff n\dff]$\sss
has as $k$\dnsp-simplices\sss non-decreasing\sss maps\dss
$[\dff k\trf]\qff \ttoo\qff [\halfff n\dff]$\nnsp.\oss
Clearly,\pss 
$\bm{\Delta}\fff [\halfff n\dff]
\off =\off 
\bm{\Delta}\dff \Delta\fff [\halfff n\dff]$\nnsp.\oss
We will\dss need also\sss the simplicial\sss complex\sss 
$[\dff \infty\dff]$\sss having\sss $\nnn$\sss as its set\sss of\dss vertices
and all\dss finite subsets of\dss $\nnn$\sss as simplices,\oss
as also\sss  the $\Delta$\dnsp-set\sss $\Delta\fff [\dff \infty\dff]$\sss
and\dss the simplicial\sss set\sss
$\bm{\Delta}\fff [\dff \infty\dff]
\off =\off
\bm{\Delta}\dff \Delta\fff [\dff \infty\dff]$\nnsp.\oss

Every\sss non-decreasing\dss map\dss
$\theta\dff \colon\dff [\halfff m\dff]\qff \ttoo\qff [\halfff n\dff]$\dss
defines a simplicial\dss map\dss
$\theta_{\dff *}\dff \colon\dff
\bm{\Delta}\fff [\halfff m\dff]
\qff \ttoo\qff
\bm{\Delta}\fff [\halfff n\dff]$\dss
by\dss the rule\dss
$\theta_{\dff *}\dff(\dff \sigma\trf)
\off =\off
\theta\dff \circ\dff \sigma$\nnsp.\oss
Clearly,\pss 
$(\trf \theta\dff \circ\dff \eta\qff)_{\dff *}
\off =\off
\theta_{\dff *}\dff \circ\qff \eta_{\dff *}$\dss
and\dss the assignments 
$n\off \longmapsto\off \bm{\Delta}\fff [\halfff n\dff]$\dss
and\dss
$\theta\off \longmapsto\off \theta_{\dff *}$\dss
define a covariant\dss functor from $\bm{\Delta}$
to\sss the category\sss of\dss simplicial\sss sets.\oss 
Similarly,\oss strictly\dss increasing maps\dss
$\theta\dff \colon\dff [\halfff m\dff]\qff \ttoo\qff [\halfff n\dff]$\dss
define simplicial\dss maps\dss
$\theta_{\dff *}\dff \colon\dff
\Delta\fff [\halfff m\dff]
\qff \ttoo\qff
\Delta\fff [\halfff n\dff]$\dss
and\dss lead\dss to a functor\sss from\sss $\Delta$\sss to $\Delta$\dnsp-sets.\oss

\myuppar{Kan\dss extension and\dss lifting\dss properties.}
Let $n\qff \in\pff \nnn$ and\sss $k\qff \in\qff [\halfff n\dff]$\nnsp.\oss
The\dss \emph{$k$\dnsp-horn}\qss of\dss $[\halfff n\dff]$\sss
is\dss the simplicial\sss complex\sss $[\halfff n\dff]_{\fff k}$
having $[\halfff n\dff]$ as\sss the set\sss of\dss vertices and
subsets of\sss $[\halfff n\dff]$ not\sss containing\sss
$[\halfff n\dff]\qff \smallsetminus\qff \{\trf k\qff\}$
as simplices.\oss
Equivalently\halfff,\pss $[\halfff n\dff]_{\dff k}$ is\dss obtained\dss from
$[\halfff n\dff]$\dss by\dss removing\dss simplices $[\halfff n\dff]$
and $[\halfff n\dff]\qff \smallsetminus\qff \{\trf k\qff\}$\nnsp.\oss
The $k$\dnsp-horn $[\halfff n\dff]_{\fff k}$\dss of\dss $[\halfff n\dff]$\sss 
leads\dss to\sss
the $k$\dnsp-horn\sss
$\bm{\Lambda}_{\fff k}\fff[\halfff n\dff]
\off =\off
\bm{\Delta}\fff [\halfff n\dff]_{\fff k}$\dss
of\dss $\bm{\Delta}\fff [\halfff n\dff]$\nnsp.\oss

A simplicial\sss set\dss $K$\sss 
is\dss said\dss to have\sss the\qss
\emph{Kan\dss extension\dss property},\pss 
or\dss to be a\qss \emph{Kan\dss simplicial\sss set\halfff},\oss
if\dss every\sss simplicial\dss map\dss
$\bm{\Lambda}_{\fff k}\fff[\halfff n\dff]
\qff \ttoo\qff K$\dss
can\sss be extended\dss to a simplicial\dss map\dss
$\bm{\Delta}\dff [\halfff n\dff]\qff \ttoo\qff K$\nnsp.\oss
Let\dss $E\fff,\pff B$\dss be simplicial\sss sets.\oss 
A simplicial\sss map\dss 
$p\dff \colon\dff E\qff \ttoo\qff B$\dss 
is\dss said\dss to\sss have\qss
\emph{Kan\trs lifting\dss property},\pss 
or\dss to be a\qss \emph{Kan\dss fibration}\pss
if\dss every commutative diagram of\dss solid\dss arrows
of\trs the form\vspace{-0.125pt}
\[
\quad
\begin{tikzcd}[column sep=boom, row sep=boom]\dis
\bm{\Lambda}_{\fff k}\fff[\halfff n\dff]
\arrow[r]
\arrow[d, "\dis i\dff"']
&
E
\arrow[d, "\dis p"]
\\
\bm{\Delta}\fff[\halfff n\dff] 
\arrow[r]
\arrow[ru, dashed]
&
B\fff,
\end{tikzcd}
\]

\vspace{-9pt}
where\dss $i$\dss is\dss the inclusion,\oss
can\dss be completed\dss by\sss a dashed arrow\dss to a commutative diagram.\oss
A\sss simplicial\sss set\sss $K$\sss is\dss Kan\trs if\trs and\dss only\trs if\trs
the unique map\dss
$K\qff \ttoo\qff \bm{\Delta}\fff[\dff 0\dff]$\dss is\dss a\dss Kan\dss fibration.\oss

\myuppar{Simplices and simplicial\dss maps.}
Let\dss $\bm{\iota}_{\dff n}$\dss be\sss the identity\dss map\dss
$[\halfff n\dff]\qff \ttoo\qff [\halfff n\dff]$\dss
considered as an $n$\dnsp-simplex of\dss
$\bm{\Delta}\fff [\halfff n\dff]$\nnsp.\oss 
Then every $m$\dnsp-simplex of\dss
$\bm{\Delta}\fff [\halfff n\dff]$\sss 
is\dss equal\sss to\sss $\theta^{\fff *}\dff(\qff \bm{\iota}_{\dff n}\trf)$\sss
for a unique non-decreasing\sss map\dss
$\theta\dff \colon\dff [\halfff m\dff]\qff \ttoo\qff [\halfff n\dff]$\nnsp.\oss
It\dss follows\sss that\dss for every simplicial\sss set\dss $K$\sss
simplicial\dss maps\dss
$f\dff \colon\dff
\bm{\Delta}\fff [\halfff n\dff]
\qff \ttoo\qff
K$
are uniquely\sss determined\dss by\sss the images
$f\dff(\qff \bm{\iota}_{\dff n}\dff)$\nnsp.\oss
Conversely,\oss if\dss $\sigma\qff \in\pff K_{\dff n}$\nsp,\oss
then\dss there\dss is\dss a unique simplicial\dss map\dss\vspace{1.5pt}
\[
\quad
i_{\trf \sigma}\dff \colon\dff
\bm{\Delta}\fff [\halfff n\dff]
\qff \ttoo\qff
K
\]

\vspace{-10.5pt}
such\dss that\dss
$\sigma
\off =\off 
i_{\trf \sigma}\dff(\qff \bm{\iota}_{\dff n}\dff)$\nnsp.\oss

\myuppar{Skeletons.}
Let\dss $n\qff \in\pff \nnn$\nnsp.\oss
If\trs $D$\dss is\dss a $\Delta$\dnsp-set,\oss
then\sss the $n${\nnsp}th\qss \emph{skeleton}\qss $\sk_{\dff n}\fff D$\dss of\trs $D$\sss
is\dss the $\Delta$\dnsp-subset\sss of\trs $D$\dss consisting of\dss all $k$\dnsp-simplices
with\dss $k\qff \leq\qff n$\nnsp.\oss
If\trs $K$\dss is\dss a simplicial\sss set,\oss
then\dss $\sk_{\dff n}\fff K$\dss consists of\dss all\sss simplices of\trs the form\dss
$\theta^{\dff *}\dff(\dff \sigma\trf)$\nnsp,\oss where $\sigma$\sss is\dss a $k$\dnsp-simplex\sss
for some\dss $k\qff \leq\qff n$\nnsp.\oss 
The\qss \emph{boundary}\qss $\partial{\fff}\bm{\Delta}\fff [\halfff n\dff]$
of\sss $\bm{\Delta}\fff [\halfff n\dff]$
is\dss defined as\sss the skeleton
$\sk_{\dff n\dff -\dff 1}\fff \bm{\Delta}\fff [\halfff n\dff]$\nnsp.\oss
The simplicial\sss set $\partial{\fff}\bm{\Delta}\fff [\halfff n\dff]$
has as simplices non-de\-creas\-ing\dss maps\
$[\halfff m\dff]\qff \ttoo\qff [\halfff n\dff]$
with\dss the image\pss $\neq\off [\halfff n\dff]$\nnsp.\oss

\myuppar{Products.}
The\qss \emph{product}\pss $K\dff \times\dff L$ of\dss simplicial\sss sets and $\Delta$\dnsp-sets\sss
$K\fff,\pff L$\sss
is\dss defined dimension-wise.\oss
In\sss more details,\pss
$(\trf K\dff \times\dff L\trf)_{\dff n}
\off =\off\dff
K_{\dff n}\dff \times\trf L_{\dff n}$\sss
and\dss the structure maps of\trs $K\dff \times\dff L$\dss
are\sss the products of\trs the structure maps of\trs $K$\dss and\dss $L$\nnsp.\oss
This\qss \emph{dimension-wise}\qss product\dss is\dss hardly\dss natural\dss
for\sss finite $\Delta$\dnsp-sets,\oss 
but\dss products\sss with\sss $\Delta\fff [\dff \infty\dff]$
play a\sss key\sss role in\sss our\sss theory.\oss

\myuppar{Homotopies.}
The maps\sss 
$d\dff(\dff 0\dff),\pff d\dff(\dff 1\dff)\qff \colon\qff
[\dff 0\dff]\qff \ttoo\qff [\dff 1\dff]$\sss
take $0$ to $1$ and $0$ respectively.\oss
Let\sss
$i\dff(\dff 0\dff)\off =\off d\dff(\dff 1\dff)_{\dff *}$\sss
and\sss
$i\dff(\dff 1\dff)\off =\off d\dff(\dff 0\dff)_{\dff *}$\sss
be\sss the corresponding\sss maps\sss
$\bm{\Delta}\fff[\dff 0\dff]
\qff \ttoo\qff
\bm{\Delta}\fff[\dff 1\dff]$\nnsp.\oss
Suppose\sss that\sss
$K\fff,\pff L$\sss are simplicial\sss sets and\dss
$f,\pff g\qff \colon\trf
K\qff \ttoo\qff L$\dss
are simplicial\sss maps.\oss
A\qss \emph{homotopy}\pss between $f$ and\sss $g$\sss
is\dss a simplicial\dss map\sss
$h\dff \colon\dff
K\dff \times\dff \bm{\Delta}\fff[\dff 1\dff]
\qff \ttoo\qff 
L$\sss
such\dss that\vspace{3pt}
\[
\quad
f
\off =\off 
k\dff \circ\dff \bigl(\qff \id_{\trf K}\dff \times\dff i\dff(\dff 0\dff)\qff\bigr)
\quad\
\mbox{and}\quad\
g
\off =\off 
k\dff \circ\dff \bigl(\qff \id_{\trf K}\dff \times\dff i\dff(\dff 1\dff)\qff\bigr)
\pff,
\]

\vspace{-9pt}
where\dss $K$\sss is\dss identified\sss with\dss 
$K\dff \times\dff \bm{\Delta}\fff[\dff 0\dff]$\nnsp.\oss
\emph{Homotopy\sss equivalences}\pss are defined\sss in\dss terms of\dss homotopies
in\dss the usual\dss manner\halfff.\oss
A simplicial\sss set\sss $K$\sss is\dss said\dss to be\qss \emph{contractible}\pss
if\trs $K$\sss is\dss homotopy\sss equivalent\dss to $\bm{\Delta}\fff[\dff 0\dff]$\nnsp.\oss
In\dss this case\sss
$\id_{\trf K}$\sss is\dss homotopic\sss to\sss the composition\sss
$K\qff \ttoo\qff \bm{\Delta}\fff[\dff 0\dff]\qff \ttoo\qff K$\sss
of\trs the unique map\sss
$K\qff \ttoo\qff \bm{\Delta}\fff[\dff 0\dff]$\sss
with\sss
$i_{\dff v}\dff \colon\dff
\bm{\Delta}\fff [\halfff n\dff]
\qff \ttoo\qff
K$\sss
for some\sss $v\qff \in\pff K_{\trf 0}$\nsp.\oss
If\dss $K$\sss is\dss contractible,\oss 
then\sss the projection\sss
$K\dff \times\dff L\qff \ttoo\qff L$\sss
is\dss a\sss homotopy\sss equivalence\sss for every\sss $L$\nnsp.\oss

\myuppar{Local\sss systems of\dss coefficients.}
Let\sss $K$\sss be either a simplicial\sss or\sss a $\Delta$\dnsp-set.\oss
Let\sss $\varepsilon$\sss be a $1$\dnsp-simplex of\dss $K$\sss and\sss $v\fff,\pff w$\sss
are vertices of\dss $K$\nnsp.\oss
We say\sss that\sss $\varepsilon$\qss 
\emph{connects $v$ with $w$}\qss if\qss
$\partial_{\dff 1}\dff \varepsilon\off =\off v$\sss
and\dss
$\partial_{\dff 0}\dff \varepsilon\off =\off w$\nnsp.\oss
A\qss \emph{local\sss system of\dss coefficients},\oss or\sss simply\sss a\qss
\emph{local\sss system}\qss $\pi$ on\sss $K$\sss is\dss an assignment\sss
of\dss a\sss group $\pi_{\dff v}$\sss
to every\sss vertex\sss
$v$ of\dss $K$\sss
and\sss an\dss isomorphism\sss 
$\varepsilon^{\dff *}\dff \colon\dff
\pi_{\dff w}\qff \ttoo\qff \pi_{\dff v}$\sss
to every $1$\dnsp-simplex $\varepsilon$ 
connecting $v$ with $w$\nnsp.\oss
These groups and\sss isomorphisms are subject\dss 
to\sss the following\sss condition\halfff:\oss
if\dss $\omega\qff \in\pff K_{\trf 2}$\sss and\dss 
$\rho
\off =\off 
\partial_{\trf 2}\dff \omega\fff,\off\
\sigma
\off =\off 
\partial_{\trf 0}\dff \omega\fff,\off\
\tau
\off =\off 
\partial_{\dff 1}\dff \omega$\nnsp,\oss
then\dss 
$\tau^{\dff *}
\off =\off 
\rho^{\dff *}\dff \circ\qff \sigma^{\dff *}$\dnsp.\oss
A\sss local\sss system of\dss coefficients $\pi$\sss is\dss said\dss to be\qss
\emph{abelian}\qss if\trs all\sss groups $\pi_{\dff v}$\sss are abelian.\oss 

The\qss \emph{leading\sss vertex}\qss of\dss an $n$\dnsp-simplex $\sigma$\sss
is\dss the vertex\sss 
$v_{\dff \sigma}
\off =\off 
\theta^{\dff *}\fff(\dff \sigma\dff)$\nnsp,\oss
where\sss $\theta\dff \colon\dff [\dff 0\dff]\qff \ttoo\qff [\halfff n\dff]$\sss
is\dss the inclusion,\oss
and\dss the\qss \emph{leading\sss edge}\qss of\dss an $n$\dnsp-simplex $\sigma$\sss
is\dss the $1$\dnsp-simplex\sss
$\varepsilon_{\dff \sigma}
\off =\off 
\eta^{\dff *}\fff(\dff \sigma\dff)$\nnsp,\oss
where\sss $\eta\dff \colon\dff [\dff 1\dff]\qff \ttoo\qff [\halfff n\dff]$\sss
is\dss the inclusion.\oss
An $n$\dnsp-cochain of\trs $K$\sss with coefficients in\dss the\sss local\sss
system $\pi$\sss is\dss a map $c$\sss assigning\dss every\sss 
$\sigma\qff \in\pff K_{\dff n}$\sss an element\sss
$c\trf(\dff v\trf)\qff \in\pff \pi_{\dff v}$\nsp,\oss
where\sss $v\off =\off v_{\dff \sigma}$\sss is\dss the\sss leading\sss 
vertex of\dss $\sigma$\nnsp.\oss
The group of\dss such cochains\dss is\dss denoted\dss by\sss
$C^{\dff n}\dff(\trf K\fff,\qff \pi\trf)$\nnsp.\oss
For abelian $\pi$\sss 
the\qss \emph{coboundary\sss operators}\pss
$\partial^{\dff *}\dff \colon\dff
C^{\dff n}\dff(\trf K\fff,\qff \pi\trf)
\qff \ttoo\qff
C^{\dff n\dff +\dff 1}\dff(\trf K\fff,\qff \pi\trf)$\sss
are defined\dss by\dss the formula\vspace{4pt}
\[
\quad
\partial^{\dff *}\dff c\trf(\dff \sigma\dff)
\off =\off\dff
\varepsilon_{\dff \sigma}^{\dff *}\dff
\bigl(\trf c\trf\bigl(\qff \partial_{\trf 0}\dff \sigma\trf\bigr)\qff\bigr)
\off +\off
\sum\nolimits_{\dff i\qff =\qff 1}^{\dff n\dff +\dff 1}\off
(\qff -\qff 1\trf)^{\dff i}\pff
c\trf\bigl(\qff \partial_{\dff i}\dff \sigma\trf\bigr)
\off \in\off
\pi_{\dff v_{\dff \sigma}}
\off.
\]

\vspace{-8pt}
The\qss \emph{cocycles}\pss and\qss \emph{coboundaries}\pss
are defined\sss in\dss terms\sss $\partial^{\dff *}$\sss
in\dss the usual\dss manner\halfff.\oss

\newpage
\mysection{Postnikov\qss systems\qss and\qss minimality}{postnikov}

\myuppar{Comparing\sss simplices.}
Let\dss $K$\sss be a simplicial\sss set\sss
and\dss $q\qff \in\pff \nnn$\nnsp.\oss
Recall\dss that\dss for every $\sigma\qff \in\pff K_{\dff q}$
there exists a unique simplicial\dss map
$i_{\trf \sigma}\dff \colon\dff
\bm{\Delta}\fff [\fff q\dff]\qff \ttoo\qff K$
such\dss that\dss
$i_{\trf \sigma}\dff(\qff \bm{\iota}_{\dff q}\dff)\off =\off \sigma$\nnsp.\oss
Two simplices\dss $\sigma\fff,\pff \tau\qff \in\pff K_{\dff q}$\dss
are said\dss to be\dss \emph{$n$\dnsp-equivalent}\oss if\trs
the restrictions of\trs the maps\dss\vspace{3pt}
\[
\quad
i_{\trf \sigma}\trf,\off i_{\trf \tau}\qff \colon\qff
\bm{\Delta}\fff [\fff q\dff]
\off \ttoo\off
 K
\]

\vspace{-9pt}
to\qss
$\sk_{\dff n}\fff \bm{\Delta}\fff [\fff q\dff]$\dss
are\sss equal,\oss
or\halfff,\oss equivalently\halfff,\oss 
if\pss
$\theta^{\dff *}\dff(\dff \sigma\trf)
\off =\off
\theta^{\dff *}\dff(\dff \tau\trf)$\dss
for every\sss non-decreasing\sss map\dss
$\theta\dff \colon\dff
[\halfff n\dff]\qff \ttoo\qff [\fff q\dff]$\nnsp.\oss
We write\sss 
$\sigma\off \sim_{\dff n}\off \tau$\qss
if\trs $\sigma\fff,\pff \tau$\sss are $n$\dnsp-equivalent.\oss
Obviously,\oss if\dss $q\qff \leq\qff n$\nnsp,\oss
then\dss 
$\sigma\off \sim_{\dff n}\off \tau$\trs
if\trs and\dss only\trs if\trs
$\sigma\off =\off \tau$\nnsp.\oss

Clearly,\pss $\sim_{\dff n}$\dss is\dss an equivalence
relation on\dss the set\sss of\dss simplices.\oss
If\dss 
$\sigma\fff,\pff \tau\qff \in\pff K_{\dff q}$\dss
and\dss
$\sigma\qff \sim_{\dff n}\qff \tau$\nnsp,\oss
then\dss
$\theta^{\dff *}\dff(\dff \sigma\trf)
\qff \sim_{\dff n}\qff
\theta^{\dff *}\dff(\dff \tau\trf)$\dss
for every\sss non-decreasing\dss map\dss
$\theta\dff \colon\dff
[\halfff m\dff]\qff \ttoo\qff [\fff q\dff]$\nnsp.\oss
Therefore\sss the structural\dss maps $\theta^{\dff *}$ of\dss $K$\dss
induce maps between sets of\dss equivalence classes of\dss
$\sim_{\dff n}$\nsp.\oss
These induced\dss maps are\sss the structure maps
of\dss a canonical\sss structure of\dss a simplicial\sss set\sss
on\dss the set\sss $K\dff(\dff n\trf)$\sss 
of\dss equivalence classes with\sss respect\dss to\dss $\sim_{\dff n}$\nsp.\oss 
Clearly,\oss
there\dss is\dss a canonical\sss simplicial\dss map\dss\vspace{3pt}
\[
\quad
p_{\fff n}\dff \colon\dff
K\qff \ttoo\qff K\dff(\dff n\trf)
\pff.
\]

\vspace{-9pt} 
Also,\oss if\dss
$n\qff \leq\qff m$\dss and\dss
$\sigma\qff \sim_{\dff m}\qff \tau$\nnsp,\oss
then\dss
$\sigma\qff \sim_{\dff n}\qff \tau$\nnsp,\oss
and\dss hence\sss
there\dss is\dss a canonical\sss simplicial\dss map\vspace{3pt}
\[
\quad
p_{\fff m,\dff n}\dff \colon\dff
K\dff(\dff m\trf)\qff \ttoo\qff K\dff(\dff n\trf)
\pff.
\]

\vspace{-9pt} 
Clearly,\oss the maps\sss $p_{\fff n}$\sss and\sss $p_{\fff m,\dff n}$\sss
induce isomorphisms of\dss $n${\nnsp}th\sss skeletons.\oss
When\dss there\dss is\dss no danger of\dss confusion,\oss
we denote maps\sss $p_{\fff n}$\sss and\dss $p_{\fff m,\dff n}$\sss
simply\dss by\sss $p$\nnsp.\oss

Let\dss $M$\sss be a simplicial\sss subset\sss of\dss $K$\nnsp.\oss
Clearly,\oss for simplices\dss $\sigma\fff,\pff \tau$\dss of\trs $M$\dss the relation\dss
$\sigma\off \sim_{\dff n}\off \tau$\trs holds in\dss $M$\dss if\trs and\dss only\trs if\trs
it\dss holds in\dss $K$\nnsp.\oss
Therefore\sss there\dss is\dss a canonical\dss injective\sss map\dss
$M\dff(\dff n\trf)\qff \ttoo\qff K\dff(\dff n\trf)$\nnsp.\oss
In such a situation\dss we will\dss identify\sss $M\dff(\dff n\trf)$
with\sss its image in\sss $K\dff(\dff n\trf)$\nnsp.\oss

Two simplices\dss $\sigma\fff,\pff \tau\qff \in\pff K_{\dff q}$\dss
are said\dss to be\qss \emph{homotopic}\pss if\trs the maps\sss $i_{\trf \sigma}$\sss
and\dss $i_{\trf \tau}$\sss
are\sss homotopic
relatively\dss to\dss $\partial\dff \bm{\Delta}\fff [\fff q\dff]$\nnsp.\oss
We write\sss $\sigma\off \sim\off \tau$\sss
if\trs $\sigma\fff,\pff \tau$\sss are homotopic.\oss
If\trs $K$\dss is\dss a\dss Kan\dss simplicial\sss set\halfff,\oss
then\dss $\sim$\dss is\dss an equivalence
relation on\dss the set\sss of\dss simplices.\oss

\myuppar{Postnikov\dss systems.}
The sequence of\dss simplicial\sss sets\sss
$K\dff(\dff 0\dff)\fff,\pff
K\dff(\dff 1\dff)\fff,\pff
\ldots\fff,\pff
K\dff(\dff n\trf)\fff,\pff
\ldots$\sss
together\sss with\dss the maps\dss
$p_{\fff n}$\sss and\sss $p_{\fff m,\dff n}$\sss 
is\dss called\qss \emph{Postnikov}\qss
or\qss \emph{Moore-Postnikov\dss system}\pss of\dss $K$\nnsp.\oss
More precisely,\oss this definition\dss is\trs the
version of\trs the original\sss construction of\qss Postnikov\qss 
\cite{po1},\oss \cite{po2}\qss due\sss to\dss Moore\qss \cite{mo}.\oss
This construction and\dss the notion of\dss homotopic simplices 
are\dss useful\sss only\dss 
for\dss Kan\dss simplicial\sss sets.\oss
If\dss $K$\sss is\dss Kan,\oss
then every\dss term\dss $K\dff(\dff n\trf)$\sss of\trs the\dss
Postnikov\dss system\dss is\dss also\dss Kan,\oss
and all\dss maps\sss
$p_{\fff n}$\sss and\sss $p_{\fff m,\dff n}$\sss 
are\dss Kan\dss fibrations.\oss
See\qss \cite{m},\oss Proposition\qss 8.2.\oss

\myuppar{Minimality.}
Postnikov\dss systems are especially\dss powerful\dss
when\sss $K$\sss is\dss minimal\dss in\dss the following sense.\oss 
A\dss Kan\dss simplicial\sss set\sss $K$\sss is\dss said\dss to be\trs
\emph{minimal}\pss if\dss 
every\dss two homotopic simplices of\trs $K$\dss are equal,\oss    
i.e.\pss that\trs
$\sigma\off \sim\off \tau$\dss 
implies\dss 
$\sigma\off =\off \tau$\dss
for every\dss two simplices\sss 
$\sigma\fff,\pff \tau$\sss of\trs $M$\nnsp.\oss
This notion\dss is\dss going\dss back\sss to\dss
Eilenberg\dss and\trs Zilber\qss \cite{ez}\qss
and\dss Postnikov\qss \cite{po1}.\oss
Every\trs Kan\dss simplicial\sss set\sss $K$\sss contains a\sss
minimal\trs Kan\dss simplicial\sss subset\sss $M$\sss as a strong\sss
deformation\sss retract.\oss
Moreover\halfff,\oss every\dss two such simplicial\sss subsets\sss $M$\sss
are\sss isomorphic.\oss
See\qss \cite{m},\oss Theorems\qss 9.5\qss and\qss 9.8.\oss

Let\sss $E\fff,\pff B$\sss be simplicial\sss sets
and\dss
$p\dff \colon\dff E\qff \ttoo\qff B$\dss 
be\dss a\dss Kan\dss fibration.\oss
Two simplices\dss 
$\sigma\fff,\pff \tau\qff \in\pff E_{\dff q}$\dss
are\qss
\emph{fiberwise\dss homotopic}\oss if\trs
$p\dff(\dff \sigma\trf)
\off =\off 
p\dff(\dff \tau\trf)$\dss 
and\dss there exists\sss a\sss relative\sss to\dss
$\partial\dff \bm{\Delta}\fff [\fff q\dff]$\dss 
homotopy\vspace{1.5pt}
\[
\quad
h\dff \colon\dff
\bm{\Delta}\fff [\fff q\dff]\dff \times\dff \bm{\Delta}\fff [\dff 1\dff]
\qff \ttoo\qff
E
\]

\vspace{-10.5pt}
between\sss the maps\dss $i_{\trf \sigma}$\sss
and\dss $i_{\trf \tau}$\sss
such\dss that\dss the diagram\vspace{-2pt}
\[
\quad
\begin{tikzcd}[column sep=boomsss, row sep=boomsss]\dis
\bm{\Delta}\fff [\fff q\dff]\dff \times\dff \bm{\Delta}\fff [\dff 1\dff]
\arrow[r, "\dis h"]
\arrow[d, "\dis \pr"]
&
{\phantom{,}}E
\arrow[d, "\dis p"]
\\
\bm{\Delta}\fff [\fff q\dff] 
\arrow[r, "\dis i_{\trf \rho}"]
&
{\phantom{,}}B\dff,
\end{tikzcd}
\]

\vspace{-12pt}
where\sss
$\rho
\off =\off
p\dff(\dff \sigma\trf)
\off =\off 
p\dff(\dff \tau\trf)$\sss
and\dss $\pr$\dss is\dss the projection,\oss
is\dss commutative.\oss
In\dss particular\halfff,\oss if\trs $\sigma\fff,\pff \tau$\sss
are fiberwise homotopic,\oss then\sss $\sigma\fff,\pff \tau$\sss
are homotopic.\oss
A\dss Kan\dss fibration\dss
$p\dff \colon\dff E\qff \ttoo\qff B$\dss
is\dss said\dss to be\qss \emph{minimal}\pss
if\trs every\dss two fiberwise homotopic simplices are equal.\oss
For every\trs Kan\dss fibration\dss
$p\dff \colon\dff E\qff \ttoo\qff B$\dss
there exists a simplicial\sss subset\dss $M$\sss of\trs $E$\dss
such\dss that\trs
$p\dff|\trf M\dff \colon\dff E\qff \ttoo\qff B$\dss
is\dss a\sss minimal\trs Kan\dss fibration and\dss $M$\dss
is\dss a\qss \emph{fiberwise\sss strong\sss deformation\dss retract}\pss 
of\dss $E$\dss in\sss a\sss natural\sss sense.\oss
See\qss \cite{gj},\oss Chapter\qss I,\oss Proposition\qss 10.3,\oss
or\qss \cite{m},\oss Theorem\qss 10.9.\oss

\mypar{Theorem.}{locally-trivial}
\emph{Every\dss minimal\trs Kan\trs fibration\sss with connected\sss base\dss
is\dss a\sss locally\trs trivial\dss bundle.\oss}

\proof
See\qss \cite{gz},\oss Section\qss VI.5.4,\oss 
or\pss \cite{m},\oss Theorem\qss 11.11.\oss  \eproof

\myuppar{Postnikov\dss systems and\dss locally\dss trivial\dss bundles.}
If\dss $K$\sss is\dss a\sss minimal\trs Kan\dss simplicial\sss set,\oss
then all\dss maps\sss
$p_{\fff n}$\sss and\sss $p_{\fff m,\dff n}$\sss 
are\sss minimal\trs Kan\dss fibrations.\oss
See\qss \cite{m},\oss Lemma\qss 12.1.\oss
Therefore,\oss Theorem\qss \ref{locally-trivial}\qss implies\sss that\dss these maps
are\sss locally\dss trivial\dss bundles.\oss
In\dss particular\halfff,\vspace{1.5pt}
\[
\quad
p
\off =\off
p_{\fff n,\dff n\dff -\dff 1}\dff \colon\dff
K\dff(\dff n\trf)
\qff \ttoo\qff
K\dff(\dff n\qff -\qff 1\dff)
\]

\vspace{-10.5pt}
is\dss a\sss locally\trs trivial\dss bundle.\oss
Its\sss fiber\sss is\dss the\dss Eilenberg--MacLane\trs simplicial\dss set\dss 
$K\dff(\dff \pi\fff,\qff n\trf)$\nnsp,\oss 
where\sss
$\pi
\off =\off
\pi_{\fff n}\dff(\trf K\fff,\qff v\trf)$\sss
is\dss the $n${\dnsp}th\dss homotopy\dss groups of\dss $K$\nnsp.\oss
See\qss \cite{m},\oss the\sss beginning\sss of\trs Section\qss 25.\oss
For our\sss purposes\sss it\dss is\dss sufficient\dss to know\dss
that\dss for\sss $n\qff >\qff 1$\sss the fiber\dss is\dss an\dss
Eilenberg--MacLane\trs simplicial\dss set\sss 
$K\dff(\dff \pi\fff,\qff n\trf)$\sss
with an abelian\sss group $\pi$\nnsp,\oss
and\dss to identify\dss the simplicial\sss set\sss 
$K\dff(\dff 1\dff)$\nnsp.\oss
See\qss Lemma\qss \ref{f-group}\fff\qss for\sss the\sss latter\halfff.\oss

\newpage
\mysection{Classifying\qss spaces\qss of\pss categories\qss and\qss groups}{classifying-spaces}

\myuppar{Classifying spaces of\dss categories.}
This section\dss is\dss devoted\dss to
some classical\sss constructions of\qss Milnor\qss \cite{mi}\qss
and\trs Segal\qss \cite{s}\qss
(who attributed some of\trs the ideas of\pss \cite{s}\qss to\dss Grothendieck\halfff).\oss

A set\sss $S$\sss with a partial\sss order\sss $\leq$\sss
defines a category\dss having\sss $S$\sss as its set\sss of\dss objects.\oss
For\dss $a\fff,\pff b\qff \in\qff S$\dss
there\sss is\dss exactly\sss one morphisms\dss $a\qff \ttoo\qff b$\sss
if\dss $a\qff \leq\qff b$\nnsp,\oss an none otherwise.\oss
In\dss particular\halfff,\oss sets\sss $[\halfff n\dff]$\sss together\sss
with\dss their natural\sss order can\sss be considered as categories.\oss
From\sss this point\sss of\dss view\sss non-decreasing\sss maps\dss
$\theta\dff \colon\dff
[\halfff m\dff]\qff \ttoo\qff [\halfff n\dff]$\dss
are nothing else but\dss functors\dss
$[\halfff m\dff]\qff \ttoo\qff [\halfff n\dff]$\nnsp.\oss

Every\sss small\sss category\sss $\mathcal{C}$ defines
a simplicial\sss set $\mathit{B}\qff \mathcal{C}$\dnsp,\oss
its\qss \emph{nerve}\pss in\sss the sense of\trs G.\dss Segal\qss \cite{s},\oss
often called also\sss the\trs \emph{classifying\dss space}\qss 
of\dss $\mathcal{C}$\dnsp.\oss
The vertices of\dss $\mathit{B}\qff \mathcal{C}$\sss are\sss the objects
of\dss $\mathcal{C}$\nnsp,\oss
and\dss the $n$\dnsp-simplices are functors\dss
$\sigma\dff \colon\dff
[\halfff n\dff]\qff \ttoo\qff \mathcal{C}$\nnsp.\oss
As usual,\oss the structure maps are defined as compositions.\oss
Namely,\oss if\dss
$\theta\dff \colon\dff
[\halfff m\dff]\qff \ttoo\qff [\halfff n\dff]$\sss
is\dss a non-decreasing\sss map,\oss then\dss
$\theta^{\fff *}\dff(\dff \sigma\trf)
\off =\off 
\sigma\dff \circ\trf \theta$\nnsp,\oss
where in\sss the right\dss hand side\sss $\theta$\sss
is\dss considered as a functor.\oss
Clearly,\oss a functor\dss
$[\halfff n\dff]\qff \ttoo\qff \mathcal{C}$\dss
is\dss determined\dss by\dss its values on\sss objects
and on morphisms\dss
$i\qff \ttoo\qff i\qff +\qff 1$\nnsp,\oss
where\dss $i\qff \in\qff [\halfff n\qff -\qff 1\dff]$\nnsp.\oss
Therefore $n$\dnsp-simplices of\dss $\mathit{B}\qff \mathcal{C}$
correspond\dss to sequences of\trs morphisms of\trs the form\vspace{-0.5pt}\vspace{1.125pt}
\begin{equation}
\label{simplex-category}
\quad
\begin{tikzcd}[column sep=large, row sep=boom]\dis
v_{\trf 0}
\arrow[r, "\dis p_{\trf 1}"]
&
v_{\dff 1}
\arrow[r, "\dis p_{\trf 2}"]
&
\ldots
\arrow[r, "\dis p_{\dff n}"]
&
v_{\dff n}\pff,
\end{tikzcd}
\end{equation}

\vspace{-10.5pt}\vspace{1.125pt}
where each $v_{\dff i}$\sss is\dss an object\sss of\dss $\mathcal{C}$\sss
and each\sss $p_{\dff i}$\sss is\dss a morphism\dss
$v_{\dff i\dff -\dff 1}\qff \ttoo\qff v_{\dff i}$\nnsp.\oss
Of\dss course,\oss the objects\sss $v_{\dff i}$\sss are determined\dss by\sss
the morphisms\sss $p_{\dff k}$\sss and\dss hence $n$\dnsp-simplices 
correspond\dss to sequences\sss
$(\trf p_{\dff 1}\dff,\off p_{\trf 2}\dff,\off \ldots\dff,\off p_{\dff n}\trf)$\sss
of\dss morphisms such\dss that\dss the composition\sss
$p_{\dff i\dff +\dff 1}\dff \circ\dff p_{\dff i}$\sss
is\dss defined\sss for each\sss $i$\sss between\sss $1$\sss and\sss $n\qff -\qff 1$\nnsp.\oss
For\dss $0\qff <\qff i\qff <\qff n$\dss
the boundary\sss operators\dss $\partial_{\dff i}$\sss
acts by\dss replacing\sss $v_{\dff i}$\sss
and\sss morphisms\sss $p_{\dff i}\fff,\pff p_{\dff i\dff +\dff 1}$\sss
by\dss the composition\sss $p_{\dff i\dff +\dff 1}\dff \circ\dff p_{\dff i}$\nsp.\oss
The boundary\sss operators\dss $\partial_{\trf 0}$\sss and\dss $\partial_{\dff n}$\sss
act\dss by\sss simply\sss removing\sss $v_{\trf 0}\fff,\pff p_{\dff 1}$\sss
and\sss $p_{\dff n\dff -\dff 1}\fff,\pff v_{\dff n}$\sss respectively.\oss 
The degeneracy\sss operator\sss $s_{\dff i}$\sss acts by\sss
inserting\dss the identity\dss morphism\dss
$v_{\dff i}\qff \ttoo\qff v_{\dff i}$\nsp.\oss
Cf.\qss \cite{gj},\oss Example\qss I.1.4.\oss

Let\sss $\mathcal{C}\fff,\pff \mathcal{D}$\sss
be\sss two categories.\oss
A functor\dss 
$f\dff \colon\dff
\mathcal{C}\qff \ttoo\qff \mathcal{D}$\dss
defines,\oss in an obvious\sss way,\oss a simplicial\dss map\dss
$\mathit{B}\fff f\dff \colon\dff
\mathit{B}\qff \mathcal{C}\qff \ttoo\qff \mathit{B}\qff \mathcal{D}$\dnsp.\oss
Given\dss two functors\dss
$f,\pff g\dff \colon\dff
\mathcal{C}\qff \ttoo\qff \mathcal{D}$\dnsp,\oss
a natural\dss transformation\dss
$f\qff \ttoo\qff g$\sss
defines a homotopy\dss between\sss
$\mathit{B}\fff f$\sss and\sss $\mathit{B}\dff g$\nnsp.\oss
Indeed,\oss a natural\dss transformation\dss
$t\dff \colon\dff f\qff \ttoo\qff g$\sss
can\sss be considered as\sss a\sss functor\dss
$\mathcal{C}\dff \times\dff [\dff 1\dff]
\qff \ttoo\qff
\mathcal{D}$\nnsp,\oss
where\sss $[\dff 1\dff]$\sss is\dss considered as a category.\oss
One can easily\sss see\sss that\dss the operation\dss
$\mathcal{C}\off \longmapsto\off \mathit{B}\qff \mathcal{C}$\dss
commutes with\dss the products.\oss
Since,\oss obviously,\pss
$\mathit{B}\qff [\dff 1\dff]\off =\off \bm{\Delta}\fff [\dff 1\dff]$\nnsp,\oss
the natural\dss transformation\sss $t$\sss defines a simplicial\dss map\dss
$\mathit{B}\dff t\dff \colon\dff
\mathit{B}\qff \mathcal{C}\dff \times\dff \bm{\Delta}\fff [\dff 1\dff] 
\qff \ttoo\qff \mathit{B}\qff \mathcal{D}$\dnsp,\oss
i.e.\qss a\sss homotopy.\oss
We\sss leave\sss to\sss the reader\dss the verification\dss that\dss this\dss
is\dss a homotopy\dss between\sss
$\mathit{B}\fff f$\sss and\sss $\mathit{B}\dff g$\nnsp.\oss

A discrete group $\pi$ can\sss be considered as a category\dss
with a single object\sss and\sss $\pi$\sss being\dss the set\sss
of\dss morphisms from\sss this object\dss to\sss itself\halfff,\oss
with\dss the composition\sss being\dss the group multiplication.\oss
The classifying space\sss $\mathit{B}\trf \pi$\sss is\dss a\dss Kan\dss
simplicial\sss set.\oss
See\qss \cite{gj},\oss Lemma\qss I.3.5.\oss
Comparing\dss the definitions shows\sss that\dss the usual\sss
and\dss the bounded cohomology of\trs the group\sss $\pi$\sss are,\oss in\dss fact,\oss
cohomology of\trs the classifying space\sss $\mathit{B}\trf \pi$\nnsp.\oss

\myuppar{Milnor's\dss classifying\sss spaces.}
Another classical\sss construction of\dss classifying\sss spaces of\dss
groups\dss is\dss due\sss to\trs Milnor\qss \cite{mi}.\oss
While we will\dss not\sss use\sss it\sss directly,\oss
it\sss serves as a motivation\dss for\dss the definitions of\qss
\emph{unravelings}\pss of\dss classifying spaces and simplicial\sss sets\sss below
and\dss in\trs Section\qss \ref{unraveling}.\oss

For a discrete group\sss $\pi$\sss
let\dss $\mathcal{E}\fff \pi$\sss
be\sss the simplicial\sss complex\sss 
having\sss the product\sss $\pi\dff \times\dff \nnn$\sss 
as\sss the set\sss of\dss vertices 
and as simplices finite subsets\dss 
$\sigma\qff \subset\qff\fff \pi\dff \times\dff \nnn$\sss 
such\dss that\dss the projection\sss
$\pi\dff \times\dff \nnn\qff \ttoo\qff \nnn$\sss
is\dss injective on $\sigma$\nnsp.\oss
There\dss is\dss a\sss left\sss action\sss
$\pi\dff \times\dff \mathcal{E}\fff \pi
\qff \ttoo\qff
\mathcal{E}\fff \pi$\sss 
of\trs the group $\pi$ on\sss $\mathcal{E}\fff \pi$\sss
by\dss the rule\vspace{3pt}\vspace{0.5pt}
\[
\quad
h\dff \cdot\trf
(\trf g\fff,\qff k\trf)
\off =\off
(\qff h\fff \cdot\fff g\fff,\qff k\trf)
\pff,
\]

\vspace{-9pt}\vspace{0.5pt}
where\sss $h\qff \in\qff \pi$\sss
and\sss 
$(\qff g\fff,\qff k\trf)
\qff \in\qff
\pi\dff \times\dff \nnn$\sss
is\dss a vertex of\trs 
$\mathcal{E}\fff \pi$\nnsp.\oss
The quotient\dss 
$\mathcal{B}\dff \pi
\off =\off
\pi\fff\backslash\dff \mathcal{E}\fff \pi$\sss
of\trs the simplicial\sss complex\sss $\mathcal{E}\fff \pi$\sss
by\dss this action\dss is\dss a\sss well-defined simplicial\sss complex.\oss
We will\sss call\dss $\mathcal{B}\dff \pi$\dss the\qss
\emph{Milnor'\dss classifying\dss space}\pss of\sss $\pi$\nnsp.\oss
Milnor\dss defined directly\dss the geometric realization 
$\num{\mathcal{B}\dff \pi}$
for arbitrary\dss topological\dss group $\pi$\nnsp.\oss
By\dss this reason\sss his construction\dss is\dss different\sss
from\sss ours one. 

It\dss is\dss convenient\dss to enhance\sss the structure of\qss
Milnor's\trs classifying\sss space\sss to a $\Delta$\dnsp-set.\oss
The natural\sss order on\sss $\nnn$\sss defines\sss local\sss orders on\sss the
simplicial\sss set\sss $\mathcal{E}\fff \pi$
and\dss $\mathcal{B}\dff \pi$\sss
and\dss allows\sss to\sss turn\dss them\dss into $\Delta$\dnsp-sets,\oss
which\sss we will\sss still\sss denote by\sss 
$\mathcal{E}\dff \pi$\sss and\sss $\mathcal{B}\dff \pi$\nnsp.\oss
The\sss local\sss order on\sss $\mathcal{E}\dff \pi$\sss is\dss invariant\sss
under\sss the\sss left\sss action of\sss $\pi$\nnsp,\oss and\dss 
$\mathcal{B}\dff \pi
\off =\off
\pi\fff\backslash\dff \mathcal{E}\fff \pi$\sss
as $\Delta$\dnsp-sets also.\oss

One of\dss advantages of\dss $\mathcal{B}\dff \pi$\sss 
is\dss the existence of\dss
many\sss automorphisms.\oss
Let\sss $C^{\dff 0}\dff(\trf \nnn\fff,\pff \pi\trf)$\sss
be\sss the group of\dss all\dss maps\dss
$\nnn\qff \ttoo\qff \pi$\nnsp.\oss
Let\dss us\sss define a right\sss action of\dss
$C^{\dff 0}\dff(\trf \nnn\fff,\pff \pi\trf)$\sss
on $\mathcal{E}\fff \pi$\sss by\dss the rule\vspace{3pt}\vspace{0.5pt}
\[
\quad
(\trf g\fff,\qff k\trf)
\off \longmapsto\off
\left(\trf g\dff \cdot\dff c\trf(\dff k\trf)\fff,\qff k\trf\right)
\dff,
\]

\vspace{-9pt}\vspace{0.5pt}
where\sss
$c\qff \in\qff C^{\dff 0}\dff(\trf \nnn\fff,\pff \pi\trf)$\sss
and\sss
$(\trf g\fff,\qff k\trf)
\qff \in\qff
\pi\dff \times\dff \nnn$\sss
is\dss a vertex of\trs 
$\mathcal{E}\fff \pi$\nnsp.\oss
Clearly,\oss the right\sss action of\dss 
$C^{\dff 0}\dff(\trf \nnn\fff,\pff \pi\trf)$\sss
on $\mathcal{E}\fff \pi$\sss preserves\sss the\sss order of\dss vertices
and commutes with\sss the\sss left\sss action of\sss $\pi$\nnsp.\oss
Hence\sss this action\dss leads\sss to a right\sss action of\dss
$C^{\dff 0}\dff(\trf \nnn\fff,\pff \pi\trf)$
on\sss $\mathcal{B}\dff \pi$\nnsp.\oss
If\dss $\kappa\qff \subset\qff \pi$\sss is\dss a
subgroup of\sss $\pi$\nnsp,\oss
then\sss
$C^{\dff 0}\dff(\trf \nnn\fff,\pff \kappa\trf)$\sss
is\dss a subgroup of\sss
$C^{\dff 0}\dff(\trf \nnn\fff,\pff \pi\trf)$\nnsp,\oss
and\dss if\dss $\kappa$\sss is\dss a\sss normal\sss subgroup,\oss
then\vspace{3pt}\vspace{0.5pt}
\begin{equation}
\label{nice-quotient}
\quad
\mathcal{B}\dff (\dff \pi/\kappa\dff)
\off =\off
\mathcal{B}\dff \pi \left/\dff
C^{\dff 0}\dff(\trf \nnn\fff,\pff \kappa\trf)\right.\dff.
\end{equation}

\vspace{-9pt}\vspace{0.5pt}
This obvious property\sss
is\dss the main\sss reason of\dss our\dss interest\sss in\sss $\mathcal{B}\dff \pi$\nnsp.\oss
This property\sss strongly\sss contracts with\dss the properties of\trs the classifying\sss
spaces\sss $\mathit{B}\trf \pi$\nnsp.\oss
Namely,\oss the classifying\sss space\sss $\mathit{B}\trf (\dff \pi/\kappa\dff)$\sss
is\dss not\sss a quotient\sss of\sss $\mathit{B}\trf \pi$\nnsp,\oss
at\dss least\dss not\sss in any\dss natural\dss way.\oss\vspace{-0.125pt}

\myuppar{Unravelings\sss of\trs classifying spaces of\dss groups.}
The $\Delta$\dnsp-set\sss $\mathcal{B}\dff \pi$\sss is\dss isomorphic\sss to\sss 
$\mathit{B}\trf \pi\dff \times\dff \Delta\dff[\dff \infty\dff]$\nnsp.\oss
See\trs Lemma\qss \ref{milnor-space}.\oss
While we are not\dss going\dss to use\sss this result,\oss
it\dss motivates our\sss interest\dss to\sss the $\Delta$\dnsp-set\sss
$\mathit{B}\trf \pi\dff \times\dff \Delta\dff[\dff \infty\dff]$\nnsp,\oss
which\sss we will\sss call\dss the\qss
\emph{unraveling}\pss of\dss $\mathit{B}\trf \pi$\nnsp.\oss

As we will\sss see in a\sss moment,\oss
the action of\dss 
$C^{\dff 0}\dff(\trf \nnn\fff,\pff \pi\trf)$\sss
can\sss be defined\sss directly\dss for\sss
$\mathit{B}\trf \pi\dff \times\dff \Delta\dff[\dff \infty\dff]$\nnsp.\oss
In\sss fact,\oss it\dss is\dss easier and\dss more useful\dss to define\sss
an action of\dss 
$C^{\dff 0}\dff(\trf \nnn\fff,\pff \pi\trf)$\sss
on\dss the simplicial\sss set\sss
$\mathit{B}\trf \pi\dff \times\dff \bm{\Delta}\dff[\dff \infty\dff]$\sss
first.\oss
The\sss latter\dss is\dss the classifying\sss space of\dss a category.\oss
Indeed,\oss the simplicial\sss set 
$\bm{\Delta}\dff[\dff \infty\dff]$
has\sss non-decreasing\dss maps\sss
$[\halfff n\dff]\qff \ttoo\qff \nnn$\sss
as $n$\dnsp-simplices,\oss
with\dss the usual\sss structure maps.\oss
Let\sss 
$\bm{n}$\sss be\sss the category\dss
having\sss $\nnn$ as\sss its\sss set\sss of\dss objects,\oss 
exactly\sss one morphism\sss
$n\qff \ttoo\qff m$\sss
when\dss $n\qff \leq\qff m$\nnsp,\oss
and no morphisms\sss
$n\qff \ttoo\qff m$\sss
when\dss $n\qff >\qff m$\nnsp.\oss
Clearly,\pss
$\bm{\Delta}\dff[\dff \infty\dff]
\off =\off
\mathit{B}\dff \bm{n}$\nnsp.\oss 
It\dss follows\dss that\vspace{3pt}
\[
\quad
\mathit{B}\trf \pi\dff \times\dff \bm{\Delta}\dff[\dff \infty\dff]
\off =\off
\mathit{B}\trf \pi\dff \times\dff \mathit{B}\dff \bm{n}
\off =\off
\mathit{B}\trf(\trf \pi\dff \times\dff \bm{n}\trf)
\pff.
\]

\vspace{-9pt}
The category\sss $\pi\dff \times\dff \bm{n}$\sss has\sss $\nnn$\sss
as\sss the set\sss of\dss objects.\oss
The set\sss of\dss morphisms\sss
$n\qff \ttoo\qff m$\sss
is\dss a copy of\dss $\pi$\sss if\dss $n\qff \leq\qff m$\nnsp,\oss
and\dss is\dss empty\dss if\dss $n\qff >\qff m$\nnsp.\oss
Given\sss
$c\qff \in\qff C^{\dff 0}\dff(\trf \nnn\fff,\pff \pi\trf)$\nnsp,\oss
let\sss 
$a\dff(\dff c\trf)\dff \colon\dff
\pi\dff \times\dff \bm{n}
\qff \ttoo\qff
\pi\dff \times\dff \bm{n}$\sss
be\sss the functor\sss equal\dss to\sss the identity\sss
on objects and and acting\sss on\sss morphisms\sss
$n\qff \ttoo\qff m$\sss
identified\sss with elements of\trs the group $\pi$\sss by\dss the rule\vspace{3pt}
\[
\quad
g
\off \longmapsto\off
c\trf(\dff n\trf)^{\dff -\dff 1}\dff \cdot\dff
g\dff \cdot\dff
c\trf(\dff m\trf)
\pff.
\]

\vspace{-9pt}
Clearly,\pss $a\dff(\dff c\trf)$ is\dss an automorphism of\sss
$\pi\dff \times\dff \bm{n}$
and even an automorphism over $\bm{n}$\nnsp,\oss
in\dss the sense\sss that\sss
$\pr\dff \circ\qff a\dff(\dff c\trf)\off =\off \pr$\nnsp,\oss
where\sss
$\pr\dff \colon\dff
\pi\dff \times\dff \bm{n}
\qff \ttoo\qff
\bm{n}$\sss
is\dss the projection.\oss
Also,\oss all\dss diagrams\vspace{0pt}
\[
\quad
\begin{tikzcd}[column sep=specc, row sep=boomm]\dis
n
\arrow[r, "\dis g"]
\arrow[d, "\dis c\trf(\dff n\trf)"']
&
n
\arrow[d, "\dis c\trf(\dff m\trf)"]
\\
n 
\arrow[r, "\dis a\dff(\dff c\trf)\fff(\dff g\dff)"]
&
m
\end{tikzcd}
\]

\vspace{-7.5pt}
are commutative and\dss hence\sss morphisms\dss
$c\trf(\dff n\trf)\dff \colon\dff
n\qff \ttoo\qff n$\sss
form a natural\sss transformation\sss from\dss the identity\dss
functor\sss to\sss $a\dff(\dff c\trf)$\nnsp.\oss
It\dss follows\sss that\dss the simplicial\dss map\vspace{3pt}\vspace{0.6pt}
\[
\quad
\mathit{B}\dff a\dff(\dff c\trf)
\dff \colon\dff
\mathit{B}\trf(\trf \pi\dff \times\dff \bm{n}\trf)
\qff \ttoo\qff
\mathit{B}\trf(\trf \pi\dff \times\dff \bm{n}\trf)
\]

\vspace{-9pt}\vspace{0.6pt}
is\dss an automorphism of\dss
$\mathit{B}\trf(\trf \pi\dff \times\dff \bm{n}\trf)$\sss
homotopic\sss to\sss the identity.\oss
The map\sss
$c
\off \longmapsto\off
\mathit{B}\dff a\dff(\dff c\trf)$\sss
defines an action of\dss
$C^{\dff 0}\dff(\trf \nnn\fff,\pff \pi\trf)$
on\sss
$\mathit{B}\trf(\trf \pi\dff \times\dff \bm{n}\trf)
\off =\off
\mathit{B}\trf \pi\dff \times\dff \bm{\Delta}\dff[\dff \infty\dff]$\nnsp.\oss
Since\sss the functors $a\dff(\dff c\trf)$ are automorphisms over $\bm{n}$\nnsp,\oss
the simplicial\dss maps\sss
$\mathit{B}\dff a\dff(\dff c\trf)$\sss
are automorphisms of\dss
$\mathit{B}\trf \pi\dff \times\dff \bm{\Delta}\dff[\dff \infty\dff]$\sss 
over\sss
$\bm{\Delta}\dff[\dff \infty\dff]$\nnsp.\oss
It\dss follows\dss that\sss maps\sss
$\mathit{B}\dff a\dff(\dff c\trf)$\sss
leave\sss the $\Delta$\dnsp-subset\sss
$\mathit{B}\trf \pi\dff \times\dff \Delta\dff[\dff \infty\dff]$\sss
of\dss
$\mathit{B}\trf \pi\dff \times\dff \bm{\Delta}\dff[\dff \infty\dff]$\sss
invariant.\oss
By\dss restricting\dss these maps\sss to\sss
$\mathit{B}\trf \pi\dff \times\dff \Delta\dff[\dff \infty\dff]$\sss
we get\sss an action of\dss
$C^{\dff 0}\dff(\trf \nnn\fff,\pff \pi\trf)$
on\sss
$\mathit{B}\trf \pi\dff \times\dff \Delta\dff[\dff \infty\dff]$\nnsp.\oss
If\dss $\kappa\qff \subset\qff \pi$\sss
is\dss a normal\sss subgroup of\sss $\pi$\nnsp,\oss
then a direct\sss verification shows\sss that\vspace{3pt}\vspace{0.6pt}
\[
\quad
\mathit{B}\trf (\dff \pi/\kappa\dff)\dff \times\dff \Delta\dff[\dff \infty\dff]
\off =\off
\mathit{B}\trf \pi\dff \times\dff \Delta\dff[\dff \infty\dff] 
\left/\dff
C^{\dff 0}\dff(\trf \nnn\fff,\pff \kappa\trf)\right.\dff.
\]

\vspace{-9pt}\vspace{0.6pt}
Of\dss course,\oss this\dss is\dss simply\sss another\sss form of\trs
the property\qss (\ref{nice-quotient}).\oss
We will\dss need a slightly\sss stronger\halfff,\oss
but\sss still\sss obvious,\oss form of\trs this property.\oss
Let\sss $1\qff \in\qff \pi$\sss be\sss the unit\sss of\sss $\pi$\nnsp,\oss
and\dss let\dss 
$C_{\dff 0}\dff(\trf \nnn\fff,\pff \pi\trf)$\sss
be\sss the group of\dss maps\sss
$c\dff \colon\dff
\nnn\qff \ttoo\qff \pi$\sss
such\dss that\sss
$c\trf(\dff n\trf)\off =\off 1$\dss
for almost\sss every\sss $n$\nnsp.\oss
Then\vspace{3pt}\vspace{0.6pt}
\begin{equation}
\label{nice-quotient-b}
\quad
\mathit{B}\trf (\dff \pi/\kappa\dff)\dff \times\dff \Delta\dff[\dff \infty\dff]
\off =\off
\mathit{B}\trf \pi\dff \times\dff \Delta\dff[\dff \infty\dff] 
\left/\dff
C_{\trf 0}\dff(\trf \nnn\fff,\pff \kappa\trf)\right.\dff.
\end{equation}

\vspace{-9pt}\vspace{0.6pt}

\newpage
\mysection{Bundles\qss with\qss Eilenberg--MacLane\qss fibers}{bundles}

\myuppar{Locally\dss trivial\dss bundles.}
Let\dss $p\dff \colon\dff E\qff \ttoo\qff B$\dss be a simplicial\sss map\sss
thought\sss as a bundle,\oss
and\dss let\dss $i\dff \colon\dff A\qff \ttoo\qff B$\dss a simplicial\dss map.\oss
Let\dss $i^{\dff *}\fff E\qff \subset\pff E\dff \times\dff A$\dss
be\sss the simplicial\sss subset\sss of\trs $E\dff \times\dff A$\dss
having as $n$\dnsp-simplices pairs\sss 
$(\dff \sigma\fff,\pff \tau\trf)$\sss
such\dss that\sss 
$\sigma\qff \in\pff E_{\dff n}$\nsp,\qss
$\tau\qff \in\qff A_{\dff n}$\nsp,\oss
and\dss 
$p\dff(\dff \sigma\trf)\off =\off i\dff(\dff \tau\trf)$\nnsp.\oss
In other\dss terms,\oss\vspace{3pt}
\[
\quad
(\trf i^{\dff *}\fff E\trf)_{\dff n}
\off =\off
\left\{\pff 
(\dff \sigma\fff,\pff \tau\trf)
\qff \in\pff
E_{\dff n}\dff \times\dff A_{\dff n}
\off \bigl|\off
p\dff(\dff \sigma\trf)
\off =\off 
i\dff(\dff \tau\trf)
\pff\right\}
\pff.
\]

\vspace{-9pt}
The restriction\dss
$i^{\dff *} p\dff \colon\dff
i^{\dff *}\fff E\qff \ttoo\qff A$\dss
of\trs the projection\dss
$E\dff \times\dff A\qff \ttoo\qff A$\dss
to\sss $i^{\dff *}\fff E$\sss
is\dss called\dss the\qss  
\emph{pull-back}\pss of\dss $p$\sss by\dss $i$\nnsp,\oss
or\sss the\qss \emph{bundle\sss induced}\pss from\dss the bundle\dss 
$p\dff \colon\dff E\qff \ttoo\qff B$\trs 
by\dss $i$\nnsp.\oss
The bundle\sss $i^{\dff *} p$\sss has\sss the usual\sss universal\dss
properties of\dss pull-backs.\oss
A simplicial\dss map\dss 
$p\dff \colon\dff E\qff \ttoo\qff B$\dss 
is\dss said\dss to\sss be a\qss \emph{trivial\dss bundle}\pss
with\dss the\sss fiber $F$\trs
if\trs there exists a commutative diagram\vspace{1.5pt}
\[
\quad
\begin{tikzcd}[column sep=booms, row sep=booms]\dis
B\dff \times\dff F
\arrow[r, "\dis t\dff"]
\arrow[d, "\dis \pr\dff"']
&
E
\arrow[d, "\dis p"]
\\
B 
\arrow[r, "\dis ="]
&
B\fff,
\end{tikzcd}
\]

\vspace{-7.5pt}
such\dss that\dss 
$t\dff \colon\dss
B\dff \times\dff F\qff \ttoo\qff E$\dss is\dss an\dss isomorphism.\oss
Such\sss $t$\sss is\dss called a\qss \emph{trivialization}\pss of\dss $p$\nnsp.\oss
A map\dss 
$p\dff \colon\dff E\qff \ttoo\qff B$\dss 
is\dss a\qss \emph{locally\dss trivial\dss bundle}\pss
with\dss the fiber\sss $F$\trs
if\dss for every\sss simplex $\sigma$ of\trs $B$\dss
the pull-back\dss
$i_{\dff \sigma}^{\dff *}\dff p$\dss is\dss a\sss trivial\dss bundle
with\dss the\sss fiber $F$\dnsp.\oss
In\dss this case\sss $E$\sss is\dss called\dss the\qss \emph{total\sss space}\qss
and\dss $B$\dss the\qss \emph{base}\qss of\trs $p$\nnsp.\oss
Clearly,\oss if\dss $p$\sss is\dss a\sss locally\dss trivial\dss bundle,\oss
then\sss $p$\sss is\dss surjective.

\myuppar{Normalized\sss and\sss non-abelian cochains.}
Let\dss $n\qff \in\pff \nnn$\nnsp,\qss $n\qff >\qff 1$\nnsp,\oss 
and\dss let\sss $\pi$\sss be an abelian\sss group.\oss
A cochain of\dss a simplicial\sss set\dss $K$\sss is\dss said\dss to be\qss
\emph{normalized}\pss if\trs it\dss is\dss equal\dss to $0$ on degenerate simplices.\oss 
We will\sss denote\sss by\sss
$\mathcal{C}^{\dff n}\dff(\dff K\fff,\qff \pi\trf)$\dss
the\sss group of\dss normalized $n$\dnsp-cochains of\dss $K$\dss with\sss
coefficients\sss in\sss $\pi$
and\dss by\sss
$\mathcal{Z}^{\dff n}\dff(\dff K\fff,\qff \pi\trf)$\sss
the subgroup  of\dss normalized cocycles.\oss

\myuppar{Eilenberg--MacLane\trs simplicial\dss sets\dss $K\dff(\dff \pi\fff,\qff n\trf)$\nnsp.}
For every\dss $q\qff \in\pff \nnn$\trs let\dss us consider\sss the groups
$\mathcal{C}^{\dff n}\dff(\dff \bm{\Delta}\fff [\dff q\trf]\fff,\qff \pi\trf)$\sss
and\trs
$\mathcal{Z}^{\dff n}\dff(\dff \bm{\Delta}\fff [\dff q\trf]\fff,\qff \pi\trf)$\nnsp.\oss
Every\dss non-decreasing\sss map\dss
$\theta\dff \colon\dff
[\dff r\dff]\qff \ttoo\qff [\dff q\trf]$\dss
induces a simplicial\dss map\dss
$\theta_{\dff *}\dff \colon\dff
\bm{\Delta}\fff [\dff r\trf]\qff \ttoo\qff \bm{\Delta}\fff [\dff q\trf]$\dnsp,\oss
which,\oss in\dss turn,\oss
induces homomorphisms\vspace{3pt}\vspace{1.125pt}
\[
\quad
\theta^{\dff *}
\qff \colon\dff
\mathcal{C}^{\dff n}\dff(\dff \bm{\Delta}\fff [\dff q\trf]\fff,\qff \pi\trf)
\qff \ttoo\qff
\mathcal{C}^{\dff n}\dff(\dff \bm{\Delta}\fff [\dff r\trf]\fff,\qff \pi\trf)
\quad\hspace{0.0em}\
\mbox{and}\quad\hspace{0.0em}\
\]

\vspace{-36pt}\vspace{1.125pt}
\[
\quad
\theta^{\dff *}
\qff \colon\dff
\mathcal{Z}^{\dff n}\dff(\dff \bm{\Delta}\fff [\dff q\trf]\fff,\qff \pi\trf)
\qff \ttoo\qff
\mathcal{Z}^{\dff n}\dff(\dff \bm{\Delta}\fff [\dff r\trf]\fff,\qff \pi\trf)
\pff.
\]

\vspace{-9pt}\vspace{1.125pt}
\emph{Eilenberg--MacLane\trs simplicial\dss set}\qss 
$K\dff(\dff \pi\fff,\qff n\trf)$\sss
is\dss defined as follows.\oss
Its set\sss of\dss $q$\dnsp-simplices\dss is\vspace{3pt}\vspace{1.125pt}
\[
\quad
K\dff(\dff \pi\fff,\qff n\trf)_{\dff q}
\off =\off\dff
\mathcal{Z}^{\dff n}\dff(\dff \bm{\Delta}\fff [\dff q\trf]\fff,\qff \pi\trf)
\pff,
\]

\vspace{-9pt}\vspace{1.125pt}
and\dss the structural\dss maps\dss
$\theta^{\dff *}\dff \colon\dff
K\dff(\dff \pi\fff,\qff n\trf)_{\dff q}
\qff \ttoo\qff
K\dff(\dff \pi\fff,\qff n\trf)_{\dff r}$\dss
are\sss the above induced\dss homomorphisms\dss 
$\theta^{\dff *}$\dnsp.\oss
Let\dss 
$0_{\fff q}
\qff \in\pff 
\mathcal{Z}^{\dff n}\dff(\dff \bm{\Delta}\fff [\dff q\trf]\fff,\qff \pi\trf)
\off =\off 
K\dff(\dff \pi\fff,\qff n\trf)_{\dff q}$\dss
be\sss the zero cocycle.\oss

\myuppar{The $n$\dnsp-cocycles of\trs $K\dff(\dff \pi\fff,\qff n\trf)$\nnsp.}
Every\dss normalized\sss $n$\dnsp-cochain 
of\dss $\bm{\Delta}\fff [\halfff n\dff]$\sss is\dss a cocycle,\oss
i.e.\qss\vspace{3pt}\vspace{0.75pt}
\[
\quad
\mathcal{Z}^{\dff n}\dff(\dff \bm{\Delta}\fff [\halfff n\dff]\fff,\qff \pi\trf)
\off =\off
\mathcal{C}^{\dff n}\dff(\dff \bm{\Delta}\fff [\halfff n\dff]\fff,\qff \pi\trf)
\pff.
\]

\vspace{-9pt}\vspace{0.75pt}
Clearly,\oss a normalized $n$\dnsp-cochain $c$ of\dss $\bm{\Delta}\fff [\halfff n\dff]$\sss is\dss
determined\dss by\sss its value\dss $c\trf(\qff \bm{\iota}_{\dff n}\trf)$\sss
on\dss the unique non-de\-gen\-er\-ate $n$\dnsp-simplex $\bm{\iota}_{\dff n}$\sss
of\dss $\bm{\Delta}\fff [\halfff n\dff]$\nnsp.\oss
Therefore we can\dss identity\dss $K\dff(\dff \pi\fff,\qff n\trf)_{\dff n}$\sss
with\sss $\pi$\nnsp.\oss

A $n$\dnsp-cochain of\dss $K\dff(\dff \pi\fff,\qff n\trf)$\sss
with coefficients in $\pi$
is\dss a\sss map\dss
$K\dff(\dff \pi\fff,\qff n\trf)_{\dff n}\qff \ttoo\qff \pi$\nnsp,\oss
and\dss hence can\sss be\sss thought\sss as a map\dss 
$c\dff \colon\dff \pi\qff \ttoo\qff \pi$\nnsp.\oss
Clearly,\oss the zero cocycle\dss 
$0_{\fff n}
\qff \in\pff 
\mathcal{Z}^{\dff n}\dff(\dff \bm{\Delta}\fff [\halfff n\dff]\fff,\qff \pi\trf)$\sss
is\dss the only\sss degenerate 
$n$\dnsp-simplex of\dss $K\dff(\dff \pi\fff,\qff n\trf)$\nnsp.\oss
Therefore we can identify\dss normalized $n$\dnsp-cochains $c$ of\dss
$K\dff(\dff \pi\fff,\qff n\trf)$\sss
with\sss maps\sss
$c\dff \colon\dff \pi\qff \ttoo\qff \pi$\sss
subject\sss only\dss to\sss the condition\sss
$c\trf(\dff 0\dff)\off =\off 0$\nnsp.\oss
It\dss turns out\dss that\sss 
$c$\sss is\dss a\sss cocycle\dss if\trs and\dss only\trs if\trs
$c$\sss is\dss a\sss homomorphisms\dss
$\pi\qff \ttoo\qff \pi$\nnsp.\oss 
See\trs Lemma\qss \ref{cocycle-homomorphism}\qss for a proof\halfff.\oss
In\dss particular\halfff,\oss
the identity\dss map\dss
$\id_{\trf \pi}\dff \colon\dff
\pi\qff \ttoo\qff \pi$\dss
is\dss an $n$\dnsp-cocycle of\trs $K\dff(\dff \pi\fff,\qff n\trf)$\nnsp.\oss

\myuppar{Simplicial\dss maps\sss to\sss $K\dff(\dff \pi\fff,\qff n\trf)$\nnsp.}
Let\dss $K$\sss be a simplicial\sss set.\oss 
Let\dss us assign\dss to every\sss simplicial\dss map\dss
$f\dff \colon\dff
K\qff \ttoo\qff K\dff(\dff \pi\fff,\qff n\trf)$\dss
the $n$\dnsp-cocycle\vspace{3pt}\vspace{0.75pt}
\[
\quad
z\trf(\trf f\trf)
\off =\off
f^{\dff *}\dff\bigl(\trf \id_{\trf \pi}\trf\bigr)
\off \in\off
\mathcal{Z}^{\dff n}\dff(\trf K\fff,\qff \pi\trf)
\pff
\]

\vspace{-9pt}\vspace{0.75pt}
({\fff}the cochain\sss $z\trf(\trf f\trf)$\sss is\dss a\sss
normalized cocycle because\sss $\id_{\trf \pi}$\sss is).\oss
Unraveling\dss the definitions shows\sss that\sss $f$\sss
is\dss uniquely\sss determined\dss by\sss $z\trf(\trf f\trf)$\sss
and\dss for every\dss
$z\qff \in\qff \mathcal{Z}^{\dff n}\dff(\trf K\fff,\qff \pi\trf)$\dss
there\dss is\dss a map\dss
$f\dff \colon\dff
K\qff \ttoo\qff K\dff(\dff \pi\fff,\qff n\trf)$\dss
such\dss that\dss
$z\trf(\trf f\trf)
\off =\off
z$\nnsp.\oss
One can say\dss that\dss the definition of\trs $K\dff(\dff \pi\fff,\qff n\trf)$\dss
is\dss dictated\dss by\dss this property.\oss
If\dss $\pi$\sss is\dss abelian,\oss
then\dss two maps\dss
$f,\pff g\dff \colon\dff
K\qff \ttoo\qff K\dff(\dff \pi\fff,\qff n\trf)$\dss
are homotopic\dss
if\trs and\dss only\trs if\trs
$z\trf(\trf f\trf)\qff -\qff z\trf(\trf g\trf)$\dss
is\dss a coboundary\dss of\dss a normalized $(\dff n\qff -\qff 1\dff)$\dnsp-cochain.\oss
See,\oss for example,\oss \cite{m},\oss Lemma\qss 24.3\qss and\qss Theorem\qss 24.4,\oss
or\pss \cite{em3},\oss Theorems\qss 5.1\qss and\qss 5.2.\oss

In\dss the case of\trs
$K\off =\off K\dff(\dff \pi\fff'\fff,\qff n\trf)$\nnsp,\oss
we see\sss that\sss every\dss homomorphism\dss
$h\dff \colon\dff \pi\fff'\qff \ttoo\qff \pi$\sss
defines a\sss map\dss
$K\dff(\dff \pi\fff'\fff,\qff n\trf)\qff \ttoo\qff K\dff(\dff \pi\fff,\qff n\trf)$\nnsp.\oss
We will\sss denote\sss this map\sss by\dss $\bm{s}\trf(\dff h\trf)$\nnsp.\oss
There are no other simplicial\dss maps\dss
$K\dff(\dff \pi\fff'\fff,\qff n\trf)\qff \ttoo\qff K\dff(\dff \pi\fff,\qff n\trf)$\nnsp,\oss
and\sss $\bm{s}\trf(\dff h\trf)$\sss is\dss an\sss isomorphisms\dss
if\trs and\dss only\trs if\trs $h$\sss is.\oss

\myuppar{Maps over\dss the base.}
Let\dss $p\dff \colon\dff E\qff \ttoo\qff B$\dss 
be a\sss simplicial\dss map\sss thought\sss of\dss as a bundle.\oss
A simplicial\dss map\dss
$f\dff \colon\dff E\qff \ttoo\qff E$\dss
is\dss said\dss to be a\qss \emph{map over\dss the base},\pss
or\sss a\qss \emph{map over}\trs $B$\nnsp,\oss if\trs 
$p\trf \circ\dff f\off =\off f$\dnsp.\oss
For example,\oss the identity\dss map\dss
$\id_{\trf E}\dff \colon\dff E\qff \ttoo\qff E$\dss
is\dss a map over\dss the base.\oss\vspace{-0.125pt}

\myuppar{Simplicial\sss groups.}
A\qss \emph{simplicial\sss group}\pss is\dss a contravariant\dss functor\dss
from\dss $\bm{\Delta}$\dss to\sss the category\sss of\dss groups.\oss
In other\sss words,\oss a simplicial\dss group\dss is\dss a simplicial\sss set\sss
$G$\sss together with group structures on\sss sets\sss $G_{\dff n}$\nsp,\qss
$n\qff \in\pff \nnn$\nnsp,\oss 
such\dss that\dss  the structural\dss maps\sss $\theta^{\dff *}$\sss are homomorphisms.\oss
Let\sss $G$\sss be a simplicial\dss group.\oss
Then\sss for every\dss $q\qff \in\pff \nnn$\dss there\dss is\dss a natural\sss
action of\trs the group\sss $G_{\dff q}$\sss
on\dss $\bm{\Delta}\fff [\dff q\trf]\dff \times\dff G$\dss
defined as follows.\oss
Let\dss
$g\qff \in\qff G_{\dff q}$\nsp.\oss
The $m$\dnsp-simplices of\dss $\bm{\Delta}\fff [\dff q\trf]\dff \times\dff G$\sss
are\sss the pairs\sss $(\dff \theta\fff,\qff \tau\trf)$\nnsp,\oss
where\sss $\tau\qff \in\qff G_{\dff m}$\dss and\dss
$\theta\dff \colon\dff
[\halfff m\dff]\qff \ttoo\qff [\dff q\trf]$\sss
is\dss a $q$\dnsp-simplex of\dss $\bm{\Delta}\fff [\dff q\trf]$\nnsp.\oss
Let\vspace{3pt}\vspace{0.75pt}
\[
\quad
g\dff \cdot\dff(\dff \theta\fff,\qff \tau\trf)
\off =\off
\left(\qff
\theta\fff,\pff
\theta^{\dff *}\dff(\dff g\trf)\dff \cdot\dff \tau
\qff\right)
\pff,
\]

\vspace{-9pt}\vspace{0.75pt}
where\sss the product\dss in\dss the right\sss hand side\dss is\dss 
taken\sss in\dss $G_{\dff q}$\nsp.\oss
A routine check shows\sss that\dss the map\dss
$(\dff \theta\fff,\qff \tau\trf)
\off \longmapsto\off
g\dff \cdot\dff(\dff \theta\fff,\qff \tau\trf)$\dss
is\dss a simplicial\dss map\dss
$\bm{\Delta}\fff [\dff q\trf]\dff \times\dff G
\qff \ttoo\qff 
\bm{\Delta}\fff [\dff q\trf]\dff \times\dff G$\nnsp.\oss
We will\sss denote\sss this map\sss by\dss $\bm{t}\trf(\dff g\trf)$\nnsp.\oss
Clearly,\pss 
$\bm{t}\trf(\dff g\trf)$\sss
is\dss an automorphism of\trs the\sss bundle\dss
$\pr\dff \colon\dff
\bm{\Delta}\fff [\dff q\trf]\dff \times\dff G
\qff \ttoo\qff
\bm{\Delta}\fff [\dff q\trf]$\dss
over\dss the base.\oss
Another\sss routine check shows\sss that\dss
$g\off \longmapsto\off \bm{t}\trf(\dff g\trf)$\dss
is\dss an\sss action.\oss

The addition of\dss normalized cocycles\sss turns\dss
$K\dff(\dff \pi\fff,\qff n\trf)$\dss into a simplicial\dss group.\oss
In\dss the case of\dss
$G\off =\off K\dff(\dff \pi\fff,\qff n\trf)$\sss
we\sss get\sss an\sss automorphism\vspace{3pt}\vspace{0.75pt}
\begin{equation}
\label{ac}
\quad
\bm{t}\trf(\dff c\trf)\dff \colon\dff
\bm{\Delta}\fff [\dff q\trf]\dff \times\dff K\dff(\dff \pi\fff,\qff n\trf)
\qff \ttoo\qff 
\bm{\Delta}\fff [\dff q\trf]\dff \times\dff K\dff(\dff \pi\fff,\qff n\trf)
\pff
\end{equation}

\vspace{-9pt}\vspace{0.75pt}
of\trs
the\sss bundle\vspace{3pt}\vspace{0.75pt}
\begin{equation}
\label{pr-to-simplex}
\quad
\pr\dff \colon\dff
\bm{\Delta}\fff [\dff q\trf]\dff \times\dff K\dff(\dff \pi\fff,\qff n\trf)
\qff \ttoo\qff 
\bm{\Delta}\fff [\dff q\trf]
\pff.
\end{equation}

\vspace{-9pt}\vspace{0.75pt}
over\dss the base for every\sss normalized cocycle\dss
$c\qff \in\qff
\mathcal{Z}^{\dff n}\dff(\dff \bm{\Delta}\fff [\dff q\trf]\fff,\qff \pi\trf)$\nnsp.\oss

\mypar{Theorem.}{aut-over-simplex}
\emph{Every\dss automorphism of\oss 
\textup{(\ref{pr-to-simplex})}\qss
over\dss the base  
is\dss equal\dss to\sss the composition}\dss\vspace{3pt}\vspace{0.75pt}
\[
\quad
\left(\trf
\id_{\dff \bm{\Delta}\fff [\dff q\trf]}
\dff \times\qff
\bm{s}\trf(\dff h\trf)
\trf\right) 
\pff \circ\pff 
\bm{t}\trf(\dff c\trf)
\pff,
\]

\vspace{-9pt}\vspace{0.75pt}
\emph{where\sss 
$h\dff \colon\dff
\pi\qff \ttoo\qff \pi$\dss
is\dss an automorphism\dss
and\dss
$c\qff \in\qff
\mathcal{Z}^{\dff n}\dff(\dff \bm{\Delta}\fff [\dff q\trf]\fff,\qff \pi\trf)$\nnsp.\oss
Both\sss $h$\sss and\sss $c$ are unique\-ly\sss determined\dss by\dss the\sss
automorphism.\oss}

\proof
See\qss \cite{m},\oss Propositions\qss 25.2\qss and\qss 25.3.\oss  \eproof

\myuppar{Translations of\trs trivial\dss bundles.}
Let\dss $p\dff \colon\dff E\qff \ttoo\qff \bm{\Delta}\fff [\dff q\trf]$\dss
be a\sss trivial\dss bundle with\dss the fiber\dss $K\dff(\dff \pi\fff,\qff n\trf)$\nnsp,\qss
$n\qff >\qff 1$\nnsp,\oss
and\dss let\dss
$f\dff \colon\dff E\qff \ttoo\qff E$\sss be an automorphism over\dss the base.\oss
We will\sss say\dss that\sss $f$\sss is\dss a\qss \emph{translation}\oss
if\trs
$f\off =\off
t\dff \circ\dff \bm{t}\trf(\dff c\trf)\dff \circ\dff t^{\dff -\dff 1}$\dss
for some\sss trivialization\dss
$t\dff \colon\dff
\bm{\Delta}\fff [\dff q\trf]\dff \times\dff K\dff(\dff \pi\fff,\qff n\trf)
\qff \ttoo\qff
E$\dss 
and some normalized cocycle\dss
$c\qff \in\qff
\mathcal{Z}^{\dff n}\dff(\dff \bm{\Delta}\fff [\dff q\trf]\fff,\qff \pi\trf)$\nnsp.\oss

Theorem\qss \ref{aut-over-simplex}\qss implies\sss that\dss if\dss $f$\sss
has\sss the form\dss
$t\dff \circ\dff \bm{t}\trf(\dff c\trf)\dff \circ\dff t^{\dff -\dff 1}$\dss
for some\sss trivialization\sss $t$\nnsp,\oss
then\sss $f$\sss has such\sss form\sss for every\sss trivialization.\oss
But\dss the cocycle\sss $c$\sss depends on\dss the choice of\dss $t$\sss
because\qss (\ref{pr-to-simplex})\qss has automorphisms of\trs the form\dss
$\id_{\dff \bm{\Delta}\fff [\dff q\trf]}
\dff \times\qff
\bm{s}\trf(\dff h\trf)$\nnsp.\oss
Still,\oss there\dss is\dss a way\dss to make\sss $c$\sss to be uniquely\sss
determined\dss by\sss $f$\dnsp.\oss

Let\dss 
$i_{\trf 0}\dff \colon\dff
\bm{\Delta}\fff [\dff 0\dff]
\qff \ttoo\qff
\bm{\Delta}\fff [\dff q\trf]$\dss
be\sss the map defined\dss by\dss the inclusion\dss
$[\dff 0\dff]\qff \ttoo\qff [\dff q\trf]$\nnsp.\oss
The\sss total\sss space\sss $F$ of\trs the pull-back\dss bundle\sss
$i_{\trf 0}^{\dff *}\dff p$\dss is\dss isomorphic\sss to\dss
$K\dff(\dff \pi\fff,\qff n\trf)$\nnsp.\oss
Any\dss two isomorphisms differ\sss by\sss an automorphism of\trs
the form\sss $\bm{s}\trf(\dff h\trf)$\nnsp,\oss
where\sss $h$\sss is\dss an automorphism of\dss $\pi$\nnsp.\oss
Therefore\sss $F$\sss is\dss a simplicial\dss group
isomorphic\sss to\sss $K\dff(\dff \pi\fff,\qff n\trf)$\nnsp.\oss
In\dss particular\halfff,\oss the set\dss $F_{\dff n}$\sss of\dss $n$\dnsp-simplices of\trs
$F$\sss is\dss a group\sss isomorphic\sss to\sss $\pi$\nnsp.\oss
Let\dss us\sss denote\sss this group by\sss $\pi_{\trf 0}$\nsp.\oss
Then\dss $p$\sss is\dss also a\sss trivial\dss bundle with\sss the fiber\dss
canonically\dss isomorphic\sss to\dss
$K\dff(\dff \pi_{\trf 0}\fff,\qff n\trf)$\nnsp.\oss
Let\dss us\sss call\sss a
trivialization\sss
$t\dff \colon\dff
\bm{\Delta}\fff [\dff q\trf]\dff \times\dff K\dff(\dff \pi_{\trf 0}\fff,\qff n\trf)
\qff \ttoo\qff
E$\qss 
\emph{special}\pss if\trs the induced\dss map\dss
$K\dff(\dff \pi_{\trf 0}\fff,\qff n\trf)
\qff \ttoo\qff
F$\sss
is\dss the canonical\dss isomorphism.\oss 
By\qss Theorem\qss \ref{aut-over-simplex}\qss 
two special\dss trivializations differ\sss by\sss a\sss translation.\oss
But,\oss if\dss\vspace{3pt}\vspace{-1.25pt}
\[
\quad 
g\dff \colon\dff
\bm{\Delta}\fff [\dff q\trf]\dff \times\dff K\dff(\dff \pi_{\trf 0}\fff,\qff n\trf)
\qff \ttoo\qff 
\bm{\Delta}\fff [\dff q\trf]\dff \times\dff K\dff(\dff \pi_{\trf 0}\fff,\qff n\trf)
\]

\vspace{-9pt}\vspace{-1.25pt}
is\dss a\sss translation,\oss then\dss 
$g\dff \circ\dff \bm{t}\trf(\dff c\trf)\dff \circ\dff g^{\dff -\dff 1}
\off =\off
\bm{t}\trf(\dff c\trf)$\nnsp.\oss
It\dss follows\dss that\sss for every\sss translation\sss $f$\sss of\trs
the bundle\sss $p$\sss
there\dss is\dss a\sss well\sss defined\sss cocycle\dss\vspace{3pt}\vspace{-1.25pt}
\[
\quad 
d\trf(\trf f\trf)
\pff \in\off
\mathcal{Z}^{\dff n}\dff(\dff \bm{\Delta}\fff [\dff q\trf]\fff,\qff \pi_{\trf 0}\trf)
\]

\vspace{-9pt}\vspace{-1.25pt}
such\dss that\dss
$t^{\dff -\dff 1}\dff \circ\dff f\dff \circ\trf t
\off =\off
\bm{t}\trf(\trf d\trf(\dff f\trf)\trf)$\dss
for every\sss special\dss trivialization\sss $t$\nnsp.\oss

\myuppar{The canonical\dss local\sss system.}
Let\dss $p\dff \colon\dff E\qff \ttoo\qff B$\dss 
be a\sss locally\dss trivial\dss bundle with\dss the fiber\dss
$K\dff(\dff \pi\fff,\qff n\trf)$\sss and\dss $n\qff \geq\qff 1$\nnsp.\oss
The above discussion of\trs translations suggests\sss to associate with
each vertex\dss $v\qff \in\pff B_{\dff 0}$\sss
a\sss group\sss $\pi_{\dff v}$\sss isomorphic\sss to\sss $\pi$\nnsp.\oss
Namely,\oss the\sss total\sss space\sss $F_{\fff v}$ of\trs the pull-back\dss bundle\sss
$i_{\dff v}^{\dff *}\trf p$\dss is\dss
\dss a simplicial\dss group
isomorphic\sss to\sss $K\dff(\dff \pi\fff,\qff n\trf)$\nnsp.\oss
Let\sss $\pi_{\dff v}$\sss be\sss the group of\dss 
$n$\dnsp-simplices of\trs $F_{\fff v}$\nsp.\oss 
Then\dss $F_{\fff v}$\sss is\dss canonically\dss isomorphic\sss to\dss
$K\dff(\dff \pi_{\dff v}\fff,\qff n\trf)$\nnsp.\oss

Suppose\sss that\sss $v\fff,\pff w\qff \in\pff B_{\dff 0}$\sss
and\sss $\varepsilon$\sss is\dss a $1$\dnsp-simplex of\dss $B$\sss such\dss that\dss
$\partial_{\dff 1}\dff \varepsilon\off =\off v$\nnsp,\pss
$\partial_{\dff 0}\dff \varepsilon\off =\off w$\nnsp.\oss
Let\dss $E_{\dff \varepsilon}$\sss be\sss the\dss total\sss space
of\trs the pull-back\dss bundle\dss
$i_{\dff \varepsilon}^{\dff *}\trf p$\nnsp,\oss
and\dss let\dss
$t\dff \colon\dff
\bm{\Delta}\fff [\dff 1\dff]\dff \times\dff K\dff(\dff \pi_{\dff v}\fff,\qff n\trf)
\qff \ttoo\qff
E_{\dff \varepsilon}$\dss be a\sss special\dss trivialization.\oss 
Recall\dss that\dss two special\dss
trivializations differ\sss by\sss a\sss translation.\oss
Since\dss $n\qff >\qff 1$\nnsp,\oss
every\dss normalized $n$\dnsp-chain of\dss $\bm{\Delta}\fff [\dff 1\dff]$\sss
is\dss equal\dss to $0$\nnsp.\oss
It\dss follows\dss that\sss every\dss translation\dss is\dss
equal\dss to\sss the identity\sss and\dss hence\sss
$t$\sss is\dss uniquely\sss determined.\oss
Therefore,\oss the isomorphism\vspace{3pt}\vspace{-1.25pt}
\[
\quad
K\dff(\dff \pi_{\dff v}\fff,\qff n\trf)
\qff \ttoo\qff 
F_{\fff w}
\off =\off
K\dff(\dff \pi_{\dff w}\fff,\qff n\trf) 
\]

\vspace{-9pt}\vspace{-1.25pt}
induced\dss by\sss $t$\sss depends only on\sss $\varepsilon$\nnsp.\oss
Let\sss $\varepsilon\dff(\dff p\trf)$\sss be\sss its\dss inverse,\oss
and\dss let\dss
$\varepsilon^{\dff *}\dff \colon\dff
\pi_{\dff w}\qff \ttoo\qff \pi_{\dff v}$\sss
be\sss the unique isomorphism such\dss that\dss
$\varepsilon\dff(\dff p\trf)
\off =\off
\bm{s}\trf(\dff \varepsilon^{\dff *}\dff)$\nnsp.\oss 
By\dss using\dss trivializations of\trs the pull-back\sss bundles\sss
$i_{\dff \sigma}^{\dff *}\trf p$\sss for $2$\dnsp-simplices $\sigma$ of\dss $B$\sss
one can easily\sss check\dss that\dss the groups\sss $\pi_{\dff v}$\sss
together\sss with\sss isomorphisms\sss $\varepsilon^{\dff *}$\sss
form a\sss local\sss system of\dss coefficients on\dss $B$\nnsp,\oss
which we will\sss denote by\sss $\pi\trf(\dff p\trf)$\nnsp.

\myuppar{Translations\sss of\trs locally\dss trivial\dss bundles.}
Let\dss $p\dff \colon\dff E\qff \ttoo\qff B$\dss
be a\sss locally\dss trivial\dss bundle with\dss 
the fiber\dss $K\dff(\dff \pi\fff,\qff n\trf)$\nnsp,\qss
$n\qff \geq\qff 1$\nnsp.\oss
Let\dss
$f\dff \colon\dff
E\qff \ttoo\qff E$\dss
be an automorphism over\dss the base.\oss
For every $q$\dnsp-simplex $\sigma$ of\dss $B$\sss the automorphism\sss $f$\sss
induces an automorphism\sss $f_{\dff \sigma}$\sss of\trs the pull-back\dss bundle\sss
$i_{\dff \sigma}^{\dff *}\trf p$\sss over\dss the base.\oss
Since\sss the bundle\sss $i_{\dff \sigma}^{\dff *}\trf p$\sss is\dss trivial,\oss
it\dss make sense\sss to ask\dss if\dss $f_{\dff \sigma}$\sss
is\dss a\sss translation
and\dss to call\sss $f$\sss a\qss \emph{translation}\pss
if\dss $f_{\dff \sigma}$\sss
is\dss a\sss translation for every\sss simplex $\sigma$ of\dss $B$\nnsp.\oss
Suppose\sss that\dss
$f\dff \colon\dff
E\qff \ttoo\qff E$\dss
is\dss a\sss translation.\oss
Let\sss $\sigma$\sss be an $n$\dnsp-simplex of\trs $B$\nnsp,\oss
and\dss let\sss $\pi_{\dff \sigma}\off =\off \pi_{\dff v}$\nsp,\oss
where\dss $v\off =\off v_{\dff \sigma}$\sss is\dss the\sss
leading vertex of\dss $\sigma$\nnsp.\oss
Then\vspace{3pt}\vspace{-1.25pt}
\[
\quad
d\trf(\trf f_{\dff \sigma}\trf)
\pff \in\off
\mathcal{Z}^{\dff n}\dff(\dff \bm{\Delta}\fff [\halfff n\dff]\fff,\qff \pi_{\dff \sigma}\trf)
\pff.
\]

\vspace{-9pt}\vspace{-1.25pt}
The group\sss
$\mathcal{Z}^{\dff n}\dff(\dff \bm{\Delta}\fff [\halfff n\dff]\fff,\qff \pi_{\dff \sigma}\trf)$\sss
is\dss canonically\dss isomorphic\sss to\sss $\pi_{\dff \sigma}$\sss
and\dss hence we can consider\sss
$d\dff(\trf f_{\dff \sigma}\trf)$\sss
as an elements of\dss $\pi_{\dff \sigma}$\nsp.\oss
The map\vspace{1.5pt}
\[
\quad
D_{\fff f}\dff \colon\dff
\sigma
\off \longmapsto\off
d\trf(\trf f_{\dff \sigma}\trf)
\]

\vspace{-10.5pt}
is\dss an $n$\dnsp-cochain of\trs $B$\sss with coefficients\sss in\dss
the\sss local\sss system\sss $\pi\trf(\dff p\trf)$\nnsp.\oss

\mypar{Lemma.}{translation-cocycle}
\emph{The $n$\dnsp-cochain\trs $D_{\fff f}$\sss
is\dss a\sss normalized\sss cocycle.\oss}

\proof
If\dss $\sigma$\sss is\dss a degenerate $n$\dnsp-simplex of\trs $B$\nnsp,\oss
then $\sigma\off =\off \theta^{\dff *}\dff(\dff \tau\trf)$
for an $m$\dnsp-simplex $\tau$ such\dss that\sss
$m\qff <\qff n$\sss and a non-decreasing\sss map\dss 
$\theta\dff \colon\dff
[\halfff n\dff]\qff \ttoo\qff [\halfff m\dff]$\nnsp.\oss
Therefore
$i_{\dff \sigma}\off =\off i_{\trf \tau}\dff \circ\dff \theta_{\dff *}$\nsp,\oss
where\dss
$\theta_{\dff *}\dff \colon\dff
\bm{\Delta}\fff [\halfff n\dff]
\qff \ttoo\qff 
\bm{\Delta}\fff [\halfff m\dff]$\dss
is\dss induced\dss by\sss $\theta$\nnsp,\oss
and\dss hence\vspace{3pt}
\begin{equation}
\label{deg-ind}
\quad
i_{\dff \sigma}^{\dff *}\trf p
\off =\off
\theta_{\dff *}^{\dff *}\dff 
\left(\qff
i_{\trf \tau}^{\dff *}\trf p
\qff\right)
\end{equation}

\vspace{-9pt}
If\dss $t$\sss be a special\dss trivialization of\trs 
$i_{\trf \tau}^{\dff *}\trf p$\nnsp,\oss
then
$t^{\dff -\dff 1}\dff \circ\dff f_{\trf \tau}\dff \circ\trf t
\off =\off
\bm{t}\trf(\trf d\trf(\dff f_{\dff \tau}\trf)\trf)$\nnsp.\oss
But\sss $d\trf(\dff f_{\trf \tau}\trf)$\sss is\dss a normalized $m$\dnsp-cochain
of\dss $\bm{\Delta}\fff [\halfff m\dff]$\nnsp.\oss
Therefore $m\qff <\qff n$\sss implies\sss that\sss
$d\trf(\dff f_{\trf \tau}\trf)\off =\off 0$\sss
and\dss hence\sss $f_{\trf \tau}$\sss
is\dss equal\dss to\sss the identity.\oss
In\sss view of\qss (\ref{deg-ind})\qss this implies\sss that\sss
$f_{\dff \sigma}$\sss is\dss equal\dss to\sss the identity and\dss hence\sss
$d\trf(\trf f_{\dff \sigma}\trf)\off =\off 0$\nnsp.\oss
It\dss follows\dss that\sss $D_{\fff f}$\sss
is\dss normalized.\oss  

Let\sss $\rho$ be an $(\dff n\qff +\qff 1\dff)$\dnsp-simplex of\trs $B$\nnsp,\oss
and\dss let\sss 
$\varepsilon\off =\off \theta^{\dff *}\dff(\trf \rho\trf)$\nnsp,\oss
where\dss
$\theta\dff \colon\dff
[\dff 1\trf]\qff \ttoo\qff [\halfff n\dff]$\sss
is\dss the inclusion.\oss
Then\sss
$v\off =\off \partial_{\dff 1}\dff \varepsilon$\sss
is\dss the\sss leading\sss vertex of\dss $\rho$
and each\sss face\sss
$\partial_{\dff i}\trf \rho$ with $i\qff >\qff 0$\nnsp,\oss
and\dss
$w\off =\off \partial_{\dff 0}\dff \varepsilon$\dss
is\dss the\sss leading\sss vertex of\trs
$\tau\off =\off \partial_{\dff 0}\trf \rho$\nnsp.\oss
Let\sss $t$\sss be a special\dss trivialization of\trs
$i_{\trf \rho}^{\dff *}\trf p$\nnsp.\oss
Then\sss $t$\sss induces a\sss trivialization 
of\trs $i_{\trf \sigma}^{\dff *}\trf p$\sss
for every\sss face\dss $\sigma\off =\off \partial_{\dff i}\trf \rho$\nnsp.\oss 
If\dss $i\qff >\qff 0$\nnsp,\oss
then\dss the induced\dss trivialization\dss is\dss special.\oss  
If\dss $i\off =\off 0$\nnsp,\oss
it\sss differs\sss from a special\dss trivialization 
by\dss the isomorphism\vspace{3pt}
\[
\quad
\varepsilon\dff(\dff p\trf)\dff \colon\dff
K\dff(\dff \pi_{\dff w}\fff,\qff n\trf)
\qff \ttoo\qff 
K\dff(\dff \pi_{\dff v}\fff,\qff n\trf)
\]

\vspace{-9pt}
corresponding\sss to\sss the isomorphism\dss
$\varepsilon^{\dff *}\dff \colon\dff
\pi_{\dff w}\qff \ttoo\qff \pi_{\dff v}$\nsp.\oss
It\dss follows\dss that\vspace{3pt}\vspace{0.42pt}
\[
\quad
d\trf\left(\trf f_{\trf \rho}\trf\right)\dff
\bigl(\trf \partial_{\dff i}\qff \bm{\iota}_{\dff n\dff +\dff 1}\trf\bigr)
\off =\off\dff
D_{\fff f}\dff\left(\trf \partial_{\dff i}\trf \rho\trf\right)
\quad\
\mbox{if}\quad\
i\qff >\qff 0
\quad\
\mbox{and}\quad\
\]

\vspace{-36pt}\vspace{0.42pt}
\[
\quad
d\trf\left(\trf f_{\trf \rho}\trf\right)\dff
\bigl(\trf \partial_{\trf 0}\qff \bm{\iota}_{\dff n\dff +\dff 1}\trf\bigr)
\off =\off\dff
\varepsilon^{\dff *}\dff\left(\trf
D_{\fff f}\dff\left(\trf \tau\trf\right)
\trf\right)
\pff.
\]

\vspace{-9pt}\vspace{0.42pt}
Since\sss $d\trf(\trf f_{\trf \rho}\trf)$\sss
is\dss a cocycle with coefficients in\sss $\pi_{\trf \rho}\off =\off \pi_{\dff v}$\nsp,\oss
this\sss implies\sss that\sss $D_{\fff f}$\sss is\dss a cocycle
with coefficients in\dss the local\sss system\sss 
$\pi\trf(\dff p\trf)$\nnsp.\oss  \eproof

\mypar{Lemma.}{cocycles-translations}
\emph{For every\dss normalized $n$\dnsp-cocycle
$c\qff \in\pff
\mathcal{Z}^{\dff n}\dff(\trf B\fff,\qff \pi\trf(\dff p\trf)\trf)$
there exists a unique\sss translation 
$f\off =\off f\dff(\dff c\trf)\dff \colon\dff
E\qff \ttoo\qff E$ 
such\dss that\qss
$D_{\fff f}\off =\off c$\nnsp.\oss}

\proof
Let\sss $\sigma$\sss be a $q$\dnsp-simplex of\trs $B$\nnsp,\oss
and\dss let\sss $t$\sss be a special\dss trivialization of\trs 
$i_{\dff \sigma}^{\dff *}\trf p$\nnsp.\oss
Then\sss $t$\sss induces a\sss trivialization 
of\trs $i_{\trf \tau}^{\dff *}\trf p$\sss
for every\sss simplex $\tau$ of\trs the form\dss
$\tau\off =\off \theta^{\dff *}\dff(\dff \dff \sigma\trf)$\nnsp.\oss
In\dss general,\oss the induced\dss trivialization\dss is\dss not\sss
special,\oss but\sss differs from a special\dss one 
by\dss the isomorphism\vspace{3pt}
\[
\quad
\varepsilon\dff(\dff p\trf)\dff \colon\dff
K\dff(\dff \pi_{\dff w}\fff,\qff n\trf)
\qff \ttoo\qff 
K\dff(\dff \pi_{\dff v}\fff,\qff n\trf)
\]

\vspace{-9pt}
where $v\fff,\pff w$ are\sss the\sss leading\sss vertices of\trs
$\dff \sigma\fff,\pff \tau$ respectively,\oss
and $\varepsilon$\sss is\dss the unique $1$\dnsp-simplex of\trs the form\dss
$\varepsilon\off =\off \eta^{\dff *}\dff(\trf \rho\trf)$\sss 
such\dss that\dss
$\partial_{\dff 1}\dff \varepsilon\off =\off v$\nnsp,\pss
$\partial_{\dff 0}\dff \varepsilon\off =\off w$\nnsp.\oss
It\dss follows\dss that\vspace{3pt}
\[
\quad
\theta_{\dff *}^{\dff *}\dff 
\left(\qff
d\trf(\trf f_{\dff \dff \sigma}\trf)
\qff\right)
\off =\off
\varepsilon^{\dff *}\dff
\left(\qff
d\trf(\trf f_{\dff \tau}\trf)
\qff\right)
\pff,
\]

\vspace{-9pt}
where\sss the isomorphism\dss
$\varepsilon^{\dff *}\dff \colon\dff
\pi_{\dff w}\qff \ttoo\qff \pi_{\dff v}$\dss
is\dss applied\dss to\sss the coefficients of\dss
$d\trf(\trf f_{\dff \tau}\trf)$\nnsp.\oss
By\sss applying\dss this observation\dss to $n$\dnsp-simplices $\tau$
we see\sss that\sss $d\trf(\trf f_{\dff \sigma}\trf)$\sss
is\dss determined\dss by\dss $D_{\fff f}$\nsp,\oss
and\dss hence\vspace{3pt}
\[
\quad
f_{\dff \sigma}\dff \colon\dff
i_{\dff \sigma}^{\dff *}\trf p
\qff \ttoo\qff
i_{\dff \sigma}^{\dff *}\trf p
\]

\vspace{-9pt}
is\dss also determined\sss by\dss $D_{\fff f}$\nsp.\oss
Since\sss this\dss is\dss true for every simplex $\sigma$ of\trs $B$\nnsp,\oss
the\sss translation\sss $f$\sss is\dss determined\sss by\dss $D_{\fff f}$\nsp.\oss
This proves\sss the uniqueness.\oss
To prove\sss the existence,\oss suppose\sss that
$c\qff \in\pff
\mathcal{Z}^{\dff n}\dff(\trf B\fff,\qff \pi\trf(\dff p\trf)\trf)$
is\dss given.\oss
Let\sss $\sigma$\sss be a $q$\dnsp-simplex of\trs $B$\nnsp.\oss
The isomorphisms\sss $\varepsilon^{\dff *}$\sss from\sss the first\dss part\sss
of\trs the proof\dss establish an\sss isomorphism\sss between\dss the induced\dss
local\sss system\sss $i_{\dff \sigma}^{\dff *}\trf \pi\trf(\dff p\trf)$\sss
and\dss the constant\sss coefficients system\sss $\pi_{\dff v}$\nsp.\oss
This isomorphism\sss turns\sss the $n$\dnsp-cochain\sss
$i_{\dff \sigma}^{\dff *}\dff(\dff c\trf)$\sss of\dss 
$\bm{\Delta}\fff [\dff q\trf]$
with coefficients in\sss
$i_{\dff \sigma}^{\dff *}\trf \pi\trf(\dff p\trf)$\sss
into an $n$\dnsp-cochain\vspace{3pt}\vspace{0.5pt}
\[
\quad
c\trf(\dff \sigma\trf)
\pff \in\off
\mathcal{Z}^{\dff n}\dff(\trf \bm{\Delta}\fff [\dff q\trf]\fff,\qff \pi_{\dff v}\trf)
\pff.
\]

\vspace{-9pt}\vspace{0.5pt}
There\dss is\dss a\sss unique\sss translation\sss
$f_{\dff \sigma}\dff \colon\dff
i_{\dff \sigma}^{\dff *}\trf p
\qff \ttoo\qff
i_{\dff \sigma}^{\dff *}\trf p$
such\dss that\sss
$d\trf(\trf f_{\dff \sigma}\trf)
\off =\off
c\trf(\dff v\trf)$\nnsp.\oss
Since\sss the co\-chains $c\trf(\dff \sigma\trf)$\sss result\dss from a single cochain $c$\nnsp,\oss
the\sss translations\sss $f_{\dff \sigma}$\sss agree with each other\sss in\sss
the sense\sss that\trs if\trs
$\tau\off =\off \theta^{\dff *}\dff(\dff \sigma\trf)$\nnsp,\oss
then\sss the diagram\vspace{0pt}
\[
\quad
\begin{tikzcd}[column sep=booms, row sep=booms]\dis
i_{\dff \tau}^{\dff *}\trf E
\arrow[r, "\dis \theta_{\dff *}\dff"]
\arrow[d, "\dis f_{\dff \tau}\dff"']
&
i_{\dff \sigma}^{\dff *}\trf E
\arrow[d, "\dis \dff f_{\dff \sigma}"]
\\
i_{\dff \tau}^{\dff *}\trf E 
\arrow[r, "\dis \theta_{\dff *}"]
&
i_{\dff \sigma}^{\dff *}\trf E
\end{tikzcd}
\]

\vspace{-9pt}
is\dss commutative.\oss
It\dss follows\dss that\dss the maps\sss $f_{\dff \sigma}$\sss
together define a\sss translation\sss
$f\dff \colon\dff E\qff \ttoo\qff E$\nnsp.\oss
By\dss the construction,\pss
$D_{\fff f}\off =\off c$\nnsp.\oss
This proves\sss the existence.\oss  \eproof

\mypar{Lemma.}{translations-composition}
\emph{Let\qss $f,\pff g\qff \colon\dff E\qff \ttoo\qff E$\qss
be\sss two\sss translations.\oss
Then}\vspace{3pt}\vspace{0.5pt}
\[
\quad
D_{\fff f\dff \circ\trf g}
\off =\off
D_{\fff f}\off +\off D_{\dff g} 
\pff.
\]

\vspace{-9pt}\vspace{0.5pt}
\emph{If\pss
$c\fff,\pff d\qff \in\pff
\mathcal{Z}^{\dff n}\dff(\trf B\fff,\qff \pi\trf(\dff p\trf)\trf)$\nnsp,\oss
then\dss
$f\dff(\dff c\qff +\qff d\qff)
\off =\off
f\dff(\dff c\trf)
\trf \circ\dff
f\dff(\dff d\qff)$\nnsp.\oss}

\proof
The first\dss part\sss of\trs the\sss lemma\sss follows directly\dss
from\dss the definitions.\oss
In\sss view of\qss Lemma\qss \ref{cocycles-translations}\qss
the second\dss part\dss follows\sss from\dss the first\sss one.\oss  \eproof

\mypar{Lemma.}{translations-homotopy}
\emph{Let\qss
$c\fff,\pff d\qff \in\pff
\mathcal{Z}^{\dff n}\dff(\trf B\fff,\qff \pi\trf(\dff p\trf)\trf)$\nnsp.\oss
If\pss $c\qff -\qff d$\dss
is\dss equal\dss to\sss the coboundary of\dss a\sss normalized
cochain,\oss
then\dss the maps\dss
$f\dff(\dff c\trf)$\dss and\dss $f\dff(\dff d\qff)$\dss
are homotopic.\oss}

\proof
In\sss view of\qss Lemma\qss \ref{translations-composition}\qss
it\dss is\dss sufficient\dss to consider\dss the case when\dss $d\off =\off 0$\nnsp.\oss
Suppose\sss that\qss
$b\pff \in\off
\mathcal{C}^{\dff n\dff -\dff 1}\dff(\trf B\fff,\qff \pi\trf(\dff p\trf)\trf)$\dss
and\dss $c\off =\off\fff \partial^{\dff *}\fff b$\nnsp.\oss
In\dss this case we need\dss to prove\sss that\sss $f\dff(\dff c\trf)$\sss
is\dss homotopic\sss to\sss the identity.\oss
Let\dss us\sss consider\dss the bundle\vspace{3pt}\vspace{-0.5pt}
\[
\quad
p_{\dff 1}
\off =\off
p\dff \times\dff \id_{\dff \bm{\Delta}\fff [\dff 1\dff]}
\qff \colon\qff 
E\dff \times\dff \bm{\Delta}\fff [\dff 1\dff]
\pff \ttoo\pff
B\dff \times\dff \bm{\Delta}\fff [\dff 1\dff]
\pff.
\]

\vspace{-9pt}\vspace{-0.5pt}
Equivalently,\oss $p_{\dff 1}$\dss is\dss induced\dss from\sss $p$\sss
by\dss the projection\dss
$\pr_{\dff B}\dff \colon\dff
B\dff \times\dff \bm{\Delta}\fff [\dff 1\dff]
\qff \ttoo\qff
B$\nnsp.\oss
Clearly,\oss the\sss local\sss system\sss $\pi\trf(\dff p_{\dff 1}\trf)$\sss
is\dss induced\dss from\sss $\pi\trf(\dff p\trf)$\sss by\dss
the same projection.\oss
Recall\dss the simplicial\dss maps\dss 
$i\dff(\dff e\trf)
\qff \colon\qff
\bm{\Delta}\fff[\dff 0\dff]
\qff \ttoo\qff
\bm{\Delta}\fff[\dff 1\dff]$\nnsp,\oss
where\dss $e\off =\off 0$\sss or\dss $1$\nnsp,\oss
from\dss the definition of\dss homotopies.\oss
These maps\sss lead\dss to\sss the maps\vspace{3pt}\vspace{-0.5pt}
\[
\quad
\id_{\trf B}\dff \times\qff i\dff(\dff e\trf)
\qff \colon\qff
B\dff \times\dff \bm{\Delta}\fff[\dff 0\dff]
\qff \ttoo\qff
B\dff \times\dff \bm{\Delta}\fff[\dff 1\dff]
\pff.
\]

\vspace{-9pt}\vspace{-0.5pt}
Let\trs $B_{\dff e}$\dss be\sss the image of\trs
$\id_{\trf B}\dff \times\qff i\dff(\dff e\trf)$\nnsp.\oss
Similarly,\oss let\trs $E_{\dff e}$\dss be\sss the image of\trs
$\id_{\trf E}\dff \times\qff i\dff(\dff e\trf)$\nnsp.\oss
Let\dss us\dss identify\dss $B_{\trf 0}$\dss with\dss $B$\sss
and consider\dss the $(\dff n\qff -\qff 1\dff)$\dnsp-cochain\vspace{3pt}\vspace{-0.5pt}
\[
\quad
b_{\trf 0}
\pff \in\off
\mathcal{C}^{\dff n\dff -\dff 1}\dff(\trf B\dff \times\dff \bm{\Delta}\fff [\dff 1\dff]\fff,\qff 
\pi\trf(\dff p_{\dff 1}\trf)\trf)
\]

\vspace{-9pt}\vspace{-0.5pt}
equal\dss to\sss $b$\sss on\dss
$B_{\trf 0}\off =\off B$\dss
and\dss to\sss $0$\sss on all\sss 
$(\dff n\qff -\qff 1\dff)$\dnsp-simplices of\trs 
$B\dff \times\dff \bm{\Delta}\fff [\dff 1\dff]$\sss
not\dss in\sss $B_{\trf 0}$\nsp.\oss
Let\vspace{3pt}\vspace{-0.5pt}
\[
\quad
c_{\trf 0}
\off =\off 
\partial^{\dff *}\fff b_{\trf 0}
\quad\
\mbox{and}\quad\
\]

\vspace{-36pt}
\[
\quad
h
\off =\off 
f\dff(\dff c_{\trf 0}\trf)
\qff \colon\qff
E\dff \times\dff \bm{\Delta}\fff [\dff 1\dff]
\pff \ttoo\pff
E\dff \times\dff \bm{\Delta}\fff [\dff 1\dff]
\pff.
\]

\vspace{-9pt}\vspace{-0.5pt}
Then\sss $h$\sss is\dss a\sss translation of\dss $p_{\dff 1}$\nsp.\oss
Clearly,\oss the map\sss
$E_{\trf 0}\qff \ttoo\qff E_{\trf 0}$\sss
induced\dss by $h$
can\sss be identified\sss with\sss $f\dff(\dff c\trf)$\nnsp,\oss
and\dss the map\sss
$E_{\dff 1}\qff \ttoo\qff E_{\dff 1}$\sss
induced\dss by $h$ is\dss equal\dss to\sss the identity.\oss
Therefore,\oss the composition of\dss $h$\sss with\sss the projection\dss
$E\dff \times\dff \bm{\Delta}\fff [\dff 1\dff]
\qff \ttoo\qff
E$\sss
is\dss a\sss homotopy\dss between\sss $f\dff(\dff c\trf)$\sss
and\dss the identity.\oss
The\sss lemma\sss follows.\oss  \eproof

\myuppar{Remark.}
Since\dss
$\mathcal{Z}^{\dff n}\dff(\dff \bm{\Delta}\fff [\dff q\trf]\fff,\qff \pi_{\trf 0}\trf)
\off =\off
0$\qss
if\dss $q\qff <\qff n$\nnsp,\oss
every\dss translation\dss
$f\dff \colon\dff
E\qff \ttoo\qff E$\dss
is\dss equal\dss to\sss the identity\sss over\sss 
$\sk_{\dff n\dff -\dff 1}\fff B$\nnsp.\oss
By\dss the same reason\dss the homotopy\sss
constructed\sss in\qss Lemma\qss \ref{translations-homotopy}\qss
is\dss constant\sss over\sss 
$\sk_{\dff n\dff -\dff 2}\fff B$\nnsp.\oss\vspace{-0.125pt}

\myuppar{Remark.}
If\dss $p$\sss is\dss the\sss trivial\dss bundle\dss
$B\dff \times\dff K\dff(\dff \pi\fff,\qff n\trf)\qff \ttoo\qff B$\nnsp,\oss
then\dss the\sss translations of\dss $p$\sss correspond\dss to
maps\dss $B\qff \ttoo\qff K\dff(\dff \pi\fff,\qff n\trf)$\nnsp.\oss
So,\oss
Lemmas\qss \ref{cocycles-translations}\qss and\qss \ref{translations-homotopy}\qss
provide a\qss ``twisted''\qss version of\trs the classification of\dss
maps\dss $f\dff \colon\dff K\qff \ttoo\qff  K\dff(\dff \pi\fff,\qff n\trf)$\sss
in\dss terms of\dss cocycles\sss $z\trf(\trf f\trf)$\nnsp.\oss

\myuppar{A\sss group acting\sss on\dss $E$\nnsp.}
Let\sss 
$G\off =\off
\mathcal{C}^{\dff n\dff -\dff 1}\dff(\qff B\fff,\qff \pi\trf(\dff p\trf)\trf)$
be\sss the group of\dss normalized $(\dff n\qff -\qff 1\dff)$\dnsp-cochains of\trs $B$\sss
with\sss coefficients in\sss 
the\sss local\sss system\sss $\pi\trf(\dff p\trf)$\nnsp.\oss
By\qss Lemma\qss \ref{cocycles-translations}\qss
for every\dss 
$g\qff \in\qff
G$\dss
there exists a unique\sss automorphism\dss\vspace{3pt}
\[
\quad
a\dff(\dff c\trf)
\off =\off
f\trf(\trf \partial^{\dff *} c\trf)\dff \colon\dff 
E\qff \ttoo\qff E
\]

\vspace{-9pt}
over\trs $B$
such\dss that\dss
$D_{\fff a\dff(\dff c\trf)}\off =\off \partial^{\dff *}\dff c$\nnsp.\oss
By\qss Lemma\qss \ref{translations-composition}\qss
the map\dss
$c\off \longmapsto\off a\trf(\dff c\trf)$\dss
is\dss a\sss homomorphism.\oss
Hence\sss this map defines an action of\dss $G$\sss on\dss $E$\nnsp.\oss
By\qss Lemma\qss \ref{translations-homotopy}\qss every\sss
automorphism\dss $a\trf(\dff c\trf)$\dss is\dss homotopic\sss
to\sss the\sss identity.\oss
Moreover\halfff,\oss by\dss the\sss remark after\dss
the proof\dss of\qss Lemma\qss \ref{translations-homotopy}\qss
the homotopy\sss can\sss be chosen\dss to be constant\sss over\sss 
$\sk_{\dff n\dff -\dff 2}\fff B$\nnsp.\oss 
Therefore,\oss the group\sss $G$\sss 
acts on\sss $E$\sss
by\sss automorphisms homotopic\sss
to\sss the\sss identity\dss by\sss homotopies constant\sss over\sss 
$\sk_{\dff n\dff -\dff 2}\fff B$\nnsp.\oss

\myuppar{Free simplices.}
Let\dss us\dss say\dss  
that\sss a $q$\dnsp-simplex $\sigma$\sss is\qss 
\emph{free in dimension}\qss $m$\dss
if\trs the restriction of\trs the simplicial\dss map\dss
$i_{\trf \sigma}\dff \colon\dff
\bm{\Delta}\fff [\fff q\dff]\qff \ttoo\qff K$\dss
to\dss $\sk_{\dff m}\fff \bm{\Delta}\fff [\fff q\dff]$\sss
is\dss an\sss isomorphism onto\sss its\sss image.\oss

\mypar{Lemma.}{transitivity-bundle}
\emph{Suppose\sss that\dss $\tau\fff,\pff \tau'$\sss
are $q$\dnsp-simplices of\qss $E$\dss
such\dss that\dss
$p\dff(\dff \tau\trf)\off =\off p\dff(\dff \tau'\trf)$\nnsp.\oss
If\qss $p\dff(\dff \tau\trf)$\sss
is\qss free\sss in\sss dimension $n\qff -\qff 1$\nnsp,\oss
then\dss there exists\dss $c\qff \in\dff G$\sss such\dss that\trs
$a\trf(\dff c\trf)\dff(\dff \tau\trf)\off =\off \tau'$\nnsp.\oss}

\proof
Let\sss $\sigma\off =\off p\dff(\dff \tau\trf)$\sss
and $v$ be\sss the leading\sss vertex of\dss $\sigma$\nnsp.\oss
Let\dss\vspace{3pt}\vspace{-0.28pt}
\[
\quad 
t\dff \colon\dff
\bm{\Delta}\fff [\dff q\trf]\dff \times\dff K\dff(\dff \pi_{\fff v}\fff,\qff n\trf)
\qff \ttoo\qff
i_{\dff \sigma}^{\dff *}\qff E
\]

\vspace{-9pt}\vspace{-0.28pt}
be a special\dss trivialization of\trs the pull-back\dss bundle\sss
$i_{\dff \sigma}^{\dff *}\dff p$\nnsp.\oss
Then\vspace{3pt}\vspace{-0.28pt}
\[
\quad
t^{\dff -\dff 1}\dff(\trf \tau\trf)
\off =\off
\left(\qff \bm{\iota}_{\dff q}\fff,\qff z\qff\right)
\quad\
\mbox{and}\quad\
t^{\dff -\dff 1}\dff(\trf \tau'\trf)
\off =\off
\left(\qff \bm{\iota}_{\dff q}\fff,\qff z'\qff\right)
\pff.
\]

\vspace{-9pt}\vspace{-0.28pt}
for some $q$\dnsp-simplices\dss $z\fff,\pff z'$\dss of\trs 
$K\dff(\dff \pi_{\fff v}\dff,\qff n\trf)$\nnsp,\oss
i.e.\qss for some\dss
$z\fff,\pff z'
\qff \in\pff 
\mathcal{Z}^{\dff n}\dff(\dff \bm{\Delta}\fff [\dff q\trf]\fff,\qff \pi_{\fff v}\trf)$\nnsp.\oss
Clearly,\vspace{3pt}\vspace{-0.28pt}
\[
\quad
\bm{t}\trf(\dff z'\qff -\qff z\trf)\dff
\left(\qff \bm{\iota}_{\dff q}\fff,\qff z\qff\right)
\off =\off
\left(\qff \bm{\iota}_{\dff q}\fff,\qff z'\qff\right)
\pff.
\]

\vspace{-9pt}\vspace{-0.28pt}
Since\sss the cohomology\sss of\dss 
$\bm{\Delta}\fff [\dff q\trf]$\sss
vanish,\oss
the cocycle\dss 
$z'\qff -\qff z$\dss
is\dss the coboundary\sss of\dss some normalized
$(\dff n\qff -\qff 1\dff)$\dnsp-cochain\sss $d$\nnsp.\oss
Since\sss $\sigma$\sss
is\dss free in dimension\sss $n\qff -\qff 1$\nnsp,\oss
there exists a normalized $(\dff n\qff -\qff 1\dff)$\dnsp-cochain\sss $c$
of\dss $B$\dss such\dss that\dss
$d\off =\off i_{\trf \sigma}^{\dff *}\dff(\dff c\trf)$\qss
({\fff}the values of\dss $c$\sss on\sss non-degenerate $(\dff n\qff -\qff 1\dff)$\dnsp-simplices
not\dss belonging\dss to\sss the\sss image\sss of\dss $i_{\trf \sigma}$ are arbitrary)\qss
and\dss hence\dss\vspace{3pt}\vspace{-0.28pt}
\[
\quad
z'\qff -\qff z
\off =\off
\partial^{\dff *} d
\off =\off
i_{\trf \sigma}^{\dff *}\dff(\qff \partial^{\dff *} c\trf)
\pff.
\]

\vspace{-9pt}\vspace{-0.28pt}
Let\dss 
$f
\off =\off
a\dff(\dff c\trf)
\off =\off
f\trf(\qff \partial^{\dff *} c\trf)$\nnsp.\oss
Then\dss
$D_{\fff f}\off =\off \partial^{\dff *} c$\dss
and\dss hence\dss
$d\trf(\trf f_{\dff \sigma}\trf)
\off =\off
z'\qff -\qff z$\nnsp.\oss
It\dss follows\dss that\dss\vspace{3pt}\vspace{-0.28pt}
\[
\quad
a\trf(\dff c\trf)\dff(\dff \tau\trf)
\off =\off 
f\dff(\dff \tau\trf)
\off =\off\dff 
\tau'
\pff.
\]

\vspace{-9pt}\vspace{-0.28pt}
This completes\sss the proof\halfff.\oss  \eproof

\newpage
\mysection{Unraveling\qss simplicial\qss sets}{unraveling}

\myuppar{The unraveling.}
Let\dss 
$\Gamma\off =\off \Delta\fff [\dff \infty\dff]$\nnsp.\oss
As usual,\pss we will\sss denote by\dss $\Gamma_n$\sss the
set\sss of\dss $n$\dnsp-simplices of\trs $\Gamma$\nnsp.\oss
The\qss \emph{unravelling}\pss of\dss a simplicial\sss set\sss $K$\sss 
is\dss the dimension-wise product\dss
$\Delta\fff K\dff \times\dff \Gamma$\dnsp.\oss
We will\sss denote\sss this $\Delta$\dnsp-set\sss simply\sss by\sss
$K\dff \times\dff \Gamma$\nnsp.\oss
The goal\sss of\trs this section\dss is\dss to prove\sss that\dss the projection\dss
$p\dff \colon\dff
K\dff \times\dff \Gamma
\qff \ttoo\qff
K$\sss
induces isomorphisms in\sss bounded cohomology.\oss\vspace{0.6pt}

\myuppar{Averaging\sss operators.}
An\qss \emph{averaging\sss operator\sss on}\pss $\Gamma_n$\dss is\dss a\sss 
bounded\dss linear\sss functional\vspace{3pt}\vspace{0.6pt}
\[
\quad
m_{\dff n}\dff \colon\dff
B\dff(\trf \Gamma_n\trf)
\qff \ttoo\qff
\rrr
\pff,
\]

\vspace{-9pt}\vspace{0.6pt}
of\trs the norm\dss $1$\dss
equal\dss to\sss the identity\sss on\sss constant\dss functions.\oss
More precisely,\oss if\dss $f\dff(\dff n\trf)\off =\off a$\sss
for all\dss $n\qff \in\pff \nnn$\nnsp,\oss
then\sss it\dss is\dss required\dss that\dss
$m_{\dff n}\dff(\trf f\trf)\off =\off a$\nnsp.\oss
A family\dss of\dss averaging operators\dss 
$m_{\dff n}$\nsp,\oss 
where\sss $n\qff \in\pff \nnn$\nnsp,\oss 
is\dss said\dss to be\qss \emph{coherent}\oss if\dss
the operators\dss $m_{\dff n}$\dss commute with\dss the adjoints 
of\trs the face operators\dss
$\partial_{\dff i}\dff \colon\dff
\Gamma_n\qff \ttoo\qff \Gamma_{n\dff -\dff 1}$\nsp,\dff\oss
i.e.\pss if\trs\vspace{3pt}\vspace{0.6pt}
\[
\quad
m_{\dff n}\dff \circ\qff \partial_{\dff i}^{\dff *}
\off =\off\dff
m_{\dff n\dff -\dff 1}
\]

\vspace{-9pt}\vspace{0.6pt}
for every\sss 
$n\qff \geq\qff 1$\dss 
and\dss $i\qff \in\pff \nnn$\nnsp.\oss
Such a\sss family\sss defines a graded\dss map of\dss degree $0$\vspace{3pt}\vspace{0.6pt}
\[
\quad
m_{\dff *}\dff \colon\dff
B^{\dff *}\fff(\trf K\dff \times\dff \Gamma\trf)
\qff \ttoo\qff
B^{\dff *}\fff(\trf K\trf)
\pff
\]

\vspace{-9pt}\vspace{0.6pt}
by\sss averaging\sss cochains over\sss preimages\sss in\sss 
$K\dff \times\dff \Gamma$\sss
of\dss simplices of\trs $K$\nnsp.\oss
In\dss fact\halfff,\pss
$m_{\dff *}$\sss is\dss a cochain\dss map.\oss
See\qss Lemma\qss \ref{averaging-cochain}.\oss
Clearly,\oss
$m_{\dff *}\dff \circ\qff p^{\dff *}\off =\off \id$\nnsp.\oss\vspace{0.6pt}

\myuppar{Banach\dss limits.}
Given a function\sss $f\dff \colon\dff \nnn\qff \ttoo\qff \rrr$\nnsp,\oss
let\dss $s f$\sss be\sss the function\dss
$\nnn\qff \ttoo\qff \rrr$\dss defined\dss by\dss
$s f\dff(\dff n\trf)\off =\off f\dff(\dff n\qff +\qff 1\dff)$\nnsp.\oss
A\qss \emph{Banach\dss limit}\qss
is\dss a\sss linear\dss functional\dss
$l\dff \colon\dff
B\dff(\trf \nnn\trf)\qff \ttoo\qff \rrr$\dss
such\dss that\dss its norm\dss is\dss equal\dss to $1$\nnsp,\oss
$l\trf(\trf f\trf)\off =\off a$\dss if\dss $f$\sss is\dss the constant\sss
function\sss with\sss the value $a\qff \in\pff \rrr$\nnsp,\oss
and\dss
$l\trf(\trf sf\trf)\off =\off l\trf(\trf f\trf)$\sss for all\dss $f$\dnsp.\oss
It\dss is\dss well\dss known\dss that\trs Banach\dss limits exist.\oss
See,\oss for example,\oss \cite{r},\oss Exercise\qss 4\qss to\dss Chapter\qss 3.\oss
Let\dss us\dss fix a\dss Banach\dss limit\sss and\dss denote it\dss by\dss $\lim$\nnsp.\oss

Suppose now\dss that\sss $f\dff(\dff n\trf)$\sss 
is\dss a\sss bounded\dss real-valued\dss function of\trs
the natural\sss argument\sss $n$\sss defined only\sss for sufficiently\dss large
numbers $n$\nnsp.\oss
If\trs $N$\sss is\dss sufficiently\dss large,\oss
then\sss the function\sss 
$f_{\qff N}\dff(\dff n\trf)\off =\off f\dff(\dff n\qff +\qff N\trf)$\dss
is\dss defined\dss for all\dss $n\qff \in\pff \nnn$\nnsp.\oss
Clearly,\pss $f_{\qff N\qff +\qff 1}\off =\off s f_{\qff N}$\nsp.\oss
Therefore\sss we can define\dss
$\lim\dff f$\sss as\sss the common values of\trs $\lim\dff f_{\qff N}$\sss
for sufficiently\dss large natural\dss numbers\sss $N$\nnsp.\oss

Suppose now\sss that\dss
$f\dff(\trf a\fff,\pff b\fff,\pff \ldots\fff,\pff z\trf)$\sss
is\dss a\sss bounded\dss real-valued\dss function of\trs
several\dss natural\sss variable\sss
$a\fff,\pff b\fff,\pff \ldots\fff,\pff z$\sss
and\dss that\sss $k$\sss is\dss one of\trs these variables.\oss
By\sss fixing\sss values of\dss other variables and applying\sss $\lim$\sss
to\dss the resulting\sss function of\dss $k$\sss we will\sss get\sss
a bounded\sss function of\trs the other variables\sss
$a\fff,\pff \ldots\fff,\pff \widehat{k}\off \ldots\fff,\pff  z$,\oss
which\sss we will\sss denote by\sss 
$\lim_{\dff k}\dff f\dff
(\trf a\fff,\pff \ldots\fff,\pff \widehat{k}\off \ldots\fff,\pff  z\trf)$\nnsp.\oss
As above,\oss this operation applies even\trs if\dss
$f\dff(\trf a\fff,\pff b\fff,\pff \ldots\fff,\pff z\trf)$\sss
is\dss defined only\dss for sufficiently\dss large\sss $k$\nnsp.\oss

\mypar{Lemma.}{coherence}
\emph{Coherent\dss families of\dss averaging operators exist\halfff.}

\proof
The $n$\dnsp-simplices of\sss $\Gamma\nsp$ can\sss be identified\sss
with\dss the sequences
$(\qff k_{\trf 0}\fff,\pff k_{\dff 1}\fff,\pff \ldots\fff,\pff k_{\dff n}\trf)
\trf \in\trf
\nnn^{\dff n\dff +\dff 1}$\dss
such\dss that\dss
$k_{\trf 0}\qff <\qff k_{\dff 1}\qff <\qff \ldots\qff <\qff k_{\dff n}$\nsp.\oss
Given a bounded\dss function\dss
$f\dff \colon\dff
\Gamma_n\qff \ttoo\qff \rrr$\nnsp,\oss
let\sss $f^{\trf (\dff 1\dff)}$\sss be\sss the function\dss
$\Gamma_{n\dff -\dff 1}\qff \ttoo\qff \rrr$\sss
defined\dss by\vspace{3pt}
\[
\quad
f^{\trf (\dff 1\dff)}\dff
(\qff k_{\trf 0}\fff,\pff k_{\dff 1}\fff,\pff \ldots\fff,\pff k_{\dff n\dff -\dff 1}\trf)
\off =\off
\lim\nolimits_{\qff k_{\dff n}}\qff
f\dff(\qff k_{\trf 0}\fff,\pff k_{\dff 1}\fff,\pff \ldots\fff,\pff k_{\dff n}\trf)
\pff.
\]

\vspace{-9pt}
For\dss $0\qff \leq\qff m\qff \leq\qff n\qff +\qff 1$\dss let\dss
us\dss define\dss $f^{\trf (\dff m\trf)}$\sss recursively\dss by\dss
$f^{\trf (\dff 0\dff)}\off =\off f$\sss and\vspace{3pt}
\[
\quad
f^{\trf (\dff m\dff +\dff 1\trf)}
\off =\off
\left(\qff f^{\trf (\dff m\trf)}\qff\right)^{\trf (\dff 1\dff)}
\pff.
\]

\vspace{-9pt}
Then\dss $f^{\trf (\dff m\trf)}$\sss is\dss a\sss function of\dss 
$n\qff +\qff 1\qff -\qff m$\sss natural\dss variables.\oss
In\dss particular\halfff,\oss
$f^{\trf (\dff n\dff +\dff 1\dff)}$\sss is\dss a\sss function of\dss zero variables,\oss
i.e.\qss is\dss a constant.\oss
Let\dss $m_{\dff n}\dff(\trf f\trf)$\sss be\sss this constant.\oss
Clearly,\oss each $m_{\dff n}$\sss is\dss an averaging operator.\oss
We claim\dss that\dss the family\sss of\trs these operators\dss is\dss coherent,\oss
i.e.\qss that\vspace{3pt}
\begin{equation}
\label{iterated-limits}
\quad
\left(\qff \partial_{\dff i}^{\dff *}\dff f\qff\right)^{\trf (\dff n\dff +\dff 2\dff)}
\off =\off
f^{\trf (\dff n\dff +\dff 1\dff)}
\end{equation}

\vspace{-9pt}
for every\dss $i\qff \in\qff [\fff n\qff +\qff 1\dff]$\nnsp.\oss
By\dss the definition,\oss\vspace{3pt}
\[
\quad
\partial_{\dff i}^{\dff *}\dff f\dff
(\qff k_{\trf 0}\fff,\pff k_{\dff 1}\fff,\pff \ldots\fff,\pff k_{\dff n\dff +\dff 1}\trf)
\off =\off
f\dff
(\qff k_{\trf 0}\fff,\pff \ldots\fff,\off
\widehat{k_{\dff i}}\off \ldots\fff,\pff 
k_{\dff n\dff +\dff 1}\trf)
\]

\vspace{-9pt}
is\dss a function\sss independent\sss of\dss $k_{\dff i}$\nsp.\oss
By\dss consecutively\dss taking\dss limits we see\sss that\vspace{3pt}\vspace{1pt}
\[
\quad
\left(\qff \partial_{\dff i}^{\dff *}\dff f\qff\right)^{\trf (\dff m\dff)}\dff
(\qff k_{\trf 0}\fff,\pff k_{\dff 1}\fff,\pff \ldots\fff,\pff k_{\dff n\dff +\dff 1\dff -\dff m}\trf)
\off =\off
f^{\trf (\dff m\dff)}\dff
(\qff k_{\trf 0}\fff,\pff \ldots\fff,\off
\widehat{k_{\dff i}}\off \ldots\fff,\pff 
k_{\dff n\dff +\dff 1\dff -\dff m}\trf)
\]

\vspace{-9pt}\vspace{1pt}
for\dss $n\qff +\qff 1\qff -\qff m\qff \geq\qff i$\nnsp,\oss
i.e.\qss for\dss
$m\qff \leq\qff n\qff +\qff 1\qff -\qff i$\nnsp.\oss
In\dss particular\halfff,\oss\vspace{3pt}\vspace{1pt}
\[
\quad
\left(\qff \partial_{\dff i}^{\dff *}\dff f\qff\right)^{\trf (\dff n\dff +\dff 1\dff -\dff i\trf)}\dff
(\qff k_{\trf 0}\fff,\pff k_{\dff 1}\fff,\pff \ldots\fff,\pff k_{\dff i}\trf)
\off =\off
f^{\trf (\dff n\dff +\dff 1\dff -\dff i\trf)}\dff
(\qff k_{\trf 0}\fff,\pff \ldots\fff,\off
k_{\dff i\dff -\dff 1}\trf)
\pff.
\]

\vspace{-9pt}\vspace{1pt}
By\dss taking\dss the\sss limit\sss of\trs the\sss left\dss hand side,\oss
which\dss is\dss independent\sss of\dss $k_{\dff i}$\nsp,\oss
we see\sss that\vspace{3pt}\vspace{1pt}
\[
\quad
\left(\qff \partial_{\dff i}^{\dff *}\dff f\qff\right)^{\trf (\dff n\dff +\qff 2\dff -\dff i\trf)}\dff
(\qff k_{\trf 0}\fff,\pff k_{\dff 1}\fff,\pff \ldots\fff,\pff k_{\dff i\dff -\dff 1}\trf)
\off =\off
f^{\trf (\dff n\dff +\dff 1\dff -\dff i\trf)}\dff
(\qff k_{\trf 0}\fff,\pff \ldots\fff,\off
k_{\dff i\dff -\dff 1}\trf)
\pff.
\]

\vspace{-9pt}\vspace{1pt}
Taking\dss the\sss limits\sss $i$\sss more\sss times shows\sss that\dss the equality\qss
(\ref{iterated-limits})\qss holds.\oss
Therefore our\sss family\sss of\dss averaging\sss operators\dss is\dss
indeed\sss coherent.\oss  \eproof

\myuppar{Acyclicity\sss of\dss 
$\bm{\Delta}\dff[\halfff n\dff]\dff \times\dff \Gamma$\dnsp.}
Recall\dss that\sss an $m$\dnsp-chain of\dss
a $\Delta$\dnsp-set\sss $D$\sss is\dss
a finite formal\sss sum of\dss $m$\dnsp-simplices of\dss $D$\sss
with coefficients in some abelian\sss group.\oss
A\qss \emph{vertex}\pss of\dss a chain\sss is\dss defined as a vertex of\dss
some simplex entering\sss into\sss this sum\sss with\dss non-zero coefficient.\oss 
If\dss an $m$\dnsp-chain $c$ of\dss $\bm{\Delta}\dff[\halfff n\dff]\dff \times\dff \Gamma$\sss 
is\dss a cycle,\oss 
then\sss $c$\sss is\dss a\sss boundary\sss in\dss 
$\bm{\Delta}\dff[\halfff n\dff]\dff \times\dff \Gamma$\dnsp.\oss
Indeed,\oss since $c$\sss is\dss a finite sum,\oss
there exists\dss $m\qff \in\pff \nnn$\sss
such\dss that\dss for every\sss vertex\sss $(\dff v\fff,\qff k\trf)$\dss
of\sss $c$\sss the inequality\dss $k\qff <\qff m$\sss holds.\oss
Let\sss $w$\sss be a vertex of\dss $\bm{\Delta}\dff[\halfff n\dff]$\nnsp,\oss
and\dss let\dss us consider\dss the cone\sss $b$\sss over $c$ with\dss 
the apex\dss $(\dff w\fff,\pff m\trf)$\nnsp.\oss
In\sss order\dss to ensure\sss that\dss this cone\dss is\dss indeed
a chain of\dss 
$\bm{\Delta}\dff[\halfff n\dff]\dff \times\dff \Gamma$\sss
one needs\sss to\sss build\dss the cone by\sss adding\dss the apex
as\sss the\dss last\sss vertex of\dss every\dss simplex of\dss $c$\nnsp.\oss
Then\vspace{3pt}
\[
\quad
c
\off =\off 
(\qff -\qff 1\trf)^{\dff m\dff +\dff 1}\qff 
\partial\dff b
\pff.
\]

\vspace{-9pt}
The sign\dss is\dss caused\dss by\sss adding\dss the apex as\sss the\dss last\dss vertex.\oss

\myuppar{The method of\dss acyclic models.}
For a simplicial\sss or $\Delta$\dnsp-set\dss $K$\dss
let\dss $C_{\dff *}\dff(\trf K\trf)$\sss
be\sss the complex of\dss chains\sss in\sss $K$\sss 
with coefficients in some abelian\dss group.\oss
The method of\dss acyclic models applied\dss to\sss the functors
$K\off \longmapsto\off C_{\dff *}\dff(\trf K\trf)$
and
$K\off \longmapsto\off C_{\dff *}\dff(\trf K\dff \times\dff \Gamma\trf)$
from simplicial\sss sets\sss to chain complexes\sss
implies\sss that\dss
$p_{\dff *}\dff \colon\dff
C_{\dff *}\dff(\trf K\dff \times\dff \Gamma\trf)
\qff \ttoo\qff
C_{\dff *}\dff(\trf K\trf)$\sss
is\dss a chain\sss homotopy\sss equivalence.\oss

We will\sss adapt\dss the method of\dss acyclic models\sss
to prove\sss that\dss
$p^{\dff *}\dff \colon\dff
B^{\fff *}\dff(\trf K\trf)
\qff \ttoo\qff
B^{\fff *}\dff(\trf K\dff \times\dff \Gamma\trf)$\sss
is\dss a cochain\sss homotopy\sss equivalence.\oss
Let\sss $m_{\dff n}$\sss be a coherent\dss family\sss
of\dss averaging operators.\oss
Since\sss $m_{\dff *}\dff \circ\qff p^{\dff *}\off =\off \id$\nnsp,\oss 
it\dss is\dss sufficient\dss to prove\sss that
$p^{\dff *}\dff \circ\qff m_{\dff *}
\dff \colon\dff
B^{\fff *}\dff(\trf K\dff \times\dff \Gamma\trf)
\qff \ttoo\qff
B^{\fff *}\dff(\trf K\dff \times\dff \Gamma\trf)$
is\dss cochain\sss homotopic\sss to\sss the~identity.\oss

\myuppar{Some special\sss chains.}
Recall\dss that\dss
$d\trf(\dff i\trf)\dff \colon\dff
[\halfff n\qff -\qff 1\dff]
\qff \ttoo\qff
[\halfff n\dff]$\dss
is\dss the unique strictly\dss increasing\sss map not\dss having $i$ as\sss a value.\oss
Let\dss
$\delta_{\dff i}
\off =\off
d\trf(\dff i\trf)_{\dff *}
\dff \colon\dff
\bm{\Delta}\dff[\halfff n\qff -\qff 1\dff]
\qff \ttoo\qff 
\bm{\Delta}\dff[\halfff n\dff]$\dss
be\sss the simplicial\dss map induced\dss by\sss $d\trf(\dff i\trf)$\nnsp.\oss
We will\sss use\sss the following abbreviated\dss 
notation\dss for sums\fff:\vspace{6pt}
\[
\quad
\sum\nolimits_{\qff i}'\off \bullet
\off =\off
\sum\nolimits_{\qff i}\qff (\trf -\qff 1\dff)^{\dff i}\off \bullet
\off\off.
\]

\vspace{-6pt}
For\sss every\dss two\sss simplices\dss
$\tau\fff,\pff \tau'\qff \in\pff \Gamma_n$\dss
we are going\dss to define an
$(\dff n\qff +\qff 1\dff)$\dnsp-chain\sss
$c_{\dff n}\dff(\dff \tau\fff,\qff \tau'\trf)$\sss
of\dss $\bm{\Delta}\dff[\halfff n\dff]\dff \times\dff \Gamma$\sss
with\sss integer coefficients in such a way\dss that\vspace{7.5pt}\vspace{0.875pt}
\begin{equation}
\label{cn-boundary}
\quad
\partial\dff
c_{\dff n}\dff(\dff \tau\fff,\qff \tau'\trf)
\off =\off
(\trf \bm{\iota}_{\dff n}\fff,\qff \tau\trf)
\qff -\qff 
(\trf \bm{\iota}_{\dff n}\fff,\qff \tau'\trf)
\off -\off\qff
\sum'\nolimits_{\qff i}\pff
\bigl(\trf 
\delta_{\dff i}\dff \times\trf \id_{\qff \Gamma} 
\trf\bigr)_{\dff *}\trf
\left(\trf 
c_{\dff n\dff -\dff 1}\dff
\left(\qff 
\partial_{\dff i}\dff \tau\fff,\pff \partial_{\dff i}\dff \tau'
\qff\right)
\qff\right)
\pff.
\end{equation}

\vspace{-4.5pt}\vspace{0.875pt}
In addition,\oss we will\dss require\sss that\dss the 
$l_{\dff 1}$\dnsp-norm of\dss $c_{\dff n}\dff(\dff \tau\fff,\qff \tau'\trf)$\qss
({\halfff}i.e.\qss the sum of\trs the absolute values of\trs the coefficients)\qss
can\sss be bounded\dss by constants
depending\sss only\sss on $n$\nnsp.\oss
For\dss $n\off =\off 0$\dss the condition\qss (\ref{cn-boundary})\qss 
simplifies\sss to\vspace{4.5pt}
\[
\quad
\partial\dff
c_{\trf 0}\dff(\dff \tau\fff,\qff \tau'\trf)
\off =\off
(\trf \bm{\iota}_{\trf 0}\fff,\qff \tau\trf)
\qff -\qff 
(\trf \bm{\iota}_{\trf 0}\fff,\qff \tau'\trf)
\pff.
\]

\vspace{-7.5pt}
We will\sss construct\sss such chains using a recursion\sss by\sss $n$\nnsp.\oss
The chain\dss
$(\trf \bm{\iota}_{\trf 0}\fff,\qff \tau\trf)
\qff -\qff 
(\trf \bm{\iota}_{\trf 0}\fff,\qff \tau'\trf)$\dss
has\sss the augmentation\sss $0$\sss and\dss is\dss a boundary\sss in\dss
$\bm{\Delta}\dff[\dff 0\dff]\dff \times\dff \Gamma$\dss
if\qss $N\qff \geq 2$\nnsp.\oss
Assuming\dss that\trs $N\qff \geq\qff 2$\nnsp,\oss
let\dss us choose a vertex $v$ of\trs $\Gamma$\sss strictly\dss larger\dss that\dss
$\tau\fff,\pff \tau'$\sss in\dss the natural\sss order\qss
(recall\dss that\dss $\Gamma_0\off =\off \nnn$\nsp)\qss 
and\dss take as\sss 
$c_{\trf 0}\dff(\dff \tau\fff,\qff \tau'\trf)$\sss
the cone with\dss the apex\sss $(\trf \bm{\iota}_{\trf 0}\fff,\qff v\trf)$\dss
over\dss the cycle\dss
$(\trf \bm{\iota}_{\trf 0}\fff,\qff \tau\trf)
\qff -\qff 
(\trf \bm{\iota}_{\trf 0}\fff,\qff \tau'\trf)$\nnsp.\oss
Then\sss the
$l_{\dff 1}$\dnsp-norm of\dss $c_{\dff n}\dff(\dff \tau\fff,\qff \tau'\trf)$\sss
is\dss $\leq\qff 2$\nnsp.\oss

Suppose\sss that\dss the chains\sss
$c_{\dff m}\dff(\dff \tau\fff,\qff \tau'\trf)$\sss
are already\sss defined\dss for\dss $m\qff \leq\qff n\qff -\qff 1$\nnsp,\oss
the condition\qss (\ref{cn-boundary})\qss 
holds for\dss them,\oss 
and\dss there are\sss required\dss bounds
on\dss the $l_{\dff 1}$\dnsp-norms.\oss
In order\dss to define\sss the chains\sss
$c_{\dff n}\dff(\dff \tau\fff,\qff \tau'\trf)$\sss
we need\dss to verify\dss that\dss the right\dss hand side of\qss
(\ref{cn-boundary})\qss is\dss a cycle.\oss
The boundary\sss of\trs the right\dss hand side\dss is\vspace{5.5pt}
\[
\quad
\partial\dff (\trf \bm{\iota}_{\dff n}\fff,\qff \tau\trf)
\qff -\qff 
\partial\dff (\trf \bm{\iota}_{\dff n}\fff,\qff \tau'\trf)
\off -\off\qff
\sum'\nolimits_{\qff i}\pff
\partial\qff \bigl(\trf 
\delta_{\dff i}\dff \times\trf \id_{\qff \Gamma} 
\trf\bigr)_{\dff *}\trf
\left(\trf 
c_{\dff n\dff -\dff 1}\dff
\left(\qff 
\partial_{\dff i}\dff \tau\fff,\pff \partial_{\dff i}\dff \tau'
\qff\right)
\qff\right)
\]

\vspace{-30pt}
\[
\quad
=\off
\sum'\nolimits_{\qff i}\pff
\bigl(\qff 
\partial_{\dff i}\dff \bm{\iota}_{\dff n}\fff,\qff \partial_{\dff i}\dff \tau
\qff\bigr)
\off -\off 
\left(\qff 
\partial_{\dff i}\dff \bm{\iota}_{\dff n}\fff,\qff \partial_{\dff i}\dff \tau'
\qff\right)
\off -\off
\bigl(\trf 
\delta_{\dff i}\dff \times\trf \id_{\qff \Gamma} 
\trf\bigr)_{\dff *}\trf
\left(\qff
\partial\dff  
c_{\dff n\dff -\dff 1}\dff
\left(\qff 
\partial_{\dff i}\dff \tau\fff,\pff \partial_{\dff i}\dff \tau'
\qff\right)
\qff\right)
\]

\vspace{-30pt}
\[
\quad
=\off
\sum'\nolimits_{\qff i}\pff
\bigl(\trf 
\delta_{\dff i}\dff \times\trf \id_{\qff \Gamma} 
\trf\bigr)_{\dff *}\trf
\Bigl(\qff
\bigl(\qff 
\bm{\iota}_{\dff n\dff -\dff 1}\fff,\qff \partial_{\dff i}\dff \tau
\qff\bigr)
\pff -\pff
\bigl(\qff 
\bm{\iota}_{\dff n\dff -\dff 1}\fff,\qff \partial_{\dff i}\dff \tau'
\qff\bigr)
\pff -\pff
\partial\dff  
c_{\dff n\dff -\dff 1}\dff
\left(\qff 
\partial_{\dff i}\dff \tau\fff,\pff \partial_{\dff i}\dff \tau'
\qff\right)
\qff\Bigr)
\off.
\]

\vspace{-6pt}
By\sss applying\qss (\ref{cn-boundary})\qss with\dss $n\qff -\qff 1$\sss
in\dss the role of\dss $n$\sss and cancelling\dss two occurrences of\vspace{3pt}
\[
\quad
\bigl(\qff
\bm{\iota}_{\dff n\dff -\dff 1}\fff,\qff \partial_{\dff i}\dff \tau
\qff\bigr)
\pff -\pff
\bigl(\qff 
\bm{\iota}_{\dff n\dff -\dff 1}\fff,\qff \partial_{\dff i}\dff \tau'
\qff\bigr)
\]

\vspace{-9pt}
we conclude\sss that\dss the boundary of\trs 
the right\dss hand side of\qss (\ref{cn-boundary})\qss is\dss equal\dss to\vspace{6pt}
\[
\quad
\sum'\nolimits_{\qff i}\pff
\bigl(\trf 
\delta_{\dff i}\dff \times\trf \id_{\qff \Gamma} 
\trf\bigr)_{\dff *}\trf
\left(\qff
\sum'\nolimits_{\qff k}\pff
\bigl(\trf 
\delta_{\dff k}\dff \times\trf \id_{\qff \Gamma} 
\trf\bigr)_{\dff *}\trf
\left(\trf 
c_{\dff n\dff -\dff 1}\dff
\left(\qff 
\partial_{\dff k}\trf \partial_{\dff i}\dff \tau\fff,\pff 
\partial_{\dff k}\trf \partial_{\dff i}\dff \tau'
\qff\right)
\qff\right)
\qff\right)
\]

\vspace{-30pt}
\[
\quad
=\off
\sum'\nolimits_{\qff i}\pff
\sum'\nolimits_{\qff k}\pff
\bigl(\trf 
\delta_{\dff i}\dff \times\trf \id_{\qff \Gamma} 
\trf\bigr)_{\dff *}
\trf \circ\trf
\bigl(\trf 
\delta_{\dff k}\dff \times\trf \id_{\qff \Gamma} 
\trf\bigr)_{\dff *}\trf
\left(\trf 
c_{\dff n\dff -\dff 1}\dff
\left(\qff 
\partial_{\dff k}\trf \partial_{\dff i}\dff \tau\fff,\pff 
\partial_{\dff k}\trf \partial_{\dff i}\dff \tau'
\qff\right)
\qff\right)
\]

\vspace{-30pt}
\[
\quad
=\off
\sum'\nolimits_{\qff i}\pff
\sum'\nolimits_{\qff k}\pff
\bigl(\trf 
\delta_{\dff i}\dff \circ\dff \delta_{\dff k}\trf \times\trf \id_{\qff \Gamma} 
\trf\bigr)_{\dff *}\trf
\left(\trf 
c_{\dff n\dff -\dff 1}\dff
\left(\qff 
\partial_{\dff k}\trf \partial_{\dff i}\dff \tau\fff,\pff 
\partial_{\dff k}\trf \partial_{\dff i}\dff \tau'
\qff\right)
\qff\right)
\off.
\]

\vspace{-6pt}
As in\dss the proof\dss of\trs the identity\dss 
$\partial\dff \circ\dff \partial\off =\off 0$\nnsp,\oss
all\sss summands\sss in\dss the\sss last\sss double sum cancel\qss
(recall\dss that\dss $\sum'$\sss denotes an alternating sum).\oss
It\dss follows\dss that\dss the right\dss hand side 
of\qss (\ref{cn-boundary})\qss is\dss a cycle.\oss
Clearly,\oss
the $l_{\dff 1}$\dnsp-norm of\trs
the right\dss hand side can be bounded\sss in\sss terms of\sss $n$ and\dss
the $l_{\dff 1}$\dnsp-norm of\sss
$c_{\dff n\dff -\dff 1}\dff(\qff \rho\fff,\qff \rho'\trf)$\nnsp.\oss
By\dss the inductive assumption\dss this implies\sss that\dss these norms
can\sss bounded\sss in\sss terms of\dss $n$\sss only.\oss
Hence one can\dss take as\sss $c_{\dff n}\dff(\dff \tau\fff,\qff \tau'\trf)$\sss
the cone over\dss the right\dss hand side with an appropriate apex.\oss
Then\dss the $l_{\dff 1}$\dnsp-norms of\sss
$c_{\dff n}\dff(\dff \tau\fff,\qff \tau'\trf)$\sss
and of\trs
the right\dss hand side are equal\sss and can\dss 
be bounded\sss in\dss terms of\dss $n$\nnsp.\oss
This completes\sss the construction of\trs the chains\sss
$c_{\dff n}\dff(\dff \tau\fff,\qff \tau'\trf)$\nnsp.\oss

\myuppar{Partial\sss averaging.}
We will\sss deal\dss with\dss the functions of\dss several\sss variables such as bounded\sss
cochains\dss
$(\dff \sigma\fff,\qff \tau\fff,\qff \tau'\trf)
\off \longmapsto\off
g\trf(\dff \sigma\fff,\qff \tau\fff,\qff \tau'\trf)$\dss
and apply\dss the averaging operators\sss to only\sss one of\trs the variables.\oss
If\halfff,\oss say,\pss $\tau'$\sss runs over $\Gamma_n$\nsp,\oss
we will\sss denote by\vspace{1.5pt}
\[
\quad
(\dff \sigma\fff,\qff \tau\trf)
\off \longmapsto\off
m_{\dff n}\dff\sco{\tau'\fff}\trf g\trf(\dff \sigma\fff,\qff \tau\fff,\qff \tau'\trf)
\pff
\]

\vspace{-10.5pt}
the function of\dss variables $(\dff \sigma\fff,\qff \tau\trf)$
resulting\sss from applying $m_{\dff n}$\sss to functions\vspace{1.5pt}
\[
\quad
\tau'
\pff \longmapsto\pff
g\trf(\dff \sigma\fff,\qff \tau\fff,\qff \tau'\trf)
\pff.
\]

\vspace{-10.5pt}
With\dss these notations\sss the coherence condition\dss
takes\sss the form\vspace{4.5pt}
\[
\quad
m_{\dff n}\dff\sco{\tau'\fff}\trf 
g\trf(\dff \sigma\fff,\qff \tau\fff,\qff \partial_{\dff i}\dff \tau'\trf)
\off =\off
m_{\dff n\dff -\dff 1}\dff\sco{\fff\rho}\trf 
g\trf(\dff \sigma\fff,\qff \tau\fff,\pff \rho\trf)
\pff,
\]

\vspace{-7.5pt}
where $\tau'$\sss runs over\sss $\Gamma_n$\sss
and\sss $\rho$\sss runs over\sss $\Gamma_{n\dff -\dff 1}$\nsp.\oss

\myuppar{Constructing\sss cochain\sss homotopies.}
Let\dss $K$\sss be a simplicial\sss set.\oss
Every $n$\dnsp-simplex of\trs $K\dff \times\dff \Gamma\dff \times\dff \Gamma$\dss
has\sss the form $(\dff \sigma\fff,\qff \tau\fff,\qff \tau'\trf)$\nnsp,\oss
where $\sigma\qff \in\pff K_{\dff n}$\dss and\dss
$\tau\fff,\pff \tau'\qff \in\pff \Gamma_n$\nsp.\oss
Clearly,\vspace{4.5pt}
\[
\quad
(\dff \sigma\fff,\qff \tau\fff,\qff \tau'\trf)
\off =\off
\bigl(\trf i_{\trf \sigma}\dff \times\trf \id_{\qff \Gamma\dff \times\dff \Gamma} \trf\bigr)_{\dff *}\trf
(\qff \bm{\iota}_{\dff n}\fff,\qff \tau\fff,\qff \tau'\trf)
\pff.
\]

\vspace{-7.5pt}
Let\dss 
$k_{\dff n}\dff \colon\dff
C_{\dff n}\dff(\trf K\dff \times\dff \Gamma\dff \times\dff \Gamma\trf)
\qff \ttoo\qff
C_{\dff n}\dff(\trf K\dff \times\dff \Gamma\trf)$\dss
be\sss the unique homomorphism such\dss that\vspace{4.5pt}
\[
\quad 
k_{\dff n}\dff(\dff \sigma\fff,\qff \tau\fff,\qff \tau'\trf)
\off =\off
\bigl(\trf i_{\trf \sigma}\dff \times\trf \id_{\qff \Gamma\dff \times\dff \Gamma} \trf\bigr)_{\dff *}\trf
\left(\trf c_{\dff n}\dff(\dff \tau\fff,\qff \tau'\trf)\trf\right)
\pff.
\]

\vspace{-7.5pt}
for every $n$\dnsp-simplex $(\dff \sigma\fff,\qff \tau\fff,\qff \tau'\trf)$\nnsp.\oss
The condition\qss (\ref{cn-boundary})\qss implies\sss that\vspace{6pt}
\begin{equation}
\label{special-homotopy}
\quad
\partial\dff
k_{\dff n}\dff(\dff \sigma\fff,\qff \tau\fff,\qff \tau'\trf)
\off =\off
(\trf \sigma\fff,\qff \tau\trf)
\qff -\qff 
(\trf \sigma\fff,\qff \tau'\trf)
\off -\off\qff
\sum'\nolimits_{\qff i}\pff
k_{\dff n\dff -\dff 1}\dff(\qff \partial_{\dff i}\dff \sigma\fff,\qff 
\partial_{\dff i}\dff \tau\fff,\qff 
\partial_{\dff i}\dff \tau'\qff)
\pff.
\end{equation} 

\vspace{-6pt}
The next\sss step\sss is\dss to apply\dss the averaging operators $m_{\dff n}$\nsp.\oss
If\dss
$f\qff \in\pff B^{\fff n\dff +\dff 1}\dff(\trf K\dff \times\dff \Gamma\trf)$\nnsp,\oss
then\vspace{4.5pt}
\[
\quad
(\dff \sigma\fff,\qff \tau\fff,\qff \tau'\trf)
\off \longmapsto\off
f\qff
\bigl(\qff
k_{\dff n}\dff(\dff \sigma\fff,\qff \tau\fff,\qff \tau'\trf)
\qff\bigr)
\]

\vspace{-7.5pt}
is\dss a bounded\dss function\dss because,\oss 
together\sss with\dss the $l_{\dff 1}$\dnsp-norm of\dss
$c_{\dff n}\dff(\dff \tau\fff,\qff \tau'\trf)$\nnsp,\oss
the $l_{\dff 1}$\dnsp-norm of\dss 
$k_{\dff n}\dff(\dff \sigma\fff,\qff \tau\fff,\qff \tau'\trf)$\sss
can\sss be bounded\sss in\dss terms of\sss $n$ only.\oss
Let\vspace{4.5pt}
\[
\quad
h_{\dff n\dff +\dff 1}\dff(\trf f\trf)\trf
(\dff \sigma\fff,\qff \tau\trf)
\off =\off
m_{\dff n}\dff\sco{\tau'\fff}\trf f\qff
\bigl(\qff
k_{\dff n}\dff(\dff \sigma\fff,\qff \tau\fff,\qff \tau'\trf)
\qff\bigr)
\pff.
\]

\vspace{-7.5pt}
Then\dss 
$h_{\dff n\dff +\dff 1}\dff(\trf f\trf)
\qff \in\pff
B^{\fff n}\dff(\trf K\dff \times\dff \Gamma\trf)$\dss
and,\oss moreover\halfff,\oss\vspace{4.5pt}
\[
\quad
h_{\dff n\dff +\dff 1}
\dff \colon\dff
B^{\fff n\dff +\dff 1}\dff(\trf K\dff \times\dff \Gamma\trf)
\qff \ttoo\qff
B^{\fff n}\dff(\trf K\dff \times\dff \Gamma\trf)
\]

\vspace{-7.5pt}
is\dss a\sss bounded operator\halfff.\oss

\mypar{Lemma.}{main-homotopy}
\emph{The operators\sss $h_{\dff n}$\qss form\sss a cochain\dss
homotopy\dss between $p^{\dff *}\dff \circ\qff m^{\dff *}$
and\dss the identity.}

\proof
Let\dss 
$f\qff \in\pff B^{\fff n}\dff(\trf K\dff \times\dff \Gamma\trf)$\nnsp.\oss
By\sss appling\sss $f$\sss 
to\qss (\ref{special-homotopy})\qss we see\sss that\vspace{6pt}
\[
\quad
f\qff
\bigl(\qff
\partial\dff
k_{\dff n}\dff(\dff \sigma\fff,\qff \tau\fff,\qff \tau'\trf)
\qff\bigr)
\off =\off\qff
f\trf(\dff \sigma\fff,\qff \tau\trf)
\qff -\qff 
f\trf(\dff \sigma\fff,\qff \tau'\trf)
\off -\off
\qff
\sum'\nolimits_{\qff i}\pff
f\qff
\bigl(\qff
k_{\dff n\dff -\dff 1}\dff 
(\qff \partial_{\dff i}\dff \sigma\fff,\qff 
\partial_{\dff i}\dff \tau\fff,\qff 
\partial_{\dff i}\dff \tau'\qff) 
\qff\bigr)
\pff,
\]

\vspace{-6pt}
or\halfff,\oss equivalently,\vspace{4.5pt}
\[
\quad
\partial^{\dff *}\nsp f\qff
\bigl(\qff
k_{\dff n}\dff(\dff \sigma\fff,\qff \tau\fff,\qff \tau'\trf)
\qff\bigr)
\off =\off\qff
f\trf(\dff \sigma\fff,\qff \tau\trf)
\qff -\qff 
f\trf(\dff \sigma\fff,\qff \tau'\trf)
\off -\off
\qff
\sum'\nolimits_{\qff i}\pff
f\qff
\bigl(\qff
k_{\dff n\dff -\dff 1}\dff 
(\qff \partial_{\dff i}\dff \sigma\fff,\qff 
\partial_{\dff i}\dff \tau\fff,\qff 
\partial_{\dff i}\dff \tau'\qff) 
\qff\bigr)
\pff.
\]

\vspace{-6pt}
Next,\oss let\dss us\sss apply\dss $m_{\dff n}\dff\sco{\tau'\fff}$\sss to\sss the\sss terms of\trs
this equality.\oss
By\dss the definition of\dss $h_{\dff n\dff +\dff 1}$\nsp,\vspace{4.5pt}
\[
\quad
m_{\dff n}\dff\sco{\tau'\fff}\qff\partial^{\dff *}\nsp f\qff
\bigl(\qff
k_{\dff n}\dff(\dff \sigma\fff,\qff \tau\fff,\qff \tau'\trf)
\qff\bigr)
\off =\off
h_{\dff n\dff +\dff 1}\dff\bigl(\trf \partial^{\dff *}\nsp f\qff\bigr)\trf
(\dff \sigma\fff,\qff \tau\trf)
\pff.
\]

\vspace{-7.5pt}
Since\sss $f\trf(\dff \sigma\fff,\qff \tau\trf)$\sss does not\sss depend
on\sss $\tau'$\nnsp,\vspace{3.5pt}
\[
\quad
m_{\dff n}\dff\sco{\tau'\fff}\trf f\trf(\dff \sigma\fff,\qff \tau\trf)
\off =\off
f\trf(\dff \sigma\fff,\qff \tau\trf)
\pff.
\]

\vspace{-8.5pt}
By\dss the definition of\dss $m^{\dff *}\dff(\trf f\trf)$\nsp,\vspace{4.5pt}
\[
\quad
m_{\dff n}\dff\sco{\tau'\fff}\trf f\trf(\dff \sigma\fff,\qff \tau'\trf)
\off =\off
m^{\dff *}\dff(\trf f\trf)\trf(\dff \sigma\trf)
\off =\off
p^{\dff *}\dff \circ\qff m^{\dff *}\dff(\trf f\trf)\trf(\dff \sigma\fff,\qff \tau\trf)
\pff.
\]

\vspace{-7.5pt}
Finally,\oss the coherence of\trs the family\sss $m_{\dff n}$\dss 
implies\sss that\vspace{4.5pt}
\[
\quad
m_{\dff n}\dff\sco{\tau'\fff}\trf 
\sum'\nolimits_{\qff i}\pff
f\qff
\bigl(\qff
k_{\dff n\dff -\dff 1}\dff(\qff \partial_{\dff i}\dff \sigma\fff,\qff 
\partial_{\dff i}\dff \tau\fff,\qff 
\partial_{\dff i}\dff \tau'\qff)
\qff\bigr)
\off
=\off 
\sum'\nolimits_{\qff i}\pff
m_{\dff n}\dff\sco{\tau'\fff}\trf
f\qff
\bigl(\qff
k_{\dff n\dff -\dff 1}\dff(\qff \partial_{\dff i}\dff \sigma\fff,\qff 
\partial_{\dff i}\dff \tau\fff,\qff 
\partial_{\dff i}\dff \tau'\qff)
\qff\bigr)
\pff
\]

\vspace{-29pt}
\[
\quad
\off
=\off 
\sum'\nolimits_{\qff i}\pff
m_{\dff n\dff -\dff 1}\dff\sco{\fff\rho}\trf
f\qff
\bigl(\qff
k_{\dff n\dff -\dff 1}\dff(\qff \partial_{\dff i}\dff \sigma\fff,\qff 
\partial_{\dff i}\dff \tau\fff,\pff 
\rho\qff)
\qff\bigr)
\pff
\]

\vspace{-29pt}
\[
\quad
\off
=\off 
\sum'\nolimits_{\qff i}\pff
h_{\dff n}\dff(\trf f\trf)\trf
(\qff \partial_{\dff i}\dff \sigma\fff,\qff 
\partial_{\dff i}\dff \tau\qff)
\off =\off\dff 
\partial^{\dff *}\dff 
\bigl(\qff
h_{\dff n}\dff(\trf f\trf)
\qff\bigr)\trf
(\qff \sigma\fff,\qff \tau\qff)
\pff.
\]

\vspace{-7.5pt}
By\sss collecting\sss all\dss these observations\sss together\halfff,\oss
we see\sss that\vspace{3pt}
\[
\quad
h_{\dff n\dff +\dff 1}\dff\bigl(\trf \partial^{\dff *}\nsp f\qff\bigr)
\off =\off
f
\qff -\qff
p^{\dff *}\dff \circ\qff m^{\dff *}\dff(\trf f\trf)
\qff -\qff
\partial^{\dff *}\dff 
\bigl(\qff
h_{\dff n}\dff(\trf f\trf)
\qff\bigr)
\]

\vspace{-9pt}
for every\sss
$f\qff \in\pff B^{\fff n}\dff(\trf K\dff \times\dff \Gamma\trf)$\sss
and\dss hence\vspace{3pt}
\[
\quad
h_{\dff n\dff +\dff 1}\dff\circ\qff \partial^{\dff *}
\off =\off\qff
\id
\pff -\pff
p^{\dff *}\dff \circ\qff m^{\dff *}
\pff -\pff
\partial^{\dff *}\dff \circ\qff h_{\dff n}
\pff,
\]

\vspace{-9pt}
or\halfff,\oss equivalently,\pss
$\id
\pff -\pff
p^{\dff *}\dff \circ\qff m^{\dff *}
\off =\off
h_{\dff n\dff +\dff 1}\dff\circ\qff \partial^{\dff *}
\pff +\pff
\partial^{\dff *}\dff \circ\qff h_{\dff n}$\nsp.\oss  \eproof\vspace{0.25pt}

\mypar{Theorem.}{blow-up-isomorphism}
\emph{The projection\dss
$p\dff \colon\dff
K\dff \times\dff \Gamma
\qff \ttoo\qff
K$\dss
induces\sss isometric\sss isomorphisms\sss in\dss the  
bounded cohomology\dss groups.\oss}

\proof
Lemma\qss \ref{main-homotopy}\qss together\dss with\dss
$m^{\dff *}\dff \circ\qff p^{\dff *}\off =\off \id$\dss
implies\sss that\dss the induced\sss homomorphisms are isomorphisms.\oss
Since\sss the norms of\dss $m^{\dff *}$\sss and\sss $p^{\dff *}$
are\dss $\leq\qff 1$\nnsp,\oss
these\sss induced\dss homomorphisms are isometries.\oss  \eproof

\newpage
\mysection{Isometric\qss isomorphisms\qss in\qss bounded\qss cohomology}{main-theorems}

\mypar{Theorem.}{bundle-isomorphism}
\emph{Let\dss $p\dff \colon\dff E\qff \ttoo\qff B$\dss
be a\sss locally\dss trivial\dss bundle with\dss 
the\dss fiber\dss $K\dff(\dff \pi\fff,\qff n\trf)$\nnsp.\oss
If\qss $n\qff >\qff 1$\nnsp,\oss
then\dss the map induced\dss by $p$\sss in\sss bounded cohomology\dss 
is\dss an\dss isometric\sss isomorphism.\oss}

\proof
Let\dss 
$\bm{\Gamma}\off =\off \bm{\Delta}\fff [\dff \infty\dff]$\nnsp.\oss
Let\dss us\sss consider\dss the diagram\oss\vspace{3pt}
\[
\quad
\begin{tikzcd}[column sep=booom, row sep=booom]\dis
E
\arrow[d, "\dis p\phantom{\dff \times\dff \id}"]
&
E\dff \times\dff \Gamma
\arrow[l]
\arrow[r]
\arrow[d, "\dis p\dff \times\dff \id"]
&
E\dff \times\dff \bm{\Gamma}\phantom{\dff,}
\arrow[d, "\dis p\dff \times\dff \id"]
\\
B 
&
B\dff \times\dff \Gamma
\arrow[l]
\arrow[r]
&
B\dff \times\dff \bm{\Gamma}\dff,
\end{tikzcd}
\]

\vspace{-9pt}
where\sss the\sss left\dss horizontal\sss arrows are projections,\oss
and\dss the right\dss horizontal\sss arrows are inclusions.\oss
The unravellings\sss $E\dff \times\dff \Gamma$\sss and\sss $B\dff \times\dff \Gamma$\sss
are only $\Delta$\dnsp-sets,\oss
and\sss the arrows of\trs this diagram are simplicial\dss maps of\dss $\Delta$\dnsp-sets,\oss
except\sss of\dss 
$p\dff \colon\dff 
E\qff \ttoo\qff B$\dss
and\dss
$p\dff \times\dff \id\dff \colon\dff 
E\dff \times\dff \bm{\Gamma}\qff \ttoo\qff B\dff \times\dff \bm{\Gamma}$\dnsp,\oss
which are maps of\dss simplicial\sss sets.\oss
Clearly,\oss this diagram\dss is\dss commutative.\oss
By\qss Theorem\qss \ref{blow-up-isomorphism}\qss 
the\sss left\dss horizontal\sss arrows induce\sss isometric\sss isomorphisms\sss in\sss bounded cohomology.\oss
Therefore,\oss it\dss is\dss
sufficient\dss to prove\sss that\dss the simplicial\dss map\dss
$q
\off =\off
p\dff \times\dff \id\dff \colon\dff 
E\dff \times\dff \Gamma\qff \ttoo\qff B\dff \times\dff \Gamma$\dss
induces\sss isometric\sss isomorphisms\sss in\sss bounded cohomology.\oss
Since\dss $\bm{\Gamma}$\sss is\dss contractible,\oss
the projections\dss
$E\dff \times\dff \bm{\Gamma}\qff \ttoo\qff E$\sss
and\sss
$B\dff \times\dff \bm{\Gamma}\qff \ttoo\qff B$\sss
are homotopy\sss equivalences.\oss
In\dss particular\halfff,\oss
these projections induce\sss isometric\sss isomorphisms\sss in\sss bounded cohomology.\oss
It\dss follows\sss that\dss the right\dss horizontal\sss arrows 
of\trs the diagram\sss induce isometric\sss isomorphisms\sss in\sss bounded cohomology.\oss

Clearly,\oss the simplicial\dss map\sss
$\bm{q}
\off =\off
p\dff \times\dff \id
\dff \colon\dff
E\dff \times\dff \bm{\Gamma}\qff \ttoo\qff B\dff \times\dff \bm{\Gamma}$\dss
is\dss a\sss locally\dss trivial\dss bundle.\oss
The main\dss part\sss of\trs the proof\trs is\dss an application of\trs
the\sss theory\sss developed\sss in\trs Section\qss \ref{bundles}\qss 
to $\bm{q}$\sss in\sss the role of\dss $p$\nnsp.\oss
Let\dss 
$G
\off =\off
\mathcal{C}^{\dff n\dff -\dff 1}\dff
(\qff B\dff \times\dff \bm{\Gamma}\fff,\qff \pi\trf(\dff \bm{q}\trf)\trf)$
be\sss the group of\dss normalized 
$(\dff n\qff -\qff 1\dff)$\dnsp-cochains of\trs 
$B\dff \times\dff \bm{\Gamma}$\sss
with\sss coefficients in\sss 
the\sss local\sss system\sss $\pi\trf(\dff \bm{q}\trf)$\nnsp.\oss
The group\sss $G$\sss 
acts on\sss $E\dff \times\dff \bm{\Gamma}$\sss
by\sss homotopic\sss to\sss the\sss identity\sss automorphisms 
over\sss 
$B\dff \times\dff \bm{\Gamma}$\dnsp.\oss 
The group $G$\sss is\dss abelian and\dss hence\dss is\dss amenable.\oss

Clearly,\oss a simplex of\dss $E\dff \times\dff \bm{\Gamma}$\sss
belongs\sss to\dss $E\dff \times\dff \Gamma$\sss
if\trs and\dss only\trs if\trs its image in\sss
$B\dff \times\dff \bm{\Gamma}$\sss
belongs\sss to\dss $B\dff \times\dff \Gamma$\dnsp.\oss
It\dss follows\dss that\dss the action of\dss $G$\sss leaves\sss the $\Delta$\dnsp-subset\sss
$E\dff \times\dff \Gamma$\sss invariant.\oss
Obviously,\oss every\sss simplex of\trs $\Gamma$\sss is\dss free in every\sss dimension.\oss
It\dss follows\dss that\sss every\sss simplex of\trs 
$B\dff \times\dff \Gamma$\sss is\dss free in every\sss dimension.\oss
Hence\trs Lemma\qss \ref{transitivity-bundle}\qss implies\sss that\dss
$E\dff \times\dff \Gamma/\fff G\off =\off B\dff \times\dff \Gamma$\sss 
and\vspace{3pt}
\[
\quad
q^{\dff *}\qff \colon\trf 
B^{\fff *}\dff(\trf B\dff \times\dff \Gamma\trf)
\pff \ttoo\off
B^{\fff *}\dff(\trf E\dff \times\dff \Gamma\trf)
\]

\vspace{-9pt}
induces isomorphism\dss from 
$B^{\fff *}\dff(\trf B\dff \times\dff \Gamma\trf)$
to\sss the space $G$\dnsp-invariant\sss cochains in\sss 
$B^{\fff *}\dff(\trf E\dff \times\dff \Gamma\trf)$\nnsp.\oss
Let\vspace{3pt}
\[
\quad
q^{\dff **}\qff \colon\trf 
\widehat{H}^{\fff *}\dff(\trf B\dff \times\dff \Gamma\trf)
\pff \ttoo\off
\widehat{H}^{\fff *}\dff(\trf E\dff \times\dff \Gamma\trf)
\]

\vspace{-9pt}
be\sss the map induced\sss by\sss $q^{\dff *}$\dnsp.\oss
Let\dss us prove\sss that\sss $q^{\dff **}$\sss
is\dss surjective.\oss
Let 
$\gamma\qff \in\pff B^{\fff m}\dff(\trf E\dff \times\dff \Gamma\trf)$
be a cocycle.\oss 
Since
$E\dff \times\dff \Gamma
\qff \ttoo\qff
E\dff \times\dff \bm{\Gamma}$\sss
induces\sss isomorphisms in\sss bounded cohomology,\oss
there exists
a cocycle
$c\qff \in\pff B^{\fff m}\dff(\trf E\dff \times\dff \bm{\Gamma}\trf)$
such\dss that\dss the restriction of\sss $c$\sss to
$E\dff \times\dff \Gamma$\sss
is\dss cohomologous\sss to $\gamma$\nnsp.\oss
Since $G$\sss is\dss amenable
and\sss acts on\sss $E$\sss by\sss automorphisms\sss homotopic\sss
to\sss the identity,\oss
the cocycle $c$\sss is\dss
cohomologous\sss to a $G$\dnsp-invariant\dss bounded cocycle $b$\nnsp.\oss
See\trs Lemma\qss \ref{invariant-cochains}.\oss
Let\sss $\beta$\sss be\sss the restriction of\dss $b$\sss to\sss
$E\dff \times\dff \bm{\Gamma}$\nnsp.\oss
Then\sss $\beta$\sss is\dss cohomologous\sss to $\gamma$
and\dss is\sss $G$\dnsp-invariant.\oss
Since\sss $\beta$\sss is $G$\dnsp-invariant,\pss
$\beta\off =\off q^{\dff *}\dff(\dff \alpha\trf)$\dss
for some\dss
$\alpha
\qff \in\pff 
B^{\fff m}\dff(\trf B\dff \times\dff \bm{\Gamma}\trf)$\nnsp.\oss
It\dss follows\sss that\sss $q^{\dff **}$\sss 
is\dss surjective.\oss

Let\dss us prove\sss now\sss that\sss $q^{\dff **}$\sss
is\dss injective.\oss
Since $G$\sss is\dss amenable,\oss
there exists a\sss $G$\dnsp-invariant\dss mean\sss
$\mu\dff \colon\dff
B\dff(\trf G\trf)\qff \ttoo\qff \rrr$\nnsp.\oss
For each $m\qff \in\pff \nnn$\sss let\dss us define a map\vspace{3pt}\vspace{-1.2pt}
\[
\quad
\mu_{\dff *}\dff \colon\dff
B^{\fff m}\dff(\trf E\dff \times\dff \Gamma\trf)
\pff \ttoo\pff
B^{\fff m}\dff(\trf B\dff \times\dff \Gamma\trf)
\]

\vspace{-9pt}\vspace{-1.2pt}
as follows.\oss
Let\dss $c\qff \in\qff B^{\fff m}\dff(\trf E\dff \times\dff \Gamma\trf)$
and $\sigma$\sss is\dss an $m$\dnsp-simplex of\dss 
$B\dff \times\dff \Gamma$\nnsp.\oss
Let\dss us\sss choose an $m$\dnsp-simplex\sss 
$\sigma\fff'$\sss of\dss $E\dff \times\dff \Gamma$
such\dss that\sss $p\trf(\dff \sigma\fff'\trf)\off =\off \sigma$\sss
and consider\dss the function\dss 
$g
\off \longmapsto\off
c\trf(\trf g\dff \cdot\dff \sigma\fff'\trf)$\sss
on\sss $G$\nnsp.\oss
Let\dss the value of\trs the cochain\sss 
$\mu_{\dff *}\dff(\dff c\trf)$\sss
on $\sigma$\sss be\dss equal\dss to\sss the value of\dss $\mu$ on\sss this function.\oss
Since $\mu$\sss is\sss $G$\dnsp-invariant,\oss this value\dss is\dss
independent\sss on\sss the choice of\dss $\sigma\fff'$\nnsp.\oss
Therefore\sss $\mu_{\dff *}$ is\dss well-defined.\oss
Clearly,\oss the composition\sss 
$\mu_{\dff *}\dff \circ\qff q^{\dff *}$\sss
is\dss equal\dss to\sss the identity.\oss
Since $G$ acts by\sss automorphisms,\pss
$\mu_{\dff *}$ commutes with\sss the duals\sss $\partial_{\dff i}^{\dff *}$\sss of\trs
the face operators\sss $\partial_{\dff i}$ and\dss hence with\sss the coboundary\sss
operator\sss $\partial^{\dff *}$\dnsp.\oss
Hence\sss $\mu_{\dff *}$\sss leads\sss to maps\vspace{3pt}\vspace{-1.2pt}
\[
\quad
\mu_{\dff **}\dff \colon\dff
\widehat{H}^{\fff m}\dff(\trf E\dff \times\dff \Gamma\trf)
\pff \ttoo\pff 
\widehat{H}^{\fff m}\dff(\trf B\dff \times\dff \Gamma\trf)
\pff.
\]

\vspace{-9pt}\vspace{-1.2pt}
Since\sss the composition\sss 
$\mu_{\dff *}\dff \circ\qff q^{\dff *}$\sss
is\dss equal\dss to\sss the identity,\oss
the composition\sss 
$\mu_{\dff **}\dff \circ\qff q^{\dff **}$\sss
is\dss also equal\dss to\sss the identity.\oss
It\dss follows\dss that\sss $q^{\dff **}$ is\dss injective.\oss 

We see\sss that\sss $q^{\dff **}$ is\dss an\sss isomorphism.\oss
It\dss remains\sss to prove\sss that\sss $q^{\dff **}$ is\dss an\sss isometric isomorphism.\oss
Since\sss $q^{\dff **}$\sss is\dss induced\dss by\sss
a simplicial\dss map,\oss 
the norm of\dss $q^{\dff **}$\sss is\qss $\leq\qff 1$\nnsp.\oss
On\sss the other\dss hand,\oss the norm of\dss $\mu$\sss is\qss $\leq\qff 1$\sss
and\dss hence\sss the norms of\dss $\mu_{\dff *}$\sss and\sss
$\mu_{\dff **}$\sss are also\qss $\leq\qff 1$\nnsp.\oss
Since\sss $q^{\dff **}$\sss is\dss an\sss isomorphism
and\sss $\mu_{\dff **}\dff \circ\dff q^{\dff **}$\sss
is\dss the identity,\pss
$\mu_{\dff **}$\sss is\dss the inverse of\dss $q^{\dff **}$\dnsp.\oss
So,\oss the norms of\dss $q^{\dff **}$\sss and of\dss its\sss inverse
are\qss $\leq\qff 1$\nnsp.\oss
It\dss follows\dss that\dss $q^{\dff **}$\sss 
is\dss an\sss isometry.\oss  \eproof

\myuppar{Remark.}
One can\sss prove\sss that\sss $q^{\dff **}$ is\dss
an\sss isometry\sss in a different\dss way.\oss
Let\sss us\sss return\sss to\sss the proof\dss of\trs the surjectivity\dss
of\dss $q^{\dff **}$\snsp.\dff\oss
Since
$E\dff \times\dff \Gamma
\qff \ttoo\qff
E\dff \times\dff \bm{\Gamma}$\sss
induces\sss isometric\sss isomorphisms in\sss bounded cohomology,\oss
one can choose\sss $c$\sss in such a way\dss that\dss
$\norm{\gamma}\qff \geq\qff
\norm{c}$\nnsp.\oss
By\qss Lemma\qss \ref{invariant-cochains}\qss
one can choose\sss $b$\sss in such a way\dss that\dss
$\norm{c}\qff \geq\qff
\norm{b}$\nnsp.\oss
Clearly,\pss
$\norm{b}\qff \geq\qff
\norm{\beta}$\dss
and\dss
$\norm{\beta}
\off =\off
\norm{\alpha}$\nnsp.\oss
It\dss follows\dss that\dss
$\norm{\gamma}\qff \geq\qff \norm{\alpha}$\nnsp.\oss
Therefore
for every\sss cohomology\sss class\dss
$\bm{\gamma}
\qff \in\pff 
\widehat{H}^{\fff m}\dff(\trf E\dff \times\dff \Gamma\trf)$\sss
there exists a cohomology\sss class\dss
$\bm{\alpha}
\qff \in\pff 
\widehat{H}^{\fff m}\dff(\trf B\dff \times\dff \Gamma\trf)$\sss
such\dss that\dss
$\bm{\gamma}
\off =\off
q^{\dff **}\dff(\dff \bm{\alpha}\trf)$\sss
and\dss
$\norm{\bm{\alpha}}\off \leq\off \norm{\bm{\gamma}}$\nnsp.\oss
Since,\oss at\dss the same\sss time,\oss
$\norm{q^{\dff **}\dff(\dff \bm{\alpha}\trf)}
\off \leq\off
\norm{\bm{\alpha}}$\nnsp,\oss
the isomorphism\sss $q^{\dff **}$\sss is\dss 
an\sss isometric\sss isomorphism.

\mypar{Theorem.}{quotient-by-amenable}
\emph{Let $\pi$ be a discrete group
and\dss $\kappa\qff \subset\qff \pi$ be a normal\sss
amenable subgroup of\dss $\pi$\nnsp.\oss
Let\dss 
$p\dff \colon\dff
\pi\qff \ttoo\qff \pi/\kappa$\dss
be\sss the quotient\dss homomorphism.\oss
Then\dss
$\mathit{B}\fff p\dff \colon\dff
\mathit{B}\trf \pi
\qff \ttoo\qff
\mathit{B}\trf (\dff \pi/\kappa\dff)$\dss
induces\sss isometric\sss isomorphism\sss in\sss bounded cohomology.\oss}

\proof
The proof\dss is\dss completely\sss similar\dss to\sss the proof\dss of\qss
Theorem\qss \ref{bundle-isomorphism},\oss
with\dss the\sss map\dss\vspace{1.5pt}\vspace{-0.375pt}
\[
\quad
\mathit{B}\fff p\dff \colon\dff
\mathit{B}\trf \pi
\qff \ttoo\qff
\mathit{B}\trf (\dff \pi/\kappa\dff)
\]

\vspace{-10.5pt}\vspace{-0.375pt}
playing\dss the role of\dss
$p\dff \colon\dff E\qff \ttoo\qff B$\nnsp.\oss
As we saw\sss in\trs Section\qss \ref{classifying-spaces},\oss
the group\dss
$G
\off =\off
C_{\trf 0}\dff(\trf \nnn\fff,\pff \kappa\trf)$\sss acts on\sss
$\mathit{B}\trf \pi\dff \times\dff \Gamma$\dss
and\dss
$\mathit{B}\trf \pi\dff \times\dff \Gamma\halfff/\halfff G
\off =\off
\mathit{B}\trf (\dff \pi/\kappa\dff)\dff \times\dff \Gamma$\dnsp.\oss
See\qss (\ref{nice-quotient-b}).\oss
While,\oss to\sss the best\sss of\dss author's\sss
knowledge,\oss the group 
$C^{\trf 0}\dff(\trf V\fff,\pff \kappa\trf)$
is\dss not\dss known\dss to be amenable,\oss
the group
$C_{\dff 0}\dff(\trf V\fff,\qff \kappa\trf)$
is\dss a direct\sss sum 
of\dss copies of\sss $\kappa$
and\dss hence\dss is\dss amenable.\oss
Hence one can argue as\sss in\dss the proof\dss of\qss
Theorem\qss \ref{bundle-isomorphism}\qss and conclude\sss
that\sss $\mathit{B}\fff p$\sss induces\sss isometric\sss isomorphism 
in\dss bounded cohomology.\oss  \eproof

\myuppar{The fundamental\dss group.}
Let\dss $K$\sss be a connected\trs Kan\dss simplicial\sss set.\oss
Suppose\sss that\sss $K$\sss has only\sss one vertex,\oss
which\sss we will\sss denote by $v$\nnsp.\oss
Let\dss us\dss interpret\sss a $1$\dnsp-simplex $\sigma\qff \in\pff K_{\dff 1}$
as a\sss loop based at\sss $v$\nnsp.\oss 
The\dss Kan\dss extension\sss property\dss implies\sss that\dss
for every\dss two $1$\dnsp-simplices\sss $\rho\fff,\pff \sigma$\sss
there exists a $2$\dnsp-simplex $\omega$ such\dss that\dss
$\rho\off =\off \partial_{\trf 2}\dff \omega$\sss
and\sss
$\sigma\off =\off \partial_{\trf 0}\dff \omega$\nnsp.\oss
One can easily\sss check\dss that\sss up\sss to homotopy\sss
$\tau\off =\off \partial_{\dff 1}\dff \omega$\sss
does not\sss depends on\sss the choice of\dss $\omega$\nnsp,\oss
and,\oss moreover\halfff,\oss up\sss to homotopy\sss $\tau$\sss
depends only\sss on\sss the homotopy\sss classes of\dss
$\sigma\fff,\pff \tau$\nnsp.\oss
One can\sss take\sss the homotopy\sss class of\sss $\tau$\sss
as\sss the\qss \emph{product}\pss $r\dff \cdot\dff s$\sss of\trs the homotopy\sss classes\sss
$r\fff,\pff s$\sss 
of\trs $\rho\fff,\pff\sigma$\sss respectively.\oss
The set\sss of\trs homotopy\sss classes of\dss $1$\dnsp-simplicies\sss
together\sss with\dss this product\dss is\dss the\qss
\emph{fundamental\dss group}\qss
$\pi_{\fff 1}\dff(\trf K\fff,\qff v\trf)$ of\trs $K$\nnsp.\oss
If\halfff,\oss in\sss addition,\pss $K$\sss is\dss minimal,\oss
then every\dss two homotopic $1$\dnsp-simplices are equal.\oss
In\dss this case\sss $\pi_{\fff 1}\dff(\trf K\fff,\qff v\trf)$\sss
can\sss be identified\sss with\sss $K_{\dff 1}$\sss as a set.\oss\vspace{-0.125pt}

\mypar{Lemma.}{f-group}
\emph{Suppose\sss that\qss $K$\dss is\dss a connected\dss minimal\trs Kan\dss
simplicial\sss set.\oss
Then\trs $K\dff(\dff 1\dff)$
is\dss canonically\dss isomorphic\sss to\sss
$\mathit{B}\qff \pi_{\fff 1}\dff(\trf K\fff,\qff v\trf)$\nnsp,\oss
where $v$ is\dss the unique vertex of\qss $K$\nnsp.\oss}

\proof
For\dss $i\fff,\pff j\fff,\pff n\qff \in\pff \nnn$\dss 
such\dss that\dss 
$0\qff \leq\qff i\qff <\qff j\qff \leq\qff n$\dss let\dss
$\theta_{\dff i\fff,\dff j}\qff \colon\dff [\dff 1\dff]\qff \ttoo\qff [\halfff n\dff]$\dss
be\sss the map\vspace{1.5pt}\vspace{-0.375pt}
\[
\quad
\theta_{\dff i\fff,\dff j}\dff\colon\dff
0\off \longmapsto\off i\qff,\quad\
1\off \longmapsto\off j
\pff.
\]

\vspace{-10.5pt}\vspace{-0.375pt}
Suppose\sss that\dss 
$\rho_{\dff 1}\dff,\off \rho_{\trf 2}\dff,\off \ldots\dff,\off \rho_{\dff n}$\dss
are $1$\dnsp-simplices of\trs $K$\nnsp.\oss
If\vspace{1.5pt}\vspace{-0.375pt}
\[
\quad
\theta^{\dff *}_{\dff i\dff -\dff 1\fff,\dff i}\qff(\trf \sigma\trf)
\off =\off
\rho_{\dff i}
\]

\vspace{-10.5pt}\vspace{-0.375pt}
for some $n$\dnsp-simplex $\sigma$ of\dss $K$
and every $i$ between $1$ and $n$\nnsp,\oss then\vspace{1.5pt}\vspace{-0.375pt}
\[
\quad
\theta^{\dff *}_{\dff i\fff,\dff j}\qff(\trf \sigma\trf)
\off =\off
\rho_{\dff i}\dff \cdot\off \ldots\off \cdot\dff \rho_{\dff j}
\]

\vspace{-10.5pt}\vspace{-0.375pt}
for every\dss $i\qff <\qff j$\nnsp.\oss
This\dss follows\dss from\dss the definition of\trs the product\dss together\dss
with an\dss induction\dss by\dss $j\qff -\qff i$\nnsp.\oss
In\dss turn,\oss this implies\sss that\dss the restriction of\dss $i_{\dff \sigma}$\sss
to\sss $\sk_{\trf 1}\fff \bm{\Delta}\fff [\halfff n\dff]$\dss
is\dss uniquely\sss determined\dss by\dss
$\rho_{\dff 1}\dff,\off \rho_{\trf 2}\dff,\off \ldots\dff,\off \rho_{\dff n}$\nsp.\oss
On\dss the other\dss hand,\oss Kan\dss extension\sss property\dss
implies\sss that\sss such a simplex $\sigma$ exists for every
$n$\dnsp-tuple\dss
$\rho_{\dff 1}\dff,\off \rho_{\trf 2}\dff,\off \ldots\dff,\off \rho_{\dff n}$\nsp.\oss
It\dss follows\dss that\sss one can\dss identify\sss 
$n$\dnsp-simplices of\trs $K\dff(\dff 1\dff)$\sss
with sequences\dss
$(\trf \rho_{\dff 1}\dff,\off \rho_{\trf 2}\dff,\off \ldots\dff,\off \rho_{\dff n}\trf)$\dss
of\dss elements of\dss
$\pi_{\fff 1}\dff(\trf K\fff,\qff v\trf)$\nnsp,\oss
i.e.\qss with $n$\dnsp-simplices of\trs
$\mathit{B}\qff \pi_{\fff 1}\dff(\trf K\fff,\qff v\trf)$\nnsp.\oss
A direct\sss check\sss shows\sss that\dss this identification\sss
respects\sss the boundary\sss and\sss degeneracy\sss operators.\oss 
The\sss lemma\sss follows.\oss \eproof

\mypar{Theorem.}{bounded-and-pione}
\emph{Let\qss $K$\sss be\dss a connected\trs Kan\dss
simplicial\sss set\sss
and\dss
$f\dff \colon\dff
K\qff \ttoo\qff
\mathit{B}\qff \pi_{\fff 1}\dff(\trf K\fff,\qff v\trf)$\nnsp,\oss
where $v$\sss is\dss a vertex of\pss $K$\nnsp,\oss
be a simplicial\dss map\sss inducing isomorphism of\trs fundamental\dss groups.\oss
Then\sss $f$\sss induces an\sss isomorphism\dss in\dss bounded cohomology.\oss}\vspace{-0.125pt}

\proof
The proof\dss is\dss based on\dss the\sss theory of\qss Postnikov\dss systems.\oss
See\trs Section\qss \ref{postnikov}\qss for a review of\trs the definitions and\dss
the\sss theorems used\sss in\dss this proof\halfff.\oss

Let\dss $\pi_{\dff 1}$\sss be\sss the fundamental\dss group of\dss $K$\nnsp.\oss
Since\sss $\mathit{B}\qff \pi_{\dff 1}$\sss is\dss a\dss Kan\sss simplicial\sss set,\oss
every\dss two map\dss
$K\qff \ttoo\qff \mathit{B}\qff \pi_{\dff 1}$\dss
inducing\sss isomorphism of\trs the fundamental\dss groups are homotopic.\oss
Hence\sss it\dss is\dss sufficient\dss to prove\sss to prove\sss the\sss theorem\sss
for one such\sss map.\oss
Let\sss $M$\sss be a\sss minimal\dss Kan\dss simplicial\sss subset\sss of\dss $K$\sss
which\dss is\dss a strong deformation\sss retract\sss of\dss $K$\nnsp.\oss
Since\sss $M$\sss is\dss minimal\sss and connected,\pss
$M$\sss has only one vertex,\oss
which we denote by\sss $v$\nnsp.\oss
Let\dss
$M\dff(\dff 0\dff)\fff,\pff
M\dff(\dff 1\dff)\fff,\pff
\ldots\fff,\pff
M\dff(\dff n\trf)\fff,\pff
\ldots$\sss
and\dss the maps\dss
$p_{\fff n}$\sss and\sss $p_{\fff m,\dff n}$\sss
be\sss the\dss Postnikov\dss system of\sss $M$\nnsp.\oss 
Then every\sss map\vspace{1.875pt}
\[
\quad
p_{\fff n,\dff n\dff -\dff 1}\dff \colon\dff
M\dff(\dff n\trf)
\qff \ttoo\qff
M\dff(\dff n\qff -\qff 1\dff)
\]

\vspace{-10.125pt}
is\dss a\sss locally\trs trivial\dss bundle\sss
with\sss the fiber\sss 
$K\dff(\dff \pi_{\dff n}\dff,\qff n\trf)$\nnsp,\oss 
where\sss
$\pi_{\dff n}
\off =\off
\pi_{\fff n}\dff(\trf M\fff,\qff v\trf)$\sss
is\dss the $n${\dnsp}th\dss homotopy\dss group of\dss $M$\nnsp.\oss
If\dss $n\qff >\qff 1$\nnsp,\oss 
then\sss $p_{\fff n,\dff n\dff -\dff 1}$\sss induces\sss isometric\sss isomorphism\sss
in\sss bounded cohomology.\oss
It\dss follows\sss that\sss for every\sss $n\qff >\qff 1$\sss 
the map\vspace{1.875pt}
\[
\quad
p_{\fff n,\dff 1}\dff \colon\dff
M\dff(\dff n\trf)
\qff \ttoo\qff
M\dff(\dff 1\dff)
\]

\vspace{-10.125pt}
induces\sss isometric\sss isomorphism\sss in\sss bounded cohomology.\oss
On\sss the other\sss hand,\sss for\sss $n\qff \geq\qff m$\sss
the $m${\nnsp}th\dss skeletons of\dss $M$ and\sss $M\dff(\dff n\trf)$
are\sss the same by\dss the very\sss definition of\dss $M\dff(\dff n\trf)$\nnsp.\oss
By\dss the definition,\oss the bounded cohomology\sss group\sss 
$\widehat{H}^{\fff m}\dff(\trf M\trf)$\sss depends only on\sss the
$(\dff m\qff +\qff 1\dff)${\nnsp}th\dss skeleton\sss
$\sk_{\dff m\dff +\dff 1}\fff M$\sss of\dss $M$\nnsp.\oss
It\dss follows\sss that\sss the map\vspace{1.875pt}
\[
\quad
p_{\fff 1}\dff \colon\dff
M\qff \ttoo\qff M\dff(\dff 1\trf)
\pff
\]

\vspace{-10.125pt}
induces\sss isometric\sss isomorphism\sss in\sss bounded cohomology.\oss 
By\trs Lemma\qss \ref{f-group}\qss the simplicial\sss set\sss $M\dff(\dff 1\trf)$
is\dss canonically\dss isomorphic\sss to\sss
$\mathit{B}\qff \pi_{\fff 1}$\nsp.\oss
Moreover\halfff,\oss the description of\dss fundamental\dss groups preceding\qss
Lemma\qss \ref{f-group}\qss shows\sss that\sss
$p_{\fff 1}\dff \colon\dff
M\qff \ttoo\qff M\dff(\dff 1\trf)$\sss
induces isomorphism of\dss fundamental\sss groups.\oss
If\dss
$r\dff \colon\dff K\qff \ttoo\qff M$\sss
is\dss a strong deformation\sss retraction,\oss 
then $r$ induces isomorphism of\dss fundamental\dss groups
and\sss isometric\sss isomorphism\sss in\sss bounded cohomology.\oss
It\dss follows\sss that\sss\vspace{1.875pt}
\[
\quad
p_{\fff 1}\dff \circ\trf r\trf \colon\trf
K\qff \ttoo\qff
M\dff(\dff 1\trf)
\off =\off
\mathit{B}\qff \pi_{\fff 1}
\]

\vspace{-10.125pt}
also has\sss this property.\oss
This proves\sss the\sss theorem\sss for\sss
$f\off =\off p_{\fff 1}\dff \circ\trf r$\nnsp.\oss
As was pointed out\sss above,\oss
any\sss special\sss case of\trs the\sss theorem\dss
implies\sss the general\sss one.\oss
This completes\sss the proof\halfff.\oss  \eproof

\mypar{Corollary.}{isomorphism-in-pione}
\emph{Let\qss $K\fff,\pff L$\sss be\dss connected\trs Kan\dss
simplicial\sss sets.\oss
If\qss
$f\dff \colon\dff
K\qff \ttoo\qff L$\dss
is\dss a simplicial\sss map
inducing isomorphism of\trs fundamental\dss groups,\oss
then\sss $f$\sss induces isomorphism\dss in\dss bounded cohomology.\oss}  \eproof

\mypar{Theorem.}{iso-amenable-kernel}
\emph{Let\qss $K\fff,\pff L$\dss be connected\dss Kan\qss simplicial\sss sets
and\dss let\sss $v$ be a vertex of\qss $K$\nnsp.\oss
Let\qss $f\dff \colon\dff K\qff \ttoo\qff L$\qss be a simplicial\dss map.\oss
If\qss 
$f_{\dff *}\dff \colon\dff
\pi_{\fff 1}\dff(\trf K\fff,\qff v\trf)
\qff \ttoo\qff
\pi_{\fff 1}\dff(\trf L\fff,\qff f\dff(\dff v\trf)\trf)$\dss
is\dss surjective and\dss has amenable kernel,\oss
then\dss $f$\dss induces an\sss isometric\sss isomorphism\sss in\dss bounded cohomology.\oss}

\proof
In order\sss not\dss to clutter\sss the notations,\oss
we will\dss not\dss mention\dss the base points\sss
$v\fff,\pff f\dff(\dff v\trf)$\sss anymore.\oss
Let\dss us\dss consider\dss the diagram\oss\vspace{-1.5pt}
\[
\quad
\begin{tikzcd}[column sep=booom, row sep=booom]\dis
K
\arrow[r, "\dis f\dff"]
\arrow[d, "\dis p_{\trf K}\dff"']
&
L
\arrow[d, "\dis \dff p_{\trf L}"]
\\
\mathit{B}\qff \pi_{\fff 1}\dff(\trf K\trf)
\arrow[r, "\dis \mathit{B}\fff f_{\dff *}"]
&
\mathit{B}\qff \pi_{\fff 1}\dff(\trf L\trf)\dff,
\end{tikzcd}
\]

\vspace{-9pt}
where\sss
$p_{\trf K}\dff,\pff p_{\trf L}$\sss
are some maps inducing\sss isomorphisms of\dss fundamental\sss groups.\oss
By\qss Theorem\qss \ref{quotient-by-amenable}\qss the map\sss
$\mathit{B}\fff f_{\dff *}$\sss 
induces\sss isometric\sss isomorphism\sss in\dss bounded cohomology.\oss
By\qss Theorem\qss \ref{bounded-and-pione}\qss the maps\sss
$p_{\trf K}\dff,\pff p_{\trf L}$\sss
induce isometric\sss isomorphisms\sss in\dss bounded cohomology.\oss
Since\sss the above diagram\dss is\dss commutative up\sss to homotopy,\oss
it\dss follows\sss that\sss
$f$\sss induces\sss isometric\sss isomorphism\sss 
in\dss bounded cohomology.\oss  \eproof

\newpage
\myappend{The\qss constructions\qss of\pss Milnor\qss and\qss Segal}{milnor-segal}

\myapar{Lemma.}{milnor-space}
\emph{The $\Delta$\dnsp-set\dss $\mathcal{B}\dff \pi$\sss
is\dss canonically\dss isomorphic\sss to\sss the product\dss
$\mathit{B}\qff \pi\dff \times\dff \Delta\dff[\dff \infty\dff]$\nnsp.\oss}

\proof
To begin\sss with,\oss we will\sss give an explicit\sss description of\dss
simplices of\dss $\mathcal{B}\dff \pi$\nnsp.\oss
The $n$\dnsp-simplices of\dss $\mathcal{E}\dff \pi$\sss
can\sss be identified\sss with\dss pairs of\dss sequences\dss\vspace{3pt}
\[
\quad
(\qff g_{\trf 0}\fff,\pff g_{\dff 1}\fff,\pff \ldots\fff,\pff g_{\dff n}\trf)
\qff \in\qff
\pi^{\dff n\dff +\dff 1}\dff,\quad
(\qff k_{\trf 0}\fff,\pff k_{\dff 1}\fff,\pff \ldots\fff,\pff k_{\dff n}\trf)
\qff \in\qff
\nnn^{\dff n\dff +\dff 1}
\dff,
\]

\vspace{-9pt}
such\dss that\dss
$k_{\trf 0}\qff <\qff k_{\dff 1}\qff <\qff \ldots\qff <\qff k_{\dff n}$\nsp,\oss
and\dss $g\qff \in\qff \pi$\dss acts by\dss the rules\vspace{3pt}
\[
\quad
g\dff \cdot\trf
(\qff g_{\trf 0}\fff,\pff g_{\dff 1}\fff,\pff \ldots\fff,\pff g_{\dff n}\trf)
\off =\off
(\qff g\dff g_{\trf 0}\fff,\pff g\dff g_{\dff 1}\fff,\pff \ldots\fff,\pff g\dff g_{\dff n}\trf)
\quad\
\mbox{and}\quad\
\]

\vspace{-36pt}
\[
\quad
g\dff \cdot\dff
(\qff k_{\trf 0}\fff,\pff k_{\dff 1}\fff,\pff \ldots\fff,\pff k_{\dff n}\trf)
\off =\off
(\qff k_{\trf 0}\fff,\pff k_{\dff 1}\fff,\pff \ldots\fff,\pff k_{\dff n}\trf)
\pff.
\]

\vspace{-9pt}
In order\dss to give a direct\sss description of\dss simplices 
of\dss $\mathcal{B}\dff \pi$\nnsp,\oss
let\dss us\sss set
use\sss the\qss \emph{bar\dss notations}\vspace{3pt}
\[
\quad
g_{\trf 0}\qff [\trf g_{\dff 1}\fff\mid\dff g_{\dff 2}\fff\mid\dff \ldots\fff\mid\dff g_{\dff n}\trf]
\off =\off
(\qff g_{\trf 0}\fff,\pff g_{\trf 0}\dff g_{\dff 1}\fff,\pff \ldots\fff,\pff 
g_{\trf 0}\dff g_{\dff 1}\dff \ldots\dff g_{\dff n}\trf)
\pff.
\]

\vspace{-9pt}
In\dss these notations\sss the action of\dss 
$\pi$\sss takes\sss the form\vspace{3pt}
\[
\quad
g\dff \cdot\trf
\left(\qff
g_{\trf 0}\qff [\trf g_{\dff 1}\fff\mid\dff g_{\dff 2}\fff\mid\dff \ldots\fff\mid\dff g_{\dff n}\trf]
\qff\right)
\off =\off
g\dff g_{\trf 0}\qff [\trf g_{\dff 1}\fff\mid\dff g_{\dff 2}\fff\mid\dff \ldots\fff\mid\dff g_{\dff n}\trf]
\pff.
\]

\vspace{-9pt}
Therefore $n$\dnsp-simplices of\dss $\mathcal{B}\dff \pi$\sss
can\sss be identified\sss with\dss pairs of\dss sequences\dss\vspace{3pt}
\[
\quad
[\trf g_{\dff 1}\fff\mid\dff g_{\dff 2}\fff\mid\dff \ldots\fff\mid\dff g_{\dff n}\trf]
\qff \in\qff
\pi^{\dff n}\dff,\quad
(\qff k_{\trf 0}\fff,\pff k_{\dff 1}\fff,\pff \ldots\fff,\pff k_{\dff n}\trf)
\qff \in\qff
\nnn^{\dff n\dff +\dff 1}
\dff,
\]

\vspace{-9pt}
such\dss that\dss
$k_{\trf 0}\qff <\qff k_{\dff 1}\qff <\qff \ldots\qff <\qff k_{\dff n}$\nsp.\oss\vspace{-0.125pt}

The boundary\sss operators  act\sss independently\sss
on\sss these sequences.\oss
Namely,\oss the action of\trs the boundary\sss operator\dss
$\partial_{\dff i}$\sss on\dss the sequences\sss
$(\qff k_{\trf 0}\fff,\pff k_{\dff 1}\fff,\pff \ldots\fff,\pff k_{\dff n}\trf)$\sss
is\dss the same as\sss in $\Delta\dff[\dff \infty\dff]$\nnsp,\oss
and\dss the action on\dss the sequences\dss
$[\trf g_{\dff 1}\fff\mid\dff g_{\dff 2}\fff\mid\dff \ldots\fff\mid\dff g_{\dff n}\trf]$\sss
is\dss given\dss
by\dss the rules\vspace{4pt}
\[
\quad
\partial_{\trf 0}\qff 
[\trf g_{\dff 1}\fff\mid\dff g_{\dff 2}\fff\mid\dff \ldots\fff\mid\dff g_{\dff n}\trf]
\off =\off
[\trf g_{\dff 2}\fff\mid\dff g_{\dff 3}\fff\mid\dff \ldots\fff\mid\dff g_{\dff n}\trf]
\qff,
\]

\vspace{-35.5pt}
\[
\quad
\partial_{\dff n}\qff 
[\trf g_{\dff 1}\fff\mid\dff g_{\dff 2}\fff\mid\dff \ldots\fff\mid\dff g_{\dff n}\trf]
\off =\off
[\trf g_{\dff 1}\fff\mid\dff g_{\dff 2}\fff\mid\dff \ldots\fff\mid\dff g_{\dff n\dff -\dff 1}\trf]
\qff,\quad
\mbox{and}
\pff
\]

\vspace{-35.5pt}
\[
\quad
\partial_{\dff i}\qff 
[\trf g_{\dff 1}\fff\mid\dff g_{\dff 2}\fff\mid\dff \ldots\fff\mid\dff g_{\dff n}\trf]
\off =\off
[\trf g_{\dff 1}\fff\mid\dff 
\ldots\fff\mid 
g_{\dff i}\dff g_{\dff i\dff +\dff 1}\fff\mid\dff 
\ldots\fff\mid\dff 
g_{\dff n}\trf]
\quad
\mbox{for}\quad
0\qff <\qff i\qff <\qff n
\pff.
\]

\vspace{-8pt}
This differs\sss from\sss the definition of\trs $\partial_{\dff i}$ for
$\mathit{B}\trf \pi$ only\sss in\sss notations\qss
({\fff}the product $g_{\dff i}\dff g_{\dff i\dff +\dff 1}$
is\dss interpreted as\sss the composition $g_{\dff i\dff +\dff 1}\dff \circ\dff g_{\dff i}$\nsp).\oss
It\dss follows\sss that\sss 
$\mathcal{B}\dff \pi
\off =\off
\mathit{B}\trf \pi\dff \times\dff \Delta\dff[\dff \infty\dff]$\nnsp.\oss  \eproof

\myuppar{Unravelings\sss of\trs classifying spaces of\dss categories.}
For a category\sss $\mathcal{C}$\sss 
let\sss 
$\mathcal{C}_{\dff \bm{n}}$\sss 
be\sss the subcategory\sss of\dss
$\mathcal{C}\dff \times\dff \bm{n}$\sss
obtained\dss by\sss deleting all\dss morphisms of\trs the form\sss 
$(\dff c\fff,\qff n\trf)
\qff \ttoo\qff
(\dff c'\fff,\qff n\trf)$\sss
where $c\fff,\pff c'$ are objects of\dss $\mathcal{C}$ and
$n\qff \in\pff \nnn$\nnsp,\oss
except\dss identity\dss morphisms.\oss
This construction\dss is\dss due\sss to\trs
Segal\qss \cite{s},\oss
who called\sss $\mathcal{C}_{\dff \bm{n}}$\sss
the\qss \emph{unraveling}\pss of\sss $\mathcal{C}$\sss
over\sss the ordered\sss set\dss $\nnn$\sss
and\dss pointed\sss out\dss that\sss
for a group $\pi$\sss
the geometric realizations of\dss $\mathcal{B}\dff \pi$\sss
and\sss $\mathit{B}\trf \pi_{\dff \bm{n}}$\sss
are homeomorphic.\oss

This result\sss can\sss be interpreted\sss in\dss terms
of\dss simplicial\sss sets and extended\sss to arbitrary\sss categories.\oss
Namely,\oss Lemma\qss \ref{milnor-space}\qss suggests\sss
that\dss the $\Delta$\dnsp-set\sss 
$\mathcal{B}\fff \mathcal{C}
\off =\off
\mathit{B}\trf \mathcal{C}\dff \times\dff 
\Delta\dff[\dff \infty\dff]$\sss
is\dss an analogue of\dss $\mathcal{B}\dff \pi$\nnsp.\oss
In\sss contrast\dss with\dss the case groups,\oss
in\dss general\sss $\mathcal{B}\fff \mathcal{C}$\sss is\dss not\sss arising\dss
from a simplicial\sss complex.\oss
It\dss turns out\dss that\dss the simplicial\sss set\sss
$\bm{\Delta}\dff \mathcal{B}\fff \mathcal{C}$\dss
is\dss isomorphic\sss to\dss
$\mathit{B}\trf \mathcal{C}_{\dff \bm{n}}$\nsp.\oss
Before proving\dss this,\oss it\dss is\dss convenient\dss
to introduce\sss the notion of\trs the\qss \emph{core}\pss
of\dss a simplicial\sss set.\oss

\myuppar{The core of\dss a simplicial\sss set.}
Following\dss 
Rourke\dss and\dss Sanderson\qss \cite{rs},\oss
let\dss us\sss define\sss the\qss \emph{core}\pss of\dss a simplicial\sss set\sss $K$\sss 
as\sss the $\Delta$\dnsp-subset $\core(\trf K\trf)$ of\dss $K$\sss
consisting\sss of\trs simplices of\trs the form
$\theta^{\fff *}\dff(\dff \sigma\trf)$
with\dss non-degenerate $\sigma$ and strictly\dss increasing $\theta$\nnsp.\oss
The simplicial\sss set\sss $K$\sss is\dss said\dss to have\qss \emph{non-degenerate core}\pss 
if\trs non-degenerate simplices of\sss $K$
form a $\Delta$\dnsp-subset\sss of\sss $K$\nnsp.\oss
Clearly,\oss this $\Delta$\dnsp-subset\dss is\dss 
equal\dss to $\core(\trf K\trf)$\nnsp.\oss
There\dss is\dss a canonical\sss simplicial\dss map
$\Theta\dff \colon\dff
\bm{\Delta} \core(\trf K\trf)
\qff \ttoo\qff
K$\dss
defined\dss by\dss
$\Theta\trf(\dff \sigma\fff,\qff \rho\trf)
\off =\off
\rho^{\dff *}\dff(\dff \sigma\trf)$\nnsp.\oss

\myapar{Lemma.}{ndc}
\emph{If\pss $K$\sss has\dss non-degenerate core,\oss
then\sss $\Theta$\sss is\dss an\dss isomorphism.\oss}\vspace{-0.125pt}

\proof
Every\sss simplex $\sigma$ of\dss a simplicial\sss set\sss admits a unique presentation\sss
$\sigma\off =\off \theta^{\dff *}\dff(\dff \tau\trf)$\sss
with\dss non-degenerate $\tau$ and surjective $\theta$\nnsp.\oss
See\trs Lemma\qss \ref{eilenberg-zilber}.\oss
This\sss implies\sss that\sss $\Theta$\sss is\dss surjective.\oss
If\trs $K$\sss has non-degenerate core and\sss
$\Theta\trf(\dff \sigma_{\dff 1}\fff,\qff \rho_{\dff 1}\trf)
\off =\off
\Theta\trf(\dff \sigma_{\dff 2}\fff,\qff \rho_{\dff 2}\trf)$\nnsp,\oss
then\dss
$\rho_{\dff 1}^{\dff *}\dff(\dff \sigma_{\dff 1}\trf)
\off =\off
\rho_{\dff 2}^{\dff *}\dff(\dff \sigma_{\dff 2}\trf)$\dss
and\dss
$\sigma_{\dff 1}\fff,\pff \sigma_{\dff 2}$\dss
are non-degenerate.\oss
Therefore\sss the uniqueness part\sss of\qss Lemma\qss \ref{eilenberg-zilber}\qss
implies\sss that\dss
$(\dff \sigma_{\dff 1}\fff,\qff \rho_{\dff 1}\trf)
\off =\off
(\dff \sigma_{\dff 2}\fff,\qff \rho_{\dff 2}\trf)$\nnsp.\oss  \eproof

\myapar{Theorem.}{segal-space}
\emph{The simplicial\sss set\trs 
$\bm{\Delta}\dff \mathcal{B}\fff \mathcal{C}$\dss
is\dss isomorphic\sss to\dss
$\mathit{B}\trf \mathcal{C}_{\dff \bm{n}}$\nsp.\oss}

\proof
By\dss restricting\dss the projection\dss
$\mathcal{C}\dff \times\dff \bm{n}
\qff \ttoo\qff
\bm{n}$\sss
to\sss the subcategory\sss $\mathcal{C}_{\dff \bm{n}}$\sss
we get\sss a functor\dss
$p\dff \colon\dff
\mathcal{C}_{\dff \bm{n}}
\qff \ttoo\qff
\bm{n}$\nnsp.\oss
This\dss functor\sss induces a simplicial\dss map\dss\vspace{1.5pt}
\[
\quad
\mathit{B}\dff p\dff \colon\dff
\mathit{B}\qff \mathcal{C}_{\dff \bm{n}}
\qff \ttoo\qff 
\mathit{B}\dff \bm{n}
\off =\off
\bm{\Delta}\dff[\dff \infty\dff]
\pff.
\]

\vspace{-10.5pt}
By\dss the definition of\trs the category\sss 
$\mathcal{C}_{\dff \bm{n}}$\sss
a morphism $f$ of\trs this category\sss is\dss an\sss identity\sss morphism\dss
if\trs and\dss only\trs if\dss
$p\trf(\trf f\trf)$\sss is\dss an\dss identity\sss morphism.\oss
It\dss follows\dss that\sss a simplex $\sigma$ of\dss
$\mathit{B}\qff \mathcal{C}_{\dff \bm{n}}$\sss
is\dss non-degenerate\dss
if\trs and\dss only\trs if\dss
$\mathit{B}\dff p\dff(\dff \sigma\dff)$\sss is\dss non-degenerate.\oss
This implies\sss that\vspace{1.5pt}
\[
\quad
\core(\trf \mathit{B}\qff \mathcal{C}_{\dff \bm{n}}\trf)
\off =\off
\mathit{B}\qff \mathcal{C}\dff \times\dff \core\bm{\Delta}\dff[\dff \infty\dff]
\off =\off
\mathit{B}\qff \mathcal{C}\dff \times\dff \Delta\dff[\dff \infty\dff]
\off =\off
\mathcal{B}\fff \mathcal{C},
\]

\vspace{-10.5pt}
where\sss $\mathit{B}\qff \mathcal{C}$\sss is\dss considered as a $\Delta$\dnsp-set.\oss
Also,\oss since $\bm{\Delta}\dff[\dff \infty\dff]$ has\sss non-degenerate core,\oss
this implies\sss that\sss 
$\mathit{B}\qff \mathcal{C}_{\dff \bm{n}}$ 
has non-degenerate core.\oss
Therefore\dss
$\mathit{B}\qff \mathcal{C}_{\dff \bm{n}}
\off =\off
\bm{\Delta} \core(\trf \mathit{B}\qff \mathcal{C}_{\dff \bm{n}}\trf)
\off =\off
\bm{\Delta}\dff \mathcal{B}\fff \mathcal{C}$\nnsp.\oss  \eproof

\newpage
\myappend{Few\qss technical\pss lemmas}{lemmas}

\myapar{Lemma.}{averaging-cochain}
\emph{If\qss the averaging\dss maps\dss $m_{\dff n}$\nsp,\dss $n\qff \in\pff \nnn$\qss 
form\dss a coherent\dss family,\oss
then\sss $m_{\dff *}$\sss is\dss a cochain\dss map.\oss}

\proof
For an $n$\dnsp-cochain\dss 
$f\qff \in\qff B^{\dff n}\dff(\trf K\dff \times\dff \Gamma\trf)$\dss
and\sss an $n$\dnsp-simplex\sss 
$\sigma\qff \in\qff K_{\dff n}$\dss 
let\dss 
$f_{\dff \sigma}\dff \colon\dff
\Gamma_n\qff \ttoo\qff \rrr$\dss
be defined\dss by\dss 
$f_{\dff \sigma}\dff(\dff \tau\trf)
\off =\off
f\dff(\dff \sigma\fff,\qff \tau\trf)$\nnsp.\oss
Then\dss
$m_{\dff *}\dff(\trf f\trf)\dff(\dff \sigma\trf)
\off =\off
m_{\dff n}\dff(\trf f_{\dff \sigma}\trf)$\nnsp.\oss

Suppose\sss that\dss
$f\qff \in\qff B^{\dff n}\dff(\trf K\dff \times\dff \Gamma\trf)$\dss
and\dss
$\rho\qff \in\pff K_{\dff n\dff +\dff 1}$\nsp.\oss 
Then\vspace{6pt}
\[
\quad
\partial^{\dff *} \left(\qff m_{\dff *}\dff\left(\trf f\trf\right)\qff\right)\dff(\trf \rho\trf)
\off =\off
\sum_{i\qff =\qff 0}^{n\dff +\dff 1}\qff
\left(\qff -\qff 1\trf\right)^i\qff
m_{\dff *}\dff(\trf f\trf)\trf
\left(\trf \partial_{\dff i}\dff \rho\trf\right)
\]

\vspace{-30pt}
\[
\quad
\phantom{\partial^{\dff *} \left(\qff m_{\dff *}\dff\left(\trf f\trf\right)\qff\right)\dff(\trf \rho\trf)
\off }
=\off
\sum_{i\qff =\qff 0}^{n\dff +\dff 1}\qff
\left(\qff -\qff 1\trf\right)^i\qff
m_{\dff n}\dff
\left(\trf f_{\pff \partial_{\dff i}\dff \rho}\trf\right) 
\]

\vspace{-30pt}
\[
\quad
\phantom{\partial^{\dff *} \left(\qff m_{\dff *}\dff\left(\trf f\trf\right)\qff\right)\dff(\trf \rho\trf)
\off }
=\off
\sum_{i\qff =\qff 0}^{n\dff +\dff 1}\qff
\left(\qff -\qff 1\trf\right)^i\qff
m_{\dff n\dff +\dff 1}\dff\left(\trf 
\partial_{\dff i}^{\dff *}\dff\left(\trf f_{\pff \partial_{\dff i}\dff \rho}\trf\right) 
\trf\right)
\pff
\]

\vspace{-6pt}
because\dss
$m_{\dff n}
\off =\off
m_{\dff n\dff +\dff 1}\dff \circ\dff \partial_{\dff i}$\nsp.\oss
If\dss $\tau\qff \in\qff \Gamma_{n\dff +\dff 1}$\nsp,\oss
then\vspace{6pt}
\[
\quad
\partial_{\dff i}^{\dff *}\dff\left(\trf f_{\pff \partial_{\dff i}\dff \rho}\trf\right)\dff(\dff \tau\trf)
\off =\off
f\dff\left(\trf \partial_{\dff i}\dff \rho\fff,\pff \partial_{\dff i}\dff \tau\trf\right)
\]

\vspace{-36pt}
\[
\quad
\phantom{\partial_{\dff i}^{\dff *}\dff\left(\trf 
f_{\pff \partial_{\dff i}\dff \rho}\trf\right)\dff
(\dff \tau\trf)
\off }
=\off
f\dff\left(\trf \partial_{\dff i}\dff (\trf \rho\fff,\pff \tau\trf) \trf\right)
\]

\vspace{-36pt}
\[
\quad
\phantom{\partial_{\dff i}^{\dff *}\dff\left(\trf 
f_{\pff \partial_{\dff i}\dff \rho}\trf\right)\dff
(\dff \tau\trf)
\off }
=\off
\partial_{\dff i}^{\dff *}\dff\left(\trf f\trf\right)\dff(\trf \rho\fff,\qff \tau\trf)
\off =\off
\partial_{\dff i}^{\dff *}\dff\left(\trf f\trf\right)_{\dff \rho}\dff(\dff \tau\trf)
\]

\vspace{-6pt}
It\dss follows\dss that\qss
$\partial_{\dff i}^{\dff *}\dff\left(\trf f_{\pff \partial_{\dff i}\dff \rho}\trf\right)
\off =\off
\partial_{\dff i}^{\dff *}\dff\left(\trf f\trf\right)_{\dff \rho}$\qss
and\dss hence\vspace{6pt}
\[
\quad
\partial^{\dff *} \left(\qff m_{\dff *}\dff(\trf f\trf)\qff\right)\dff(\trf \rho\trf)
\off =\off
\sum\nolimits_{\qff i\qff =\qff 0}^{\qff n}\qff
\left(\qff -\qff 1\trf\right)^i\qff
m_{\dff n\dff +\dff 1}\dff\left(\qff 
\partial_{\dff i}^{\dff *}\dff\left(\trf f\trf\right)_{\dff \rho} 
\trf\right)
\pff.
\]

\vspace{-6pt}
Since\sss the maps\sss $m_{\dff n\dff +\dff 1}$\sss
and\dss
$h\off \longmapsto\off h_{\trf \rho}$\dss
are\sss linear\halfff,\oss
it\dss follows\dss that\vspace{3pt}
\[
\quad
\partial^{\dff *} \left(\qff m_{\dff *}\dff(\trf f\trf)\qff\right)\dff(\trf \rho\trf)
\off =\off
m_{\dff n\dff +\dff 1}\dff\left(\pff
\sum_{i\qff =\qff 0}^{n\dff +\dff 1}\qff
\left(\qff -\qff 1\trf\right)^i\qff 
\partial_{\dff i}^{\dff *}\dff\left(\trf f\trf\right)_{\dff \rho} 
\qff\right)
\]

\vspace{-24pt}
\[
\quad
\phantom{\partial^{\dff *} \left(\qff m_{\dff *}\dff(\trf f\trf)\qff\right)\dff(\trf \rho\trf)
\off }
=\off
m_{\dff n\dff +\dff 1}\dff\left(\qff \partial^{\dff *} f\trf\right)\dff(\trf \rho\trf)
\]

\vspace{-3pt}\vspace{-1.75pt}
and\dss hence\dss
$\partial^{\dff *} \left(\qff m_{\dff *}\dff(\trf f\trf)\qff\right)
\off =\off
m_{\dff *}\dff\left(\qff \partial^{\dff *}\dff f\trf\right)$\nnsp.\oss  \eproof

\myapar{Lemma.}{cocycle-homomorphism}
\emph{Suppose\sss that\dss $n\qff \geq\qff 1$\nnsp.\oss
A\dss normalized $n$\dnsp-cochain\dss
$c\dff \colon\dff \pi\qff \ttoo\qff \pi$\dss of\oss
$K\dff(\dff \pi\fff,\qff n\trf)$\dss
is\dss a cocycle\dss if\trs and\dss only\trs if\qss
$c$\dss is\dss a\sss homomorphism\dss
$\pi\qff \ttoo\qff \pi$\nnsp.\oss}

\proof
An $n$\dnsp-cochain $u$ of\trs $\bm{\Delta}^{n\dff +\dff 1}$\sss is\dss
determined\dss by\dss its values\dss
$u_{\dff i}
\off =\off
u\trf(\trf \partial_{\dff i}\qff \bm{\iota}_{\dff n\dff +\dff 1}\trf)$\dss
on\dss the non-de\-gen\-er\-ate $n$\dnsp-simplices of\dss $\bm{\Delta}^{n\dff +\dff 1}$\dnsp.\oss
Therefore,\oss one can\dss identify\sss $u$\sss with\sss the 
$(\dff n\qff +\qff 2\dff)$\dnsp-tuple\vspace{3pt}
\[
\quad
(\trf u_{\trf 0}\dff,\off
u_{\trf 1}
\dff,\off
\ldots\dff,\off
u_{\dff n\dff +\dff 1} 
\trf)
\pff
\]

\vspace{-9pt}
of\dss elements of\dss $\pi$\nnsp.\oss
The boundary\sss operators $\partial_{\dff i}$ are given\sss by\dss the restrictions\sss to\dss
faces,\oss i.e.\vspace{3pt}
\[
\quad
\partial_{\dff i}\trf
(\trf u_{\trf 0}\dff,\off
u_{\trf 1}
\dff,\off
\ldots\dff,\off
u_{\dff n\dff +\dff 1} 
\trf)
\off =\off
u_{\dff i}
\]

\vspace{-9pt}
Suppose\sss that $\pi$ is\dss abelian.\oss
Then $u$\sss is\dss cocycle,\oss
i.e.\qss belongs\sss to\sss $K\dff(\dff \pi\fff,\qff n\trf)_{\dff n\dff +\dff 1}$\nsp,\oss
if\trs and\dss only\trs if\vspace{3pt}
\begin{equation*}
\label{em-cocycle}
\quad
\sum_{i\qff =\qff 0}^{n\dff +\dff 1}\off 
(\qff -\qff 1\dff)^{\dff i}\qff u_{\dff i}
\off =\off
0
\pff.
\end{equation*}

\vspace{-9pt}
An $n$\dnsp-cochain\dss
$c\dff \colon\dff \pi\qff \ttoo\qff \pi$\dss
is\dss a cocycle\dss if\trs and\dss only\trs if\qss
$\partial^{\dff *}\dff c\trf(\dff u\trf)\off =\off 0$\dss
for every\sss simplex\dss
$u\qff \in\pff K\dff(\dff \pi\fff,\qff n\trf)_{\dff n\dff +\dff 1}$\nsp,\oss
i.e.\qss if\trs and\dss only\trs if\trs the last\sss equality\dss implies\vspace{3pt}
\begin{equation*}
\label{em-cocycle-map}
\quad
\sum_{i\qff =\qff 0}^{n\dff +\dff 1}\off 
(\qff -\qff 1\dff)^{\dff i}\qff c\trf(\dff u_{\dff i}\trf)
\off =\off
0
\end{equation*}

\vspace{-9pt}
for every $(\dff n\qff +\qff 2\dff)$\dnsp-tuple $u$\nnsp.\oss
Clearly,\oss this\dss is\dss the case when\sss $c$\sss is\dss a\sss homomorphism.\oss
Conversely,\oss if\dss $c$\dss is\dss a cocycle,\oss
then\dss the\sss last\sss equality\dss for\sss the
$(\dff n\qff +\qff 2\dff)$\dnsp-tuples\vspace{3pt}
\[
\quad
(\trf u_{\trf 0}\dff,\off
u_{\trf 1}
\dff,\off
\ldots\dff,\off
u_{\dff n\dff +\dff 1} 
\trf)
\off =\off
(\trf v\dff,\off
v\qff +\qff w\dff,\off
w\dff,\off
0\dff,\off
\ldots\dff,\off
0 
\trf)
\pff,
\]

\vspace{-9pt}
where\dss $v\fff,\pff w\qff \in\qff \pi$\nnsp,\oss
together\dss with\dss the fact\sss that\dss
$c\trf(\dff 0\dff)\off =\off 0$\dss
implies\sss that\dss\vspace{3pt}
\[
\quad
c\trf(\dff v\trf)
\qff -\qff 
c\trf(\dff v\qff +\qff w\trf)
\qff +\qff
c\trf(\dff w\trf)
\off =\off
0
\]

\vspace{-9pt}
for every\dss $v\fff,\pff w\qff \in\qff \pi$\nnsp.\oss
This proves\sss that\sss $c$\sss is\dss a homomorphism\sss 
when\sss $\pi$\sss is\dss abelian.\oss

If\dss $\pi$\sss is\dss not\sss abelian,\oss
then\dss $n\off =\off 1$\dss and a\sss triple\sss
$u\off =\off 
(\trf u_{\trf 0}\dff,\pff u_{\trf 1}\dff,\pff u_{\trf 2} \trf)$\sss
is\dss a cocycle\dss if\trs and\dss only\trs if\dss
$u_{\trf 1}\off =\off u_{\trf 2}\dff \cdot\dff u_{\trf 0}$\nsp.\oss
It\dss follows\dss that\sss $c$\sss is\dss a cocycle\sss
if\trs and\dss only\trs if\dss 
$c\trf(\dff u_{\trf 2}\dff \cdot\dff u_{\trf 0}\dff)
\off =\off
c\trf(\dff u_{\trf 2}\dff)\dff \cdot\dff c\trf(\dff u_{\trf 0}\dff)$
for every\dss pair\sss
$u_{\trf 2}\dff,\pff u_{\trf 0}\qff \in\qff \pi$\nnsp,\oss
i.e.\qss if\trs and\dss only\trs if\dss $c$\sss
is\dss a homomorphism.\oss
In addition,\oss we see\sss that\dss when\dss $n\off =\off 1$\nnsp,\oss
every\sss cocycle\dss is\dss automatically\dss normalized.\oss   \eproof

\myapar{Lemma.}{invariant-cochains}
\emph{Suppose\sss that\trs $K$\dss is\dss a simplicial\sss set\sss
and\dss $G$\dss is\dss an amenable\sss group acting\sss on\dss $K$\dss
on\dss the\dss left\dss
by\dss automorphisms\sss homotopic\sss to\sss the identity.\oss
Then every\dss bounded cocycle\sss $c$\sss of\pss $K$\dss is\dss boundedly\dss
cohomologous\sss to a\dss $G$\dnsp-invariant\dss bounded cocycle
with\dss the norm\qss $\leq\qff \norm{c}$\nnsp.\oss}

\proof
We will\sss denote\sss the action\sss by\dss
$(\trf g\fff,\qff \sigma\trf)\off \longmapsto\off g\dff \cdot\dff \sigma$\nnsp,\oss
where\sss $g\qff \in\qff G$\sss and\sss $\sigma$\sss is\dss a simplex of\trs $K$\nnsp.\oss
Let\dss $B\dff(\trf G\trf)$\sss be\sss the space of\dss bounded\sss real-valued\dss
functions on\sss $G$\nnsp.\oss
For\dss $g\qff \in\qff G$\dss and\dss $f\qff \in\pff B\dff(\trf G\trf)$\dss
let\dss $g\dff \cdot\dff f$\dss be\sss the function\dss
$h\off \longmapsto\off f\dff(\trf h g\trf)$\nnsp.\oss
This defines an action of\dss $G$\dss on\dss $B\dff(\trf G\trf)$\nnsp.\oss

Since\sss $G$\sss is\dss amenable,\oss there exists a\trs
\emph{$G$\dnsp-invariant\dss mean}\qss on\sss $B\dff(\trf G\trf)$\nnsp,\oss
i.e.\qss a\sss linear\dss functional\dss
$\mu\dff \colon\dff
B\dff(\trf G\trf)\qff \ttoo\qff \rrr$\dss
such\dss that\dss the norm of\dss $\mu$\sss is\dss $\leq\qff 1$\nnsp,\pss
$\mu$\sss takes a constant\sss function\dss to its value,\oss
and\dss $\mu\dff(\trf g\dff \cdot\dff f\trf)\off =\off \mu\dff(\trf f\trf)$\dss
for every\dss $g\qff \in\qff G$\sss and\sss $f\qff \in\pff B\dff(\trf G\trf)$\nnsp.\oss

For\dss $g\qff \in\qff G$\dss let\dss
$a\trf(\dff g\trf)\dff \colon\dff 
K\qff \ttoo\qff K$\dss
be\sss the automorphism defined\dss by\sss $g$\nnsp.\oss
Since\sss $a\trf(\dff g\trf)$\sss is\dss homotopic\sss to\sss the identity,\oss
there exists a cochain\dss homotopy\dss
between\dss 
$a\trf(\dff g\trf)^{\dff *}\dff \colon\dff
B^{\fff *}\dff(\trf K\trf)
\qff \ttoo\qff
B^{\fff *}\dff(\trf K\trf)$\dss
and\dss the identity.\oss
In other\sss words,\oss 
for each\dss $m\qff >\qff 0$\dss a\sss homomorphism\vspace{3pt}
\[
\quad
k_{\dff m}\dff(\dff g\trf)\dff \colon\dff
B^{\fff m}\dff(\trf K\trf)
\qff \ttoo\qff
B^{\fff m\dff -\dff 1}\dff(\trf K\trf)
\]

\vspace{-10.5pt}
is\dss defined,\oss and\vspace{1.5pt}
\[
\quad
a\trf(\dff g\trf)^{\dff *}\dff(\dff c\trf)\qff -\pff c
\off =\off
k_{\dff m\dff +\dff 1}\dff(\dff g\trf)\fff \circ\dff \partial^{\dff *}\dff(\dff c\trf)
\off +\off
\partial^{\dff *} \circ\dff k_{\dff m}\dff(\dff g\trf)\dff(\dff c\trf)
\]

\vspace{-9pt}
for every\dss $c\qff \in\pff B^{\fff m}\dff(\trf K\trf)$\nnsp.\oss
Suppose\sss that\dss $c$\sss is\dss a cocycle.\oss
Then\dss this identity\sss simplifies\sss to\vspace{3pt}
\[
\quad
a\trf(\dff g\trf)^{\dff *}\dff(\dff c\trf)\qff -\pff c
\off =\off
\partial^{\dff *} \circ\dff k_{\dff m}\dff(\dff g\trf)\dff(\dff c\trf)
\pff.
\]

\vspace{-9pt}
By\sss applying\dss this equality\dss to an $m$\dnsp-simplex $\sigma$ of\dss $K$\sss
using\dss the definition of\dss $\partial^{\dff *}$\dnsp,\oss
we get\vspace{3pt}
\begin{equation}
\label{g-and-homotopy}
\quad
a\trf(\dff g\trf)^{\dff *}\dff(\dff c\trf)\trf(\dff \sigma\trf)
\qff -\pff 
c\trf(\dff \sigma\trf)
\off =\off
k_{\dff m}\dff(\dff g\trf)\dff(\dff c\trf)\trf(\trf \partial\dff \sigma\trf)
\pff.
\end{equation}

\vspace{-9pt}
We would\dss like\sss to consider all\dss terms of\trs this equality\sss
as functions of\dss $g$\sss and apply\sss $\mu$\sss to\sss them.\vspace{2.875pt}

To begin\dss with,\oss let\sss $\gamma\trf(\dff \sigma\trf)$\sss
be\sss the result\sss of\dss applying\sss $\mu$\sss to\sss the function\vspace{3pt}
\[
\quad
g
\off \longmapsto\off
a\trf(\dff g\trf)^{\dff *}\dff(\dff c\trf)\trf(\dff \sigma\trf)
\off =\off
c\trf(\trf g\dff \cdot\dff \sigma\trf)
\pff.
\]

\vspace{-9pt}
The map\dss
$\gamma\dff \colon\dff
\sigma\off \longmapsto\off \gamma\trf(\dff \sigma\trf)$\sss
is\dss a\sss bounded $m$\dnsp-cochain of\dss $K$\sss
and\dss $\norm{\gamma}\qff \leq\qff \norm{c}$\nnsp.\oss
Since $\mu$\sss is\sss $G$\dnsp-invariant,\pss
$\gamma$\sss is\dss also $G$\dnsp-invariant.\oss
Next,\oss 
the result\sss of\dss applying\sss $\mu$\sss to\sss the constant\dss map\dss
$g\off \longmapsto\off c\trf(\dff \sigma\trf)$\dss
is\sss $c\trf(\dff \sigma\trf)$\nnsp.\oss 
Let\dss $\tau$\sss be an $(\dff m\qff -\qff 1\dff)$\dnsp-simplex 
of\dss $K$\sss and\sss consider\dss the function\vspace{3pt}
\[
\quad
g
\off \longmapsto\off
k_{\dff m}\dff(\dff g\trf)\dff(\dff c\trf)\trf(\dff \tau\trf)
\pff.
\]

\vspace{-9pt}
Let\dss 
$\kappa_{\dff m}\dff(\dff c\trf)\trf(\dff \tau\trf)
\qff \in\qff \rrr$\dss be\sss the result\sss of\dss applying\sss $\mu$\sss to\sss
this function.\oss
The map\vspace{3pt}
\[
\quad
\kappa_{\dff m}\dff(\dff c\trf)\dff \colon\dff
\tau
\off \longmapsto\off
\kappa_{\dff m}\dff(\dff c\trf)\trf(\dff \tau\trf)
\]

\vspace{-9pt}
is\dss a bounded $(\dff m\qff -\qff 1\dff)$\dnsp-cochain of\dss $K$\nnsp.\oss
i.e.\qss 
$\kappa_{\dff m}\dff(\dff c\trf)
\qff \in\pff 
B^{\fff m\dff -\dff 1}\dff(\trf K\trf)$\nnsp.\oss
In\dss terms of\dss $\gamma$\dss and\dss $\kappa_{\dff m}\dff(\dff c\trf)$\dss
the result\sss of\dss applying\sss $\mu$\sss 
to\qss (\ref{g-and-homotopy})\qss can\sss be written as follows\fff:\vspace{3pt}
\[
\quad
\gamma\trf(\dff \sigma\trf)
\qff -\pff 
c\trf(\dff \sigma\trf)
\off =\off
\kappa_{\dff m}\dff(\dff c\trf)\trf(\trf \partial\dff \sigma\trf)
\pff.
\]

\vspace{-9pt}
Therefore\dss
$\gamma\qff -\qff c
\off =\off
\partial^{\dff *}\fff \kappa_{\dff m}\dff(\dff c\trf)$\nnsp.\oss
The\sss lemma\sss follows.\oss  \eproof

\myapar{Lemma.}{eilenberg-zilber}
\emph{Every\sss $n$\dnsp-simplex $\sigma$ of\dss a simplicial\sss set\dss $K$\sss
admits a unique presentation of\trs the form\dss
$\sigma\off =\off \theta^{\fff *}\dff(\dff \tau\trf)$\sss
with a surjective non-decreasing\dss map $\theta$ and a non-degenerate simplex\sss $\tau$\nnsp.\oss}

\proof
This\dss is\dss a\sss well\dss known\dss lemma of\pss 
Eilenberg\dss and\dss Zilber\qss \cite{ez}.\oss
See\qss \cite{ez},\oss (8.3).\oss

Let\dss us\dss choose among all\dss presentations\dss
$\sigma\off =\off \theta^{\fff *}\dff(\dff \tau\trf)$\dss
with surjective\dss
$\theta\dff \colon\dff [\halfff n\dff]\qff \ttoo\qff [\halfff m\dff]$\sss
and\sss an $m$\dnsp-simplex $\tau$\sss some presentation\sss
with\dss minimal\dss possible $m$\nnsp.\oss
Clearly,\oss the minimality\sss of\dss $m$\sss implies\sss that\sss
$\tau$\sss is\dss non-degenerate.\oss
This proves\sss the existence.\oss
Suppose\sss that\sss also\dss
$\sigma
\off =\off
\eta^{\fff *}\dff(\trf \rho\trf)$\nnsp,\oss
where\dss 
$\eta\dff \colon\dff [\halfff n\dff]\qff \ttoo\qff [\fff k\dff]$\sss
is\dss surjective\sss and $\rho$ is\dss a $k$\dnsp-simplex.\oss
Since\sss $\theta\fff,\pff \eta$\sss are surjective non-decreasing maps,\oss
there exist\sss strictly\sss increasing\sss maps\dss
$\alpha\dff \colon\dff [\halfff m\dff]\qff \ttoo\qff [\halfff n\dff]$\sss
and\dss
$\beta\dff \colon\dff [\fff k\dff]\qff \ttoo\qff [\halfff n\dff]$\sss
such\dss that\dss $\theta\dff \circ\dff \alpha$\sss and\dss
$\eta\dff \circ\dff \beta$\sss are\sss the identity\sss maps.\oss
Then\vspace{3pt}
\[
\quad
(\trf \eta\dff \circ\dff \alpha\trf)^{\fff *}\dff(\trf \rho\trf)
\off =\off
\alpha^{\dff *}\dff \left(\trf \eta^{\dff *}\dff (\trf \rho\trf)\trf\right)
\off =\off
\alpha^{\dff *}\dff \left(\trf \theta^{\fff *}\dff(\dff \tau\trf)\trf\right)
\off =\off
(\trf \theta\dff \circ\dff \alpha\trf)^{\fff *}\dff(\trf \tau\trf)
\off =\off
\tau
\pff.
\]

\vspace{-9pt}
Similarly,\pss
$(\trf \theta\dff \circ\dff \beta\trf)^{\fff *}\dff(\trf \tau\trf)
\off =\off
\rho$\nnsp.\oss
Since\sss $\tau$\sss and\sss $\rho$\sss are both\sss non-degenerate,\oss
both\sss $\eta\dff \circ\dff \alpha$\sss and\sss 
$\theta\dff \circ\dff \beta$\sss are strictly\dss injective.\oss
It\dss follows\sss that\dss $m\off =\off k$\dss
and\dss both\sss $\eta\dff \circ\dff \alpha$\sss and\sss 
$\theta\dff \circ\dff \beta$\sss are equal\dss to\sss the identity.\oss
In\dss turn,\oss this\sss implies\sss that\dss $\tau\off =\off \rho$\nnsp.\oss
Suppose\sss that\dss
$\theta\dff(\dff i\trf)\off \neq\off \eta\dff(\dff i\trf)$\dss
for some\dss $i\qff \in\qff [\halfff n\dff]$\nnsp.\oss
One can\dss choose\sss the map\sss $\alpha$\sss in such a way\dss that\dss
$\alpha\dff (\trf \theta\dff(\dff i\trf)\trf)\off =\off i$\nnsp.\oss
Then\vspace{3pt}
\[
\quad
(\trf \eta\dff \circ\dff \alpha\trf)\trf \bigl(\trf \theta\dff(\dff i\trf)\trf\bigr)
\off =\off
\eta\dff \circ\dff \alpha\dff \circ\dff \theta\dff(\dff i\trf)
\off =\off
\eta\dff (\dff i\trf)
\off \neq\off
\theta\dff (\dff i\trf)
\pff,
\]

\vspace{-9pt}
contrary\dss to\sss $\eta\dff \circ\dff \alpha$\sss
being\dss equal\dss to\sss the identity.\oss
Hence\dss $\theta\off =\off \eta$\nnsp.\oss
The uniqueness follows.\oss  \eproof

\begin{flushright}

October\qss 10,\oss {2020}
 
https\halfff:/\!/\hspace*{-0.06em}nikolaivivanov.com

E-mail\halfff:\oss nikolai.v.ivanov{\fff}@{\dff}icloud.com

\end{flushright}

\end{document}